\documentclass[letterpaper]{article}
\usepackage{authblk}
\usepackage{url}
\usepackage[margin=1in,footskip=0.25in]{geometry}
\usepackage[bbgreekl]{mathbbol}
\usepackage{amssymb}             % AMS Math
\DeclareSymbolFontAlphabet{\amsmathbb}{AMSb}%
\usepackage{graphicx}
\usepackage{amsmath}
\usepackage{algorithmic} 
\usepackage{algorithm}
\usepackage{float}
\usepackage{lipsum}
\usepackage{subfigure}
\usepackage{wrapfig}
\usepackage{color}
\usepackage{easylist}
\usepackage{lineno}
\usepackage{stmaryrd}
\usepackage{booktabs}
\usepackage{cleveref}
\usepackage{multirow}
\usepackage{hhline}
\usepackage{filecontents}
\usepackage{enumitem}
\usepackage{cite}
\newlist{steps}{enumerate}{1}
\setlist[steps, 1]{label = Step \arabic*:}
\def \grad{\nabla}
\def \half{\frac{1}{2}}
\def \p   {\partial}

\newcommand{\volume}{{\ooalign{\hfil$V$\hfil\cr\kern0.08em--\hfil\cr}}}

\usepackage{blindtext}

\def \II{\mathbb{I}}

\def \QQ{\mathbb{Q}}
\def \OO{\mathbb{\Omega}}

\def \sigmaf{\bbsigma_\text{f}}

\def \F{\vec{F}}
\def \disp{\vec{d}}
\def \dispdot{\dot{\vec{d}}}

\def \T{\vec{T}}

\def \U{\vec{U}}
\def \JJ{\mathbb{J}}

\def \X{\vec{\chi}}
\def \Y{\vec{\xi}}

\def \f{\vec{f}}
\def \g{\vec{g}}
\def \n{\vec{n}}
\def \s{\vec{X}}

\def \f{\vec{f}}
\def \u{\vec{u}}

\def \x{\vec{x}}

\def \grad{\nabla}
\def \half{\frac{1}{2}}
\def \p{\partial}
\def \fs{\textrm{fs}}
\def \com{\textrm{c}}
\def \Omegas{\Omega^\textrm{s}}
\def \Omegaf{\Omega^\textrm{f}}

\def \DA{{\mathrm d} A}

\def \rhos{\rho_{\textrm{s}}}
\def \Ks{K_{\textrm{s}}}
\def \Cs{C_\textrm{s}}
\def \Is{I_\textrm{s}}
\def \Ms{M_\textrm{s}}
\def \rhof{\rho_{\textrm{f}}}
\def \tauf{\vec{\tau}_{\textrm{f}}}
\def \muf{\mu_{\textrm{f}}}

\def \cL{\mathcal L}
\def \cE{\mathcal E}
\def \cI{\mathcal I}
\def \cJ{\vec{\mathcal J}}

\def \cS{\vec{\mathcal S}}
\def \cT{\mathcal T}

\def \dt{\Delta t}

\def \cI{\mathcal I}

\def \cS{\mathcal S}
\def \cT{\mathcal T}

\def \mfac{M_\text{fac}}

\def \Re{\textrm{Re}}
\def \Ured{U^{*}}
\def \Mred{m^{*}}
\def \St{\textrm{St}}

\def \L2{L^2}
\def \Linf{L^{\infty}}

\def \grad{\nabla}
\def \half{\frac{1}{2}}
\def \p   {\partial}
\def \Mfac {M_{\text{fac}}}

\renewcommand{\vec}[1]{\ensuremath\boldsymbol{#1}}
\newcommand{\norm}[1]{\left|\left|#1\right|\right|}

\title{A sharp interface Lagrangian-Eulerian method for rigid-body fluid-structure interaction}

\author[1,2]{Ebrahim M.~Kolahdouz} %\corref{cor1}}
\author[3]{Amneet P.S. Bhalla}
\author[4]{Lawrence N. Scotten}
\author[2]{Brent A.~Craven}
\author[5,6,7]{Boyce E.~Griffith} %\corref{cor2}}
\affil[1]{Department of Mathematics, University of North Carolina, Chapel Hill, NC, USA}
\affil[2]{Division of Applied Mechanics, Office of Science and Engineering Laboratories, Center for Devices and Radiological Health, United States Food and Drug Administration, Silver Spring, MD, USA}
\affil[3]{Department of Mechanical Engineering, San Diego State University, San Diego, CA, USA}
\affil[4]{Independent Consultant, Victoria, BC, Canada}
\affil[5]{Departments of Mathematics, Applied Physical Sciences, and Biomedical Engineering, University of North Carolina, Chapel Hill, NC, USA}
\affil[6]{Carolina Center for Interdisciplinary Applied Mathematics, University of North Carolina, Chapel Hill, NC, USA}
\affil[7]{McAllister Heart Institute, University of North Carolina, Chapel Hill, NC, USA}
\affil[ ]{\texttt{ebrahimk@email.unc.edu} and \texttt{boyceg@email.unc.edu}}
\begin{document}

\maketitle

\begin{abstract}
This paper introduces a sharp interface method to simulate fluid-structure interaction (FSI)
involving rigid bodies immersed in viscous incompressible fluids. The capabilities of this methodology 
are benchmarked using a range of test cases and demonstrated using large-scale models of biomedical FSI.
The numerical approach developed herein, which we refer to as an immersed Lagrangian-Eulerian (ILE) method, integrates aspects of
partitioned and immersed FSI formulations by solving separate momentum equations for the fluid and solid subdomains, 
as in a partitioned formulation, while also 
using non-conforming discretizations of the dynamic fluid and structure regions, as in an immersed formulation.
A simple Dirichlet-Neumann coupling scheme is used, in which the motion of the immersed solid is driven by fluid traction forces evaluated along 
the fluid-structure interface, and the motion of the fluid along that interface is constrained to match the solid velocity and thereby satisfy the no-slip condition.
 To develop a practical numerical method, we adopt a penalty approach that approximately imposes the no-slip condition along the fluid-structure interface.
 In the coupling strategy, a separate discretization of the fluid-structure interface is tethered to the volumetric 
 solid mesh via stiff spring-like penalty forces. 
 Our fluid-structure coupling scheme relies on an immersed interface method (IIM) for discrete geometries, which enables the accurate determination of both velocities and stresses 
 along complex internal interfaces.
Numerical methods for FSI can suffer from instabilities related to the added mass effect, 
but the computational tests indicate that the methodology introduced here remains stable for selected test cases 
across a range of solid-fluid density ratios, including extremely small, nearly equal, equal, and large density ratios.
Biomedical FSI demonstration cases include results obtained using this method
to simulate the dynamics of a bileaflet mechanical heart valve in a pulse duplicator, and to model transport of
blood clots in a patient-averaged anatomical model of the inferior vena cava.

\end{abstract}


%%%%%%%%%%%%%%%%%%%%%%%%%%%%%%%%%%%%%%%%%%%%%%%%%%%%%%%%%%%%%%%%%%%%%%%%%%%%%%%%%%%%%%%%%%%%%%%%%%%%%%%%%%%%%%%%%%%%%%%%%%%
%%%%%%%%%%%%%%%%%%%%%%%%%%%%%%%%%%%%%%%%%%%%%%%%%%%%%%%%%%%%%%%%%%%%%%%%%%%%%%%%%%%%%%%%%%%%%%%%%%%%%%%%%%%%%%%%%%%%%%%%%%%
\section{Introduction}
\label{sec:intro}

Predictive numerical models of fluid-structure interaction (FSI) have long been of major interest in the scientific computing community. 
Numerical simulations of FSI problems can be characterized by the solution approach taken for the coupled system of momentum equations associated with the fluid and structure.
\textit{Partitioned formulations} of FSI 
describe a fluid-structure system using distinct, non-overlapping fluid and structure regions.
Commonly used numerical methods of this type include well-known arbitrary Lagrangian-Eulerian (ALE) schemes \cite{takashi1992arbitrary,hu2001direct,ahn2006strongly,lee2017fluid}.
These formulations can yield outstanding resolution of flows and stresses up to the fluid-structure interface. 
Despite their high accuracy, however, ALE methods for FSI can be constrained by the difficulties of body conforming grid regeneration and mesh morphing, 
which can make it challenging to use these approaches to simulate systems involving very large structural displacements or deformations, and 
to handle transient contact between moving structures.
Overset methods using overlapping Chimera grid systems are another class of partitioned methods that have been used
to simulate moving rigid and flexible bodies \cite{zahle2009wind, nakata2012fluid, koblitz2017direct,banks2017stable,banks2012deforming}.
In these approaches, a complex domain is decomposed into multiple geometrically simple overlapping grids, and boundary information is exchanged between
these grids through interpolation.
\textit{Immersed formulations} of FSI \cite{peskin2002immersed, mittal2005immersed,hou2012numerical,griffith2020immersed} 
are alternatives to body-fitted methods. 
Many immersed approaches to FSI have been developed, including Peskin's immersed boundary (IB) method \cite{peskin2002immersed} and various sharp-interface 
IB methods \cite{fadlun2000combined,udaykumar2001sharp, gilmanov2005hybrid, xu2006, mittal2008versatile,borazjani2008curvilinear,barad2009adaptive}. These
methods commonly describe the fluid in Eulerian form (i.e.~using fixed physical coordinates) 
and the structure in Lagrangian form (i.e.~using material coordinates attached to the structure), 
and they use non-conforming
discretizations along the fluid-structure interface.  
Because these methods avoid using body-conforming discretizations of the interface, they are readily able to treat models with very large structural deformations or displacements, 
and they facilitate simulations that fundamentally involve contact or near-contact between structures \cite{fai2018lubricated, griffith2020immersed}.
The method presented herein, which we call an immersed Lagrangian-Eulerian (ILE) method, combines a partitioned approach to FSI
with an immersed coupling strategy.

The key challenge in developing immersed methods for FSI is linking the Eulerian and Lagrangian variables. 
Peskin's IB methods, for example, regularize singular forces and stress discontinuities along the 
fluid-structure interface, which enables straightforward discretization approaches but can yield low accuracy 
in the flows and stresses near those interfaces. Efforts have been made to improve the accuracy of the method, including 
the development of formally second-order IB methods
that realize second-order accuracy when applied to specific problems \cite{lai2000,griffith2005order}
and IB methods that use Cartesian grid adaptive mesh refinement to enhance spatial resolution near fluid-structure interfaces \cite{roma1999adaptive,griffith2007adaptive}.
For general FSI problems, however, formally second-order accurate IB methods 
 still only realize first-order convergence rates \cite{griffith2005order,griffith2007adaptive}. 
Motivated by the goal of improving the accuracy of the original IB method, 
the immersed interface method (IIM) was introduced by LeVeque and Li \cite{leveque1994immersed} for  
elliptic PDEs with discontinuous coefficients and singular forces. 
The IIM subsequently was extended to the 
incompressible Stokes \cite{leveque1997immersed,li2007augmented} and Navier-Stokes 
\cite{li2001immersed, lee2003immersed, le2006immersed} equations, 
and it was combined with level set methods to represent 
the interface \cite{hou1997hybrid, sethian2000structural, xu2006level}. 
When applied to the incompressible Navier-Stokes equations, the IIM sharply imposes interfacial stress discontinuities 
through an extended finite difference discretization that accounts for jump conditions induced 
by singular forces at the interface.
Modern IIMs use generalized Taylor series expansions to extend the physical jump conditions from the interface to the finite difference
discretization in the Eulerian domain while permitting the use of efficient linear solvers based on the unmodified finite difference
discretizations \cite{xu2006systematic,le2006immersed}.
The IIM has been used to simulate various phenomena, including acoustics and elastodynamic wave propagation \cite{lombard2004numerical},
fluid interfaces with insoluble surfactants \cite{xu2006level}, the osmotic swelling of a deforming 
capsule \cite{jayathilake2010deformation}, and vesicle electrohydrodynamics \cite{kolahdouz2015electrohydrodynamics, kolahdouz2015numerical, hu2016vesicle}.
Other sharp interface FSI methods have been developed 
\cite{fadlun2000combined,udaykumar2001sharp, gilmanov2005hybrid, mittal2008versatile,borazjani2008curvilinear,barad2009adaptive},
and most of these methods achieve higher-order accuracy by adopting approaches that are similar
 to body-fitted discretization methods, such as local modifications to the finite difference stencils,
 to allow for the accurate reconstruction of boundary conditions in the vicinity of the immersed interface.

Many different IB methods have been developed to treat FSI problems involving rigid bodies. 
Differences between these approaches are 
mainly related to the way that the rigidity constraint is enforced to account for the effect of the structure in the fluid region.
Previous IB approaches to rigid-body FSI 
that use regularized coupling operators include Lagrange-multiplier-based fictitious-domain methods \cite{patankar2000new,glowinski1999distributed, patankar2005physical, bhalla2013unified,glowinski2001fictitious}, 
direct forcing IB methods \cite{uhlmann2005immersed,zhang2007improved, VANELLA20096617},
projection-based methods \cite{taira2007immersed,lacis2016stable,li2016efficient},
immersed finite element methods \cite{zhang2004immersed,wang2013modified},
methods based on computing exact Lagrange multipliers for the rigidity constraint \cite{kallemov2016immersed,balboa2017hydrodynamics},
penalty immersed boundary methods \cite{kim2016penalty},
immersed boundary lattice Boltzmann methods \cite{feng2004immersed,suzuki2011effect,jianhuaefficient,qin2021immersed},
and level set based approaches \cite{coquerelle2008vortex,gibou2012efficient}.
Sharp-interface approaches designed for rigid body FSI include embedded boundary methods \cite{yang2006embedded,kim2006immersed}, cut-cell methods \cite{mittal2005immersed,schneiders2016efficient,pogorelov2018adaptive}, and
the curvilinear immersed boundary method \cite{gilmanov2005hybrid,borazjani2008curvilinear}.
Unlike the IB approach introduced by Peskin, these methods all solve the fluid equations in the exterior regions
surrounding the immersed object.
Of the sharp-interface IB methods considered here, immersed interface methods are the most similar to Peskin's IB 
method because they treat thin interfaces that are fully immersed in fluid, which enables the use of fast Cartesian 
grid fluid solvers. IIMs have been developed for bodies with prescribed motion 
\cite{xu2006, xu2008immersed, le2006immersed,kolahdouz2020immersed},
but they are more commonly used for thin flexible interfaces \cite{le2006immersed, tan2009immersed, tan2010level, thekkethil2019level,layton2009using}.
To our knowledge, the few IIM models that treat volumetric rigid structures 
do not simulate fluid-structure interaction per se, but instead prescribe the motion of the immersed body \cite{le2006immersed,xu20083d,xu2008immersed,tan2009rigid,liuimmersed2020}.
Xu and Wang used a feedback control to construct the singular force density \cite{xu2006,xu20083d}.
Recent work by some of us introduced an IIM for discrete surfaces described by a general finite element mesh 
to sharply resolve fluid dynamics for problems with prescribed motion \cite{kolahdouz2020immersed}.
The present study uses this method to develop a new sharp interface approach to rigid-body FSI.

The present ILE method introduces a partitioned approach to FSI with an immersed coupling strategy that sharply resolves flow features up to the fluid-structure interface. 
Like partitioned formulations, the present approach uses distinct momentum equations for the fluid and solid regions. However, 
like immersed methods, and unlike typical partitioned methods, our ILE approach uses a non-conforming discretization of the dynamic fluid-structure interface that
does not require any grid regeneration or mesh morphing to treat large structural motions.
The fluid and solid subproblems are solved in a partitioned manner using independent, non-conforming discretizations and are coupled only through interface conditions. 
The ILE equations of motion are first introduced using an exactly constrained formulation that exactly imposes kinematic interface conditions through 
a Lagrange multiplier force distribution applied along the fluid-structure interface.
Solving the exactly constrained equations would require 
the solution of an extended saddle-point system involving an exact Lagrange multiplier force along with
the Eulerian velocity and pressure fields.
Developing efficient linear solvers for such systems is challenging even for conventional 
IB formulations with regularized delta functions \cite{kallemov2016immersed,balboa2017hydrodynamics}. Consequently, to obtain a practical numerical method, 
we next reformulate this scheme using a penalty approach that relaxes the kinematic constraint, and we use this penalty ILE method in all of our numerical tests. Specifically, our penalty formulation uses two representations of the fluid-structure interface, including 
a surface mesh and the boundary of a bulk volumetric structural mesh, that are connected by forces that impose kinematic and dynamic interface conditions. 
The dynamics of the 
volumetric structural mesh are driven by the exterior fluid traction obtained from solving the equations of fluid dynamics.
The surface mesh moves according to the local fluid velocity and locally exerts an approximate Lagrange multiplier force distribution back 
to the fluid generated from stiff spring-like penalty forces that link the surface mesh to the boundary of its volumetric counterpart.
At least formally, in the limit of infinite spring stiffness, the two interface representations become exactly conformal in their motion.
Results demonstrate that this approach is able to control these discrepancies effectively, 
and for sufficiently large penalty spring stiffnesses, the penalty formulation has little impact on the computed dynamics.
To discretize the jump conditions, we leverage our recently developed IIM
for discrete surfaces \cite{kolahdouz2020immersed}, which allows us to impose stress jump conditions along 
complex interfaces within a Cartesian grid framework
and to use fast structured-grid solvers for the incompressible Navier-Stokes equations.
This IIM formulation describes fluid dynamics on \textit{both sides} of the fluid-structure interface. However, only the fluid forces exerted by 
the ``exterior'' fluid have a physical effect on the structural dynamics, and the motion of the structure determines the fluid velocity at the 
fluid-structure interface. 

Our approach also can be viewed as a sharp implementation of the distributed Lagrange multiplier (DLM) technique first introduced by Glowinski, Patankar, and coworkers for immersed rigid structures \cite{patankar2000new,glowinski1999distributed}. 
In the DLM approach to FSI, a Lagrange multiplier force field is introduced to 
impose the kinematic condition at the fluid-structure interface. In the DLM literature, the Lagrange multiplier field (either exact or approximate) 
has typically been smoothed using, e.g., via regularized delta functions like those used in Peskin's IB method. 
From this standpoint, our approach is different from other sharp interface immersed methods, 
in which the velocity matching condition is imposed directly, 
for example, through velocity reconstruction \cite{gilmanov2005hybrid} or through cut-cell approaches \cite{udaykumar2001sharp}. 
These approaches forgo Lagrange multipliers entirely and instead solve the fluid momentum equations only within the fluid subdomain. 
In contrast, immersed approaches, including DLM methods, typically extend the fluid domain into the solid domain, so that the fluid momentum equation 
is solved on the entire computational domain. 
The fluid velocity 
determined within the solid region can be different from the actual solid velocity, although volumetric DLM formulations 
can impose the constraint that the fictitious fluid velocity matches
the solid velocity within the overlapping region occupied by the structure.
Ultimately, however, it is necessary only for the fluid and solid velocities to match along the fluid-structure interface. 
To our knowledge, the present method is the first DLM-type formulation to sharply impose this constraint.
In addition, the current scheme uses a constant-coefficient flow solver and yet can readily treat immersed bodies that are 
lighter than the ambient fluid because it only imposes the kinematic condition along the fluid-solid interface. This is in contrast 
to volumetric DLM schemes, such as the approach of Nangia et al.~\cite{nangia2019dlm}, which approximately constrains 
fluid and solid velocities to match in the extended fluid domain and allows for solid-fluid mass density ratios less than one only 
through the use of a variable-density fluid solver~\cite{nangia2019robust}. 
 %~ At least for the numerical tests reported in this study, the mismatch of the solid and fictitious interior fluid 
 %~ velocity fields does not appear to cause any numerical instabilities, even for high density contrasts between fluid and solid.

Instabilities due to \textit{artificial added mass effect} have been 
observed in weakly coupled FSI schemes in which the fluid and solid equations are linked via explicit time stepping schemes.
Such instabilities can occur if the mass density of the solid $\rhos$ is comparable to or less than the mass density of the fluid $\rhof$.
Added mass effect instabilities have been discussed in both
sharp-interface IB-type methods for FSI ~\cite{borazjani2008curvilinear,zheng2010coupled,yang2014sharp} 
as well as body-fitted methods, including ALE methods for FSI \cite{nobile1999stability, matthies2003partitioned, lee2017fluid,dunbar2015development}.
Various approaches, including subcycling or using modified coupling conditions, have been developed to
maintain stability \cite{borazjani2008curvilinear,nobile1999stability, guidoboni2009stable,burman2009stabilization}.
Strong coupling schemes, in which the governing equations for the fluid and solid subdomains are simultaneously 
integrated in time \cite{fernandez2005newton,matthies2006algorithms,degroote2009performance}, have been shown to improve the stability of these FSI formulations.
In one common strong coupling approach, solutions are transferred between the fluid and structure 
multiple times within each time step (i.e.~through subiterations) until convergence is achieved in the forces and displacements 
\cite{hu2001direct, uhlmann2005immersed, yang2008strongly,dunbar2015development}.
Using subiterations substantially increases the computational cost per time step, however. Further, 
instabilities at low density ratios have still been reported in some situations even if using strong coupling \cite{deparis2003acceleration,kuttler2008fixed,borazjani2008curvilinear}.
Substantial work has been devoted to understanding the sources of these instabilities and to developing methods to overcome these instabilities 
\cite{nobile1999stability,causin2005added,forster2007artificial,burman2009stabilization,banks2014analysis,guidoboni2009stable,banks2018stable}. 
We computationally examine the performance of the ILE method across a broad range of density ratios and, at least 
for the examples considered herein, we do not detect added mass-related instabilities.
Specifically, as described in the following, we are able to use a simple Dirichlet-Neumann
coupling scheme  \cite{burman2009stabilization,banks2014analysis,bukac2016stability} while achieving reasonable time step sizes and 
avoiding the use of subcycling \cite{hu2001direct, uhlmann2005immersed, yang2008strongly} or other iterative techniques 
for models involving extremely small, nearly equal, equal, and large mass density ratios.

To assess the robustness and accuracy of the proposed
algorithm, results obtained using our ILE method are compared to experimental and computational data 
from benchmarks widely used to test numerical methods for FSI \cite{ahn2006strongly, robertson2003numerical, andersen2005unsteady, mordant2000velocity, ten2002particle, horowitz2010effect}. 
We consider test cases in two and three spatial dimensions involving both smooth and sharp geometries, in various fluid conditions with Reynolds numbering 
reaching up to $\Re=1147$, 
and with different numbers of translational and rotational
degrees of freedom.
Finally, cases demonstrating the application of this methodology to biomedical models are presented, including
to simulating the dynamics of a bileaflet mechanical heart valve in a pulse duplicator 
with a peak Reynolds number of 22600, and to modeling the transport of blood clots in a patient-averaged anatomical model of the inferior vena cava
with a peak Reynolds number of 1500.

%%%%%%%%%%%%%%%%%%%%%%%%%%%%%%%%%%%%%%%%%%%%%%%%%%%%%%%%%%%%%%%%%%%%%%%%%%%%%%%%%%%%%%%%%%%%%%%%%%%%%%%%%%%%%%%%%%%%%%%%%%%
%%%%%%%%%%%%%%%%%%%%%%%%%%%%%%%%%%%%%%%%%%%%%%%%%%%%%%%%%%%%%%%%%%%%%%%%%%%%%%%%%%%%%%%%%%%%%%%%%%%%%%%%%%%%%%%%%%%%%%%%%%%
\section{Continuous Equations of Motion}
\label{sec:model_formulation}

This section outlines our ILE approach to fluid-structure interaction. 
Our methodology builds on a conventional partitioned formulation for FSI, detailed in Sec.~\ref{subsec:partitioned_partitioned}, while 
leveraging a coupling scheme based on the immersed interface method.

%%%%%%%%%%%%%%%%%%%%%%%%%%%%%%%%%%%
\subsection{Partitioned formulation of FSI}
\label{subsec:partitioned_partitioned}
\begin{figure}[t!!!]
	\centering
	\includegraphics[width=0.5\textwidth]{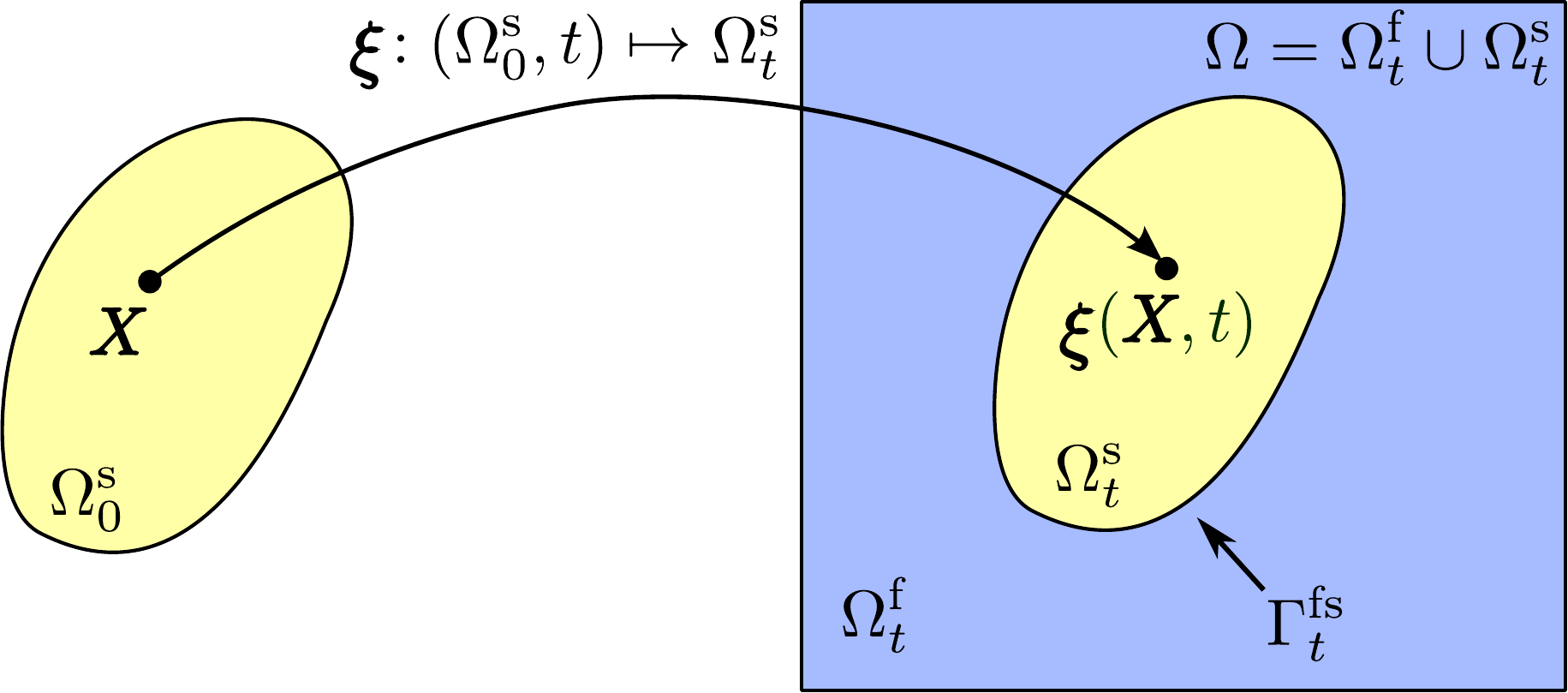}
	\caption{The computational domain $\Omega$ includes time-dependent fluid and solid subdomains, $\Omega^{\textrm{f}}_t$ and $\Omega^{\textrm{s}}_t$. 
	The solid is described using reference coordinates $\s \in \Omega^{\textrm{s}}_0$, and references and current coordinates are connected by the mapping $\Y : (\Omega^{\textrm{s}}_0 , t)\mapsto \Omega^{\textrm{s}}_t $.}
\label{fig:Lag_Eul_schematic} 
\end{figure}

Typical partitioned formulations describe the fluid-structure system occupying a computational domain $\Omega$
via moving fluid and structure subdomains, respectively $\Omega^{\textrm{f}}_t$ and $\Omega^{\textrm{s}}_t$
and indexed by time $t$, so that $\Omega=\Omegaf_t \cup \Omegas_t$; see Fig.~\ref{fig:Lag_Eul_schematic}.
The regions meet along the fluid-structure interface, $\Gamma^{\fs}_t = {\Omega^{\textrm{f}}_t} \cap {\Omega^{\textrm{s}}_t}$.
Fixed physical coordinates are $\x \in \Omega$.
We describe the structural kinematics in Lagrangian form via reference coordinates $\s \in \Omega^{\textrm{s}}_0$ attached to the solid, and we use the motion
map $\Y : (\Omega^{\textrm{s}}_0 , t)\mapsto \Omega^{\textrm{s}}_t $ to determine the physical position of 
a solid material point $\s$ at time $t$.
In the absence of additional loading terms, the equations of motion are
   \begin{align}
       \label{eq:fluid_mom_partitioned} \rhof \, \frac{{\mathrm D} \u}{\mathrm {Dt}}(\x,t)  &= \grad \cdot \sigmaf(\x,t),  &  \x \in \Omega^{\textrm{f}}_t,\\
       \label{eq:continuity_partitioned}  \grad \cdot \u(\x,t) &= 0,  &  \x \in \Omega^{\textrm{f}}_t,\\
 \label{eq:kinematic_cond} \frac{\partial \Y}{\partial t}(\s,t) &= \u(\Y(\s,t),t),  &  \s \in \Gamma^{\fs}_0,\\
 \label{eq:solid_trans_mom_partitioned}  \frac{\mathrm d}{\mathrm{dt}}\int_{\Omega^{\textrm{s}}_0} \rhos \, \frac{\partial\Y}{\partial t}(\s,t)  \, {\mathrm d} \s &=\int_{\Gamma^{\fs}_t} \tauf(\x,t) \, {\mathrm d} a,\\
 \label{eq:solid_ang_mom_partitioned}  \frac{\mathrm d}{\mathrm{dt}}\int_{\Omega^{\textrm{s}}_0} \bigl(\s - \s_{\com}\bigl) \times \bigl(\rhos \, \frac{\partial\Y}{\partial t}(\s,t)\bigl)\, {\mathrm d} \s &= \int_{\Gamma^{\fs}_t} \bigl(\x - \Y_{\com}(t)\bigl) \times  \tauf(\x,t)  \, {\mathrm d} a,
    \end{align} 
in which $\u(\x,t)$ is the fluid velocity, $\rhof$ is the mass density of the fluid, $\rhos$ is the mass density 
of the structure, 
$\Y_{\com}(t)$ and $\s_{\com}$ are, respectively, the center of mass of the structure in the current and reference configurations,
$\sigmaf(\x,t)$ is the fluid stress tensor,
\begin{equation}
 \sigmaf(\x,t) = -p(\x,t) \, \II + \muf \left( \grad\u(\x,t) + \grad\u^T(\x,t) \right),  \,  \x \in \Omega^{\textrm{f}}_t,
\end{equation}
$p(\x,t)$ is the fluid pressure, $\muf$ is the dynamic viscosity of the fluid,  $\tauf(\x,t) = \sigmaf(\x,t) \cdot \vec{n}(\x,t)$ is the fluid traction, and
$\n(\x,t)$ is the unit normal vector pointing into $\Omega^{\textrm{f}}_t$ along $\Gamma^{\fs}_t$.
Eq.~(\ref{eq:fluid_mom_partitioned}) describes the fluid momentum in Eulerian form, Eq.~(\ref{eq:continuity_partitioned}) is the incompressibility constraint, 
Eq.~(\ref{eq:kinematic_cond}) is the kinematic condition along the fluid-structure interface, which implies the no-slip and no-penetration conditions,
and Eqs.~(\ref{eq:solid_trans_mom_partitioned}) and (\ref{eq:solid_ang_mom_partitioned})
describe the dynamics of the linear and angular momentum of the immersed rigid body in Lagrangian form.
Eqs.~(\ref{eq:solid_trans_mom_partitioned}) and (\ref{eq:solid_ang_mom_partitioned}) also account for
the dynamic conditions at the fluid-structure interface because the rigid body forces are balanced by
the fluid traction.

%%%%%%%%%%%%%%%%%%%%%%%%%%%%%%%%%%%%%%%%%%%%%%%%%%%%%%%%%%%%%%%%%%%%%%%%%%%%%%%%%%%%%%%%%%%%%%%%%%%%%%%%%%%%%%%%%%%%%%%%%%%
\subsection{Immersed Lagrangian-Eulerian (ILE) formulation}
\label{subsec:immersed_formulation}

We now introduce an immersed formulation of FSI that describes the same physical model as the partitioned 
formulation detailed in Sec.~\ref{subsec:partitioned_partitioned}. As in that formulation, the
computational domain
is $\Omega$, with $\x \in \Omega$ indicating fixed physical coordinates.
As in Sec.~\ref{subsec:partitioned_partitioned}, the structural kinematics are described in Lagrangian form via reference coordinates $\s \in \Omega^{\textrm{s}}_0$ attached to the solid, and we use the 
same motion map $\Y : (\Omega^{\textrm{s}}_0 , t)\mapsto \Omega^{\textrm{s}}_t $ to determine the physical position of solid material point $\s$ at time $t$.
In the immersed formulation, however, we solve the incompressible Navier-Stokes equations on the full computational domain $\Omega$, 
including both the fluid and solid subdomains. 
We split the computational domain $\Omega$  into an exterior 
fluid region $\Omega^{\textrm{f},+}_t$ and an interior fluid region $\Omega^{\textrm{f},-}_t$,
each parameterized by time $t$, with superscripts ‘$+$’ (‘$-$’) indicating values obtained from the ‘exterior’ (‘interior’) side of the fluid-structure interface.
Using this notation, we have $ \Omega^{\textrm{f},-}_t \equiv \Omegas_t$ and $\Gamma^{\fs}_t = {\Omega^{\textrm{f},+}_t}  \cap  {\Omega^{\textrm{f},-}_t}$.
We extend the definition of the fluid velocity $\u$, pressure $p$, viscosity $\muf$, and stress tensor $\sigmaf$ 
to hold in the entire computational domain $\Omega$, so that the \textit{extended fluid stress tensor} $\sigmaf$ is
\begin{equation}
 \sigmaf(\x,t) = -p(\x,t) \, \II + \muf \left( \grad\u(\x,t) + \grad\u^T(\x,t) \right), \,  \x \in \Omega ={\Omega^{\textrm{f}}_t}  \cup  {\Omega^{\textrm{s}}_t}.
\end{equation}
Our approach applies a singular surface force density along the fluid-structure interface to impose the 
kinematic constraint, which implies a discontinuity in the traction associated with the extended fluid 
stress, $\sigmaf$, along $\Gamma^{\fs}_t$.
A jump in a scalar field $\psi(\x,t)$ at position $\x = \Y(\s,t)$ along the interface is denoted by
\begin{equation}
	\llbracket \psi(\vec{x},t) \rrbracket = \lim_{\epsilon \downarrow 0} \psi(\vec{x} + \epsilon \vec{n}(\x,t),t) - \lim_{\epsilon \downarrow 0} \psi(\vec{x} - \epsilon \vec{n}(\x,t),t)
	= \psi^{+}(\vec{x},t) - \psi^{-}(\vec{x},t),
	\label{eq:jump_definition}
\end{equation}
in which $\llbracket \cdot \rrbracket$ indicates the jump in the value across the interface, $\vec{n}(\x,t)$ is the outward unit normal vector (into the exterior fluid region) 
along $\Gamma^{\fs}_t$,
and $\psi^{+}(\vec{x},t)$ and $\psi^{-}(\vec{x},t)$ are the limiting values approaching the interface position $\vec{x}$ from the exterior 
fluid region $\Omega^{\textrm{f},+}_t$ and interior fictitious fluid region $\Omega^{\textrm{f},-}_t$, respectively.
By considering the jump in the extended fluid stress, the governing equations are %including the jump conditions in the pressure, velocity and first normal derivative of the velocity are written as \cite{kolahdouz2020immersed},
   \begin{align}
       \label{eq:iim_momentum} \rhof \, \frac{{\mathrm D} \u}{\mathrm {Dt}}(\x,t) &= -\grad p(\x,t) + \muf  \, \grad^2 \u(\x,t), &  \x \in \Omega,\\
       \label{eq:iim_continuity}  \grad \cdot \u(\x,t) &= 0, &  \x \in \Omega,\\
       \label{eq:traction_jump_codim1} \llbracket \sigmaf(\Y(\s,t),t) \cdot \vec{n}(\Y(\s,t),t) \rrbracket &= -J^{-1}(\s,t) \, \F(\s,t),  &  \s \in \Gamma^{\fs}_0,\\
		\label{eq:immersed_kinematic_cond} \frac{\partial\Y}{\partial t}(\s,t) &= \vec{u}(\Y(\s,t),t),  &  \s \in \Gamma^{\fs}_0,\\
	 \label{eq:solid_trans_mom_immersed}  \frac{\mathrm d}{\mathrm{dt}}\int_{\Omega^{\textrm{s}}_0} \rhos \, \frac{\partial\Y}{\partial t}(\s,t)  \, {\mathrm d} \s &=\int_{\Gamma^{\fs}_t} \tauf^{+}(\x,t) \, {\mathrm d} a,\\
 \label{eq:solid_ang_mom_immersed}  \frac{\mathrm d}{\mathrm{dt}}\int_{\Omega^{\textrm{s}}_0} \bigl(\s - \s_{\com}\bigl) \times \bigl(\rhos \, \frac{\partial\Y}{\partial t}(\s,t)\bigl)\, {\mathrm d} \s &= \int_{\Gamma^{\fs}_t} \bigl(\x - \Y_{\com}(t)\bigl) \times  \tauf^{+}(\x,t)  \, {\mathrm d} a,
    \end{align}
in which $J(\s,t)$ is the surface Jacobian determinant \cite{kolahdouz2020immersed}, $\F(\s,t)$ is an interfacial surface force density that 
is the Lagrange multiplier to maintain the kinematic condition for the constraint in Eq.~(\ref{eq:immersed_kinematic_cond}) applied along the fluid-solid interface $\Gamma^{\fs}_t$,
and $\tauf^{+}(\x,t)$ is the \textit{exterior} fluid traction.
Importantly, Eqs.~(\ref{eq:solid_trans_mom_immersed}) and (\ref{eq:solid_ang_mom_immersed}), which account for the dynamic interface condition,  imply that only 
the fluid momentum and stresses from the exterior fluid subregion have any physical effect
in driving the dynamics of the structure.

The jump discontinuity in Eq.~(\ref{eq:traction_jump_codim1}) can be decomposed into discontinuities in the pressure and viscous stress, which in current coordinates are
   \begin{align}
		\label{eq:pj_codim1} \llbracket p(\vec{x},t) \rrbracket &=  J^{-1}(\Y^{-1}(\x,t),t) \, \F(\Y^{-1}(\x,t),t) \cdot \vec{n}(\vec{x},t),  &  \x \in \Gamma^{\fs}_t,\\
		\label{eq:du_dn_jump_codim1}     \left\llbracket \muf\frac{\partial\vec{u}}{\partial \vec{n}}(\x,t) \right\rrbracket  &= - (\II - \vec{n}(\x,t)\otimes\vec{n}(\x,t)) \, J^{-1}(\Y^{-1}(\x,t),t) \, \F(\Y^{-1}(\x,t),t),  &  \x \in \Gamma^{\fs}_t.
    \end{align}
Higher order jump conditions, including those associated with the first normal derivative of the pressure and the second normal derivative of the velocity, 
can be also derived \cite{lai2001remark,xu2006systematic}.

\begin{figure}[t!!!]
	\centering
	\includegraphics[width=0.8\textwidth]{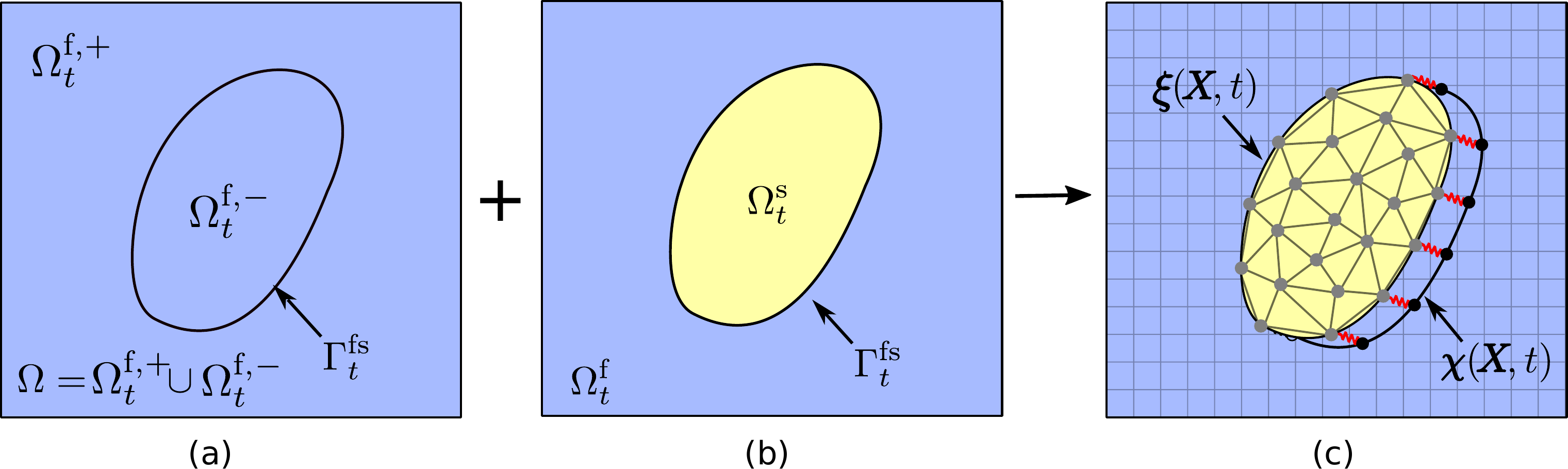}
	\caption{
(a) The immersed interface domain and (b) the partitioned fluid ($\Omega^{\textrm{f}}_t$) and solid ($\Omega^{\textrm{s}}_t$)
subdomains. (c) In the penalty-ILE method, the configuration of the explicit representation of the fluid-structure interface (determined by $\vec{\chi}$) conforms to
 the boundary of the structure (determined by $\vec{\xi}$) in an approximate sense. 
Specifically, the motion of the surface representation is determined by the local fluid velocity at the fluid-structure interface $\Gamma^{\fs}_t$
whereas the no-slip condition is satisfied in an approximate sense by spring-like forces that penalize displacements between the two representations 
of the fluid-structure interface. The displacement between the two representations is exaggerated here for illustration purposes. 
In practice, we always impose an assertion specifying that the maximum relative displacement is at least less than $0.1$ of the Cartesian grid spacing.}
\label{fig:IIM_FSI_schematic} 
\end{figure} 
%%%%%%%%%%%%%%%%%%%%%%%%%%%%%%%%%%%%%%%%%%%%%%%%%%%%%%%%%%%%%%%%%%%%%%%%%%%%%%%%%%%%%%%%%%%%%%%%%%%%%%%%%%%%%%%%%%%%%%%%%%%
\subsection{A penalty approach to the ILE formulation}
\label{subsec:ILE_approach}
The formulation introduced in Sec.~\ref{subsec:immersed_formulation} requires the solution of a saddle-point system that couples the Eulerian and Lagrangian variables \cite{kallemov2016immersed}.
To develop a practical numerical method, we relax the kinematic constraint, Eq.~(\ref{eq:immersed_kinematic_cond}), by introducing two representations of 
the fluid-structure interface 
and applying penalty forces to penalize displacements between the two representations. 
This penalty method determines an approximate Lagrange multiplier force instead of solving for the Lagrange multiplier to exactly impose the condition \cite{kallemov2016immersed}.
Specifically, along with the mapping $\Y(\s,t)$ that determines the kinematics of the
structure, we introduce an explicit representation of the fluid-structure interface that is parameterized via $\X(\s,t)$ and that moves 
with the fluid, so that $\partial\X(\s,t)/\partial t = \vec{u}(\X(\s,t),t)$; see Fig.~\ref{fig:IIM_FSI_schematic}.
In this study, we use a penalty formulation similar to one proposed
 by Goldstein et al.~\cite{goldstein1993modeling},
in which the rigidity constraint is approximately imposed through a linear spring-like force via
\begin{align}
 \F(\s,t) = \kappa\left(\Y(\s,t)-\X(\s,t)\right), \, \s \in \Gamma^{\fs}_0.
 \label{eq:penalty_force}
\end{align}
Here, $\kappa$ is the spring stiffness penalty parameter. This force penalizes deviations from
the constraint $ \Y(\s,t) =  \X(\s,t)$ and, in the discretized equations, acts to ensure that the volumetric and surface meshes are at least approximately conformal in their motion.
Note that
it is possible to control the discrepancy between the two configurations because
 as $\kappa \rightarrow \infty$, the formulation exactly imposes the constraint that the two interfaces representations move together.
 
%%%%%%%%%%%%%%%%%%%%%%%%%%%%%%%%%%%%%%%%%%%%%%%%%%%%%%%%%%%%%%%%%%%%%%%%%%%%%%%%%%%%%%%%%%%%%%%%%%%%%%%%%%%%%%%%%%%%%%%%%%%
\subsection{Rigid body dynamics}
\label{subsec:RBD}
The general formulation for description of kinematics of a rigid body includes both 
translational and rotational motions of material points in the body.
The kinematics in the current configuration coordinates of the solid can be written 
in terms of reference coordinates as
\begin{equation}
\Y(\s,t)=\QQ(t) \, \bigl(\s -\s_{\com}\bigl) \, + \, \s_{\com} \, + \, \disp_{\com}(t),
\label{eq:y_current_RBD}
\end{equation}
in which $\s_{\com}$ is the position of the center of mass in the reference coordinates, $\disp_{\com}(t)$ is the displacement of the center of mass,
and $\QQ(t)$ is the rotation matrix.
The rotation matrix $\QQ(t)$ can be expressed using the Euler angles, $\phi = \phi(t)$, $\theta = \theta(t)$, and $\psi = \psi(t)$, as
\begin{equation}
\QQ(t)=
\begin{pmatrix} 
\cos(\psi)\cos(\theta) & -\sin(\phi)\cos(\phi)+\cos(\psi)\sin(\theta)\sin(\phi) &  \sin(\phi)\sin(\phi)+\cos(\psi)\sin(\theta)\cos(\phi) \\ 
\sin(\psi)\cos(\theta) & \cos(\phi)\cos(\phi)+\sin(\psi)\sin(\theta)\sin(\phi)  & -\cos(\phi)\sin(\phi)+\sin(\psi)\sin(\theta)\cos(\phi) \\ 
-\sin(\theta) & \cos(\theta)\sin(\phi) & \cos(\theta)\cos(\phi)  
\end{pmatrix},
\end{equation}
here using the $x$-$y$-$z$~convention in the rotation order. Following the approach of Akkerman et al.~\cite{akkerman2012free}, we avoid
explicitly forming the Euler angles
except for cases where there is only one rotational degree of freedom.
As will be presented in Sec.~\ref{sec:numerical_algorithm}, 
in the time integration of the rigid body dynamics, the rotation matrix is an additional problem unknown
that is integrated in time along with the displacement and rotational degrees of freedom.
This significantly reduces the complexity of the calculation of the rotation angle for three degree of freedom (3-DOF) problems in two spatial dimensions and six degree of freedom (6-DOF) problems in three spatial dimensions.
Using Eq.~(\ref{eq:y_current_RBD}), the rigid body displacement and velocity are respectively defined as
\begin{equation}
\disp(\s,t) = (\QQ(t) - \II)\bigl(\s-\s_{\com}\bigl) \, + \, \disp_{\com}(t)
\label{eq:dis_current}
\end{equation}
and
\begin{equation}
\label{eq:rbd_vel_ref}
\dispdot(\s,t)=\dot{\QQ}(t) \, \bigl(\s-\s_{\com}\bigl) \, + \, \dispdot_{\com}(t).
\end{equation}
Denoting the rigid body center of mass in the current configuration as $\Y_{\com}(t) = \Y(\s_{\com},t)=\s_{\com} \, + \, \disp_{\com}(t)$, 
Eq.~(\ref{eq:rbd_vel_ref}) becomes
\begin{equation}
\label{eq:rbd_vel_current}
\dispdot(\s,t) =\dot{\QQ}(t)\QQ^{-1}(t)\bigl(\Y(\s,t)-\Y_{\com}(t)\bigl) \, + \, \dispdot_{\com}(t) =\OO(t) \, \bigl(\Y(\s,t)-\Y_{\com}(t)\bigl) \, + \, \dispdot_{\com}(t),
\end{equation}
in which
\begin{equation}
\OO(t)=\dot{\QQ}(t)\QQ^{-1}(t)=
\begin{bmatrix}
0 & -\omega_3(t) & \omega_2(t) \\
\omega_3(t) & 0  & -\omega_1(t) \\
-\omega_2(t) & \omega_1(t) & 0 
\end{bmatrix}
\end{equation}
is the skew-symmetric angular velocity matrix. The pseudovector $\vec{\omega}(t)=(\omega_1(t),\omega_2(t),\omega_3(t))$ can be extracted from this tensor matrix.
Using the angular velocity vector, the rigid body velocity field in Eq.~(\ref{eq:rbd_vel_current}) can be re-written in terms of $\vec{\omega}(t)$ as,
\begin{equation}
\label{eq:rbd_vel_omega}
\dispdot(\s,t) = \vec{\omega}(t) \, \times \, \bigl(\Y(\s,t)-\Y_{\com}(t)\bigl) \, + \, \dispdot_{\com}(t).
\end{equation}
For a three-dimensional problem, the three components of the translational velocity of the center of mass $\dispdot_{\com}(t)$, together
with the three components of the angular velocity $\vec{\omega}(t)$, form the six degrees-of-freedom that completely determine
the kinematics of the rigid body.

Linear and angular momentum conservation in the rigid body are described by a system of six ordinary differential equations,
\begin{align}
 \label{eq:linear_momentum_RB} \frac{\mathrm d}{\mathrm{dt}}\big(m \, \dispdot_{\com}(t)\big) &= \F^{\textrm{net}}(t) ,\\
 \label{eq:angular_momentum_RB} \frac{\mathrm d}{\mathrm{dt}}\big(\JJ(t) \, \vec{\omega}(t)\big) &=\T^{\textrm{net}}(t) ,
\end{align}
in which $m$ is the mass of the solid object, $\F^{\textrm{net}}(t)$ is the global force including the sum of all the forces exerted on the rigid body, and $\T^{\textrm{net}}(t)$ is the net torque vector.
$\JJ(t)$ is the inertia tensor of the solid body in the current configuration, which is defined in terms of the inertia tensor in the reference configuration, $\JJ_0$, via
\begin{equation}
 \label{eq:QJoQ} \JJ(t) = \QQ(t) \, \JJ_0 \, \QQ^T(t),
\end{equation}
in which $\JJ_0$ is
\begin{equation}
\JJ_0 = \int_{\Omegas_0} \rhos \, \bigl(\s - \s_{\com}\bigl)\cdot\bigl(\s-\s_{\com}\bigl) \II \,  {\mathrm d} \s -\int_{\Omegas_0} \rhos \, \bigl(\s-\s_{\com}\bigl)\otimes\bigl(\s-\s_{\com}\bigl) \, {\mathrm d} \s .
\end{equation}
In our fluid-structure interaction framework, the net force and torque vectors in Eqs.~(\ref{eq:linear_momentum_RB}) and (\ref{eq:angular_momentum_RB}) 
are computed as
\begin{align}
 \label{eq:force_RB} \F^{\textrm{net}} &= (\rhos - \rhof)V\g \, + \, \int_{\Gamma^{\fs}_t} \tauf^{+}(\x,t) \, {\mathrm d} a,\\
 \label{eq:torque_RB} \T^{\textrm{net}} &=\int_{\Gamma^{\fs}_t} \bigl(\x - \Y_{\com}(t)\bigl) \, \times \, \tauf^{+}\bigl(\x,t \bigl) \, {\mathrm d} a,
\end{align}
in which $V$ is the volume of the solid object, $\g$ is the gravity vector, and $\tauf^{+}(\x,t)$ is 
the exterior fluid traction vector exerted on the solid object by the fluid.
In the penalty formulation of the ILE method, quantities on the left hand sides of Eqs.~(\ref{eq:force_RB}) and (\ref{eq:torque_RB}) 
are evaluated on $\Gamma^{\fs}_t$, which moves with the local fluid velocity.

%%%%%%%%%%%%%%%%%%%%%%%%%%%%%%%%%%%%%%%%%%%%%%%%%%%%%%%%%%%%%%%%%%%%%%%%%%%%%%%%%%%%%%%%%%%%%%%%%%%%%%%%%%%%%%%%%%%%%%%%%%%
\section{Discrete Equations of Motion}
\label{sec:numerical_algorithm}

This section introduces numerical methods for the penalty formulation of the immersed Lagrangian-Eulerian method detailed in Sec.~\ref{subsec:ILE_approach}. 
This approach leverages our immersed interface method (IIM) for discrete surface representations \cite{kolahdouz2020immersed}. 
We include only the key aspects of this method; for additional details and benchmarking studies of problems with prescribed motion, see Kolahdouz et al. \cite{kolahdouz2020immersed}. 
Standard methods are used for the rigid body dynamics. 
A second order accurate Strang time step splitting approach \cite{strang1968construction}  is used to obtain systems of equations that can be treated via efficient linear solvers.

\subsection{Eulerian discretization}
\label{subsec:methodology-FD}
The incompressible Navier-Stokes equations are discretized on an adaptively refined Cartesian grid using
a staggered-grid finite difference discretization \cite{griffith2009} in which the pressure is approximated at the centers of the Cartesian grid cells and the components of the velocity
are approximated at the centers of the edges (in two spatial dimensions) or faces (in three spatial dimensions) of the grid cells.
Standard compact second-order accurate differencing schemes are used for the divergence, gradient, and Laplace operators. 
The discrete divergence of the velocity $\vec{D} \cdot \u$ is evaluated at the cell centers, whereas the discrete pressure gradient $\vec{G} p$ and 
the components of the discrete Laplacian of the velocity $L \vec{u}$ are 
evaluated at the cell edges (in two spatial dimensions) or faces (in three spatial dimensions).
For the nonlinear advection terms,
a staggered-grid version \cite{griffith2012volume,griffith2009} of the xsPPM7 variant \cite{rider2007accurate} of the piecewise parabolic method (PPM) \cite{colella1984piecewise} is used.
Physical boundary conditions are prescribed along the boundaries of the computational domain $\Omega$ as described previously \cite{griffith2009,griffith2012immersed}. 
Adaptive computations use a discretization approach described by Griffith \cite{griffith2012immersed} that employs Cartesian grid adaptive mesh refinement (AMR).

To account for the jump conditions along the fluid-solid interface that occur in the ILE formulation, we modify the definitions of
$\vec{G} p$  and $L \vec{u}$ for those stencils
that cross the immersed interface.
Using generalized Taylor series expansions \cite{li2001immersed,xu2006systematic}, it can be shown that if the interface
 cuts between two Cartesian grid points at location $\x_{\circ}=(x_{\circ},y_{\circ},z_{\circ})$, such that $x_{i,j,k}\leq x_{\circ} < x_{i+1,j,k}$, with $\x_{i,j,k}\in \Omega^{\textrm{f},-}$ and $\x_{i+1,j,k}\in \Omega^{\textrm{f},+}$, then for a piecewise differentiable quantity $\psi$, we have
\begin{align}
\label{eq:p_jump_imposition} \frac{\partial\psi}{\partial x}(\x_{i+\half,j,k}) &=\frac{\psi_{i+1,j,k}-\psi_{i,j,k}}{\Delta x}+\frac{\mathrm{sgn}\{n^x\}}{\Delta x}\displaystyle \sum_{m=0}^{2}\frac{(d^{+})^{m}}{m!}{\left\llbracket \frac{\partial^{m}\psi}{\partial x^{m}}\right\rrbracket}_{\x_{\circ}} + O({\Delta x}^{2}), \\
\label{eq:u_jump_imposition} \frac{\partial^{2}\psi}{\partial x^{2}}(\x_{i,j,k}) &=\frac{\psi_{i+1,j,k}-2\psi_{i,j,k}+\psi_{i-1,j,k}}{\Delta x^{2}}+\frac{\mathrm{sgn}\{n^x\}}{{\Delta x}^{2}}{\displaystyle \sum_{m=0}^{3}\frac{(d^{+})^{m}}{m!}{\left\llbracket \frac{\partial^{m}\psi}{\partial x^{m}}\right\rrbracket}_{\x_{\circ}} + O({\Delta x}^{2})},
\end{align}
in which $\Delta x$ is the grid spacing in the $x$ direction, $\psi_{i,j,k}=\psi(\x_{i,j,k})$, $d^{+}= x_{i+1,j,k}-x_{\circ} > 0$, and $n^x$ is the $x$-component of the normal vector $\vec{n}=(n^x,n^y,n^z)$ at the intersection point $\x_\circ$.
The full implementation of this approach to the three-dimensional incompressible Navier-Stokes equations,
including the application of the jump  corrections to the stencils of the pressure and the viscous terms and algorithms for identifying intersections between the finite difference stencils and 
the discrete interface representation, is detailed in our earlier work \cite{kolahdouz2020immersed}.

\subsection{Lagrangian discretization}
\label{subsec:methodology-FE}

Let $\cT_h$ be a triangulation of $\Omegas_0$, the reference configuration of the volumetric rigid body, 
composed of elements $U^e$ such that $\cT_h = \cup_{e} U^e$, in which $e$ indexes the mesh elements. We take
$\left\{\s_l\right\}_{l=1}^{M}$ to be the positions of the $M$ nodes of the mesh in the reference 
configuration, $\left\{\Y_l(t)\right\}_{l=1}^M$ 
to be the current positions of the nodes, and $\left\{\phi_l(\s)\right\}_{l=1}^{M}$ 
to be the corresponding interpolatory nodal (Lagrangian) basis functions.
A continuous description of the configuration of the structure is defined by
\begin{equation}
 \label{eq:discrete_deformation} \Y_{h}(\s,t) ={\displaystyle \sum_{l=1}^{M}\Y_{l}(t)\phi_{l}(\s)}, \, \s \in \Omegas_0.
\end{equation}

The configuration of the fluid-structure interface representation that moves with the fluid is described by the 
mapping $\X:(\Gamma_0^{\fs},t)\mapsto\Gamma_t^{\fs}$. 
To obtain a discrete representation of this interface, we use a surface mesh that corresponds to the restriction of the volumetric 
solid mesh to $\partial\Omega^s_0 \equiv \Gamma_0^{\fs}$. 
For the discrete representation of the fluid-structure interface we have
\begin{equation}
 \label{eq:discrete_deformation_surface} \X_{h}(\s,t) ={\displaystyle \sum_{l=1}^{M}\X_{l}(t)\phi_{l}(\s)},  \, \s \in \Gamma_0^{\fs},
\end{equation}
except that in practice, we only need to evaluate this sum over the lower-dimensional subset of nodes that are located on surface 
mesh, and not over all of the nodes in the volumetric representation. The reason is that the interpolatory~(nodal) basis functions associated with interior nodes 
vanish on the surface. Similarly, the surface force density is determined by
 \begin{equation}
\F_{h} (\s,t) ={\displaystyle \sum_{l=1}^{M}\F_{l}(t)\phi_{l}(\s)},  \, \s \in \Gamma_0^{\fs}.
\end{equation}
Again, this sum only needs to be evaluated using the $M^\text{fs}$ surface nodes. In an implementation, it is straightforward to use separate data structures 
for the volumetric and surface structural meshes. For the remainder of the paper, we adopt the convention that all computations involving the 
surface representation are performed using only a surface mesh with 
appropriate nodal degrees of freedom and surface-restricted basis functions.

Stress jump conditions are imposed by evaluating the correction terms from the interfacial forces and interface configuration 
(i.e.~generalizations of Eqs.~(\ref{eq:p_jump_imposition}) and (\ref{eq:u_jump_imposition})).
Geometrical quantities, including the surface normals and surface Jacobian determinant, that are needed by 
the IIM discretization, are obtained by directly differentiating Eq.~(\ref{eq:discrete_deformation}).
Note, however, that the standard nodal basis functions are $C^0$ but not $C^1$ at element boundaries, and so quantities that are obtained in terms of ${\partial \X_h}/{\partial \s}$ are discontinuous in both the reference and current configurations.
In particular, the pointwise jump conditions determined from the mesh geometry and the surface Jacobian $J$ are generally discontinuous between the elements.
Following the approach introduced in our prior work \cite{kolahdouz2020immersed}, we obtain a continuous approximation to the jump conditions through the $L^2$ projection.
Briefly, given a function $\psi \in L^{2}(\Gamma_0^{\fs})$, its $L^{2}$ projection $P_{h} \psi$ onto the subspace 
$S_h = \mathrm{span}\{\phi_l(\s)\}_{l=1}^{{M}^\text{fs}}$ is defined by requiring $P_{h} \psi$ to satisfy
\begin{equation}
   \int_{\Gamma_0^{\fs}} \big(\psi(\s) - P_{h}\psi(\s)\big) \, \phi_l(\s) \,  {\mathrm d} A = 0, \quad \forall l=1,\ldots,{M}^\text{fs}.
   \label{eq:L2_proj_def}
\end{equation}
The $L^{2}$ projection of a vector-valued quantity is determined component-wise.
Because the $L^{2}$ projection is defined via integration, the function $\psi$ does not need to be continuous or even to have well-defined nodal values.
By construction, however, $P_h \psi$ will inherit any smoothness provided by the subspace $S_h$.
In particular, for $C^0$ Lagrangian basis functions, $P_h \psi$ will be at least continuous.
In our numerical scheme, we separately compute the projection of the normal and tangential components of the surface force per unit current area, $J^{-1} \F_h(\s,t)$, onto $S_h$,
as needed to specify the conditions for the pressure and the viscous stress.
We drop the subscript ``$h$" in the remainder of the paper to simplify notation.
To solve for the projected jump conditions, linear systems of equations involving the mass matrix ${\mathcal M}$ need to be solved, in which ${\mathcal M}$ has
components ${\mathcal M}_{kl} = \int \phi_k \phi_l \, \DA$.
Eq.~\eqref{eq:L2_proj_def} is evaluated using seventh-order Gaussian quadrature. 
Notice that these projections are computed only along the fluid-solid interface and involve only surface degrees of freedom. Consequently, the computational cost of evaluating these projections is 
asymptotically smaller than the solution of the fluid equations.
Note also that similar to the conventional IB method, a force-spreading operator $\vec{\cS} = \vec{\cS}[\X]$ can be defined to evaluate and apply 
the correction terms $\vec{\cS}[\X]\F$ to the Eulerian discretization via a discrete Eulerian 
force density $\f = \vec{\cS}[\X]\F$ \cite{kolahdouz2020immersed}.

The velocity of the fluid-structure interface representation that moves with the fluid is determined by evaluating the Eulerian velocity $\u(\x,t)$ 
on the interface. 
As detailed previously \cite{kolahdouz2020immersed}, 
it is possible to interpolate the discretized Eulerian velocity field $\u$ to the Lagrangian 
interface mesh using a corrected  bilinear (or, in three spatial dimensions, trilinear) interpolation scheme that accounts for the known 
discontinuities in $\partial\vec{u}/\partial \vec{n}$. In general, however, the basic interpolation scheme will produce 
an interface velocity field that is not in the space of the nodal basis functions, which implies that it cannot be used directly 
to update the configuration of the interface. 
To obtain a suitable surface mesh velocity field, we project the interpolated velocity field 
onto the space spanned by the nodal basis functions using Gaussian quadrature. This procedure implicitly 
defines a velocity-restriction operator $\cJ = \cJ[\X,\F]$, so that $\U = \cJ[\X,\F] \, \u$.

The FSI coupling approach used herein crucially relies on the accurate evaluation of the exterior fluid traction.
This requires evaluating the exterior fluid pressure and exterior viscous shear stress.
To evaluate the exterior pressure at a position $\x \in \Gamma^{\fs}_t$, we use
\begin{equation}
 p^{+}_h(\vec{x},t) =  \llbracket p(\vec{x},t) \rrbracket + \cI[p](\vec{x}^{-},t),
 \end{equation}
in which $p^{-}=\cI[p](\vec{x}^{-},t)$ is the interior pressure interpolated to a position $\x^{-}$ away from the interface in the opposite direction 
of the normal vector $\n$ and at a distance equal to $1.3$ times the diagonal size of one grid cell. 
This factor has been chosen on an empirical basis \cite{kolahdouz2020immersed}.
Here, $\cI$ is the unmodified bilinear (or trilinear) interpolation operator involving quantities on one side of the interface.
To evaluate the exterior wall shear stress, a one-sided approximation to the normal derivative of the velocity is calculated using the same interfacial velocity reconstruction 
procedure that is used to determine the interface velocity along with the velocity value at a neighboring location in the direction of the normal vector $\vec{x}^{+}$.
As with the pressure, unmodified bilinear (or trilinear) interpolation is used to obtain the velocity away from the interface.
A one-sided finite difference formula is used to calculate the normal derivative,
\begin{equation}
 \left(\frac{\partial \u}{\partial \vec{n}}\right)^{+}_h(\vec{x},t) = \frac{\cI[\u](\vec{x}^{+},t)-\u(\vec{x},t)}{\hat{h}},
\end{equation}
in which the distance $\hat{h}$ is chosen to be slightly larger than the diagonal size of the Cartesian mesh ($1.05$ times the diagonal size), 
so that regular bilinear (or trilinear) interpolation can be used to evaluate $\cI[\u](\vec{x}^{+},t, \hat{h})$ ensuring that the interpolation only uses values on one side of the interface.
It is possible to use a second-order formula with a three-point stencil which requires interpolating 
an additional point in the normal direction, but preliminary numerical experiments (data not shown) suggest the computation using only two points 
is more stable.
Moreover, as shown previously \cite{kolahdouz2020immersed}, this simple scheme is adequate to achieve a point-wise first-order accurate approximation to the wall shear stress.
Note that as with velocity interpolation, the pressure and wall shear stress can be evaluated at arbitrary locations along the 
interface. We also use a surface-restricted $L^2$ projection to obtain nodal values of these interfacial quantities.

\subsection{Time integration scheme}
\label{subsec:time_integration}

In advancing from time step $n$ at time $t$ to time step $n+1$ at time $t + \dt$,
we define a vector of variables $\vec{\Upsilon}$ that includes all of the Eulerian and Lagrangian quantities that need to be updated
as $\vec{\Upsilon}=[\u,  p, \X, \disp_{\com}, \dispdot_{\com}, \vec{\omega}, \QQ, \JJ]$.
We use the second 
order Strang splitting scheme \cite{strang1968construction}, in which within three steps we: 
1) solve the rigid body dynamics equations over a half time step $\dt/2$, treating the fluid traction as constant in time;
2) solve the IIM/FSI equations over a full time step $\dt$, treating the configuration of the solid as constant in time; and
3) solve the rigid body dynamics equations over a final half time step $\dt/2$, treating the fluid traction as constant in time.
The details of the time stepping for the rigid body dynamics equations, the IIM/FSI equations, and the overall algorithm are given below.

\subsubsection{Rigid body time integration scheme}
Although we ultimately employ a time step splitting approach that advances the rigid body configuration in two half-steps, 
to simplify the discussion, the approach for a fixed time step size $\dt$ is detailed here.
Denote by $\cL_{\dt}$ the action of the Lagrangian rigid body dynamics solver over the time increment $\dt$ that
acts on a solution vector $\vec{\Upsilon}$. This solution vector includes all of the Eulerian and Lagrangian variables, but only advances 
the volumetric structural variables $\disp_{\com}$, $\dispdot_{\com}$, $\vec{\omega}$, $\QQ$, and $\JJ$ while keeping the remaining variables 
fixed.
Briefly, discretizations of Eqs.~(\ref{eq:angular_momentum_RB})--(\ref{eq:QJoQ}) are solved over the time increment $\dt$ to obtain
${\disp^{n+1}_{\com}}$, ${\dispdot^{n+1}_{\com}}$, 
$\vec{\omega}^{n+1,k}$, and $\QQ^{n+1,k}$ via
\begin{align}
\label{eq:time_int_step1} m \, \frac{{{\dispdot^{n+1}_{\com}}} - {\dispdot^{n}_{\com}}}{\dt}&=\int_{\Gamma^{\fs}_t} \tauf^{+}(\x) \, {\mathrm d} a,\\
\label{eq:time_int_step2} \frac{{{\disp^{n+1}_{\com}}}- {\disp^{n}_{\com}}}{\dt}&=\frac{1}{2}\big({{\dispdot^{n+1}_{\com}}} + {{\dispdot_{\com}}}^{n}\big),\\
 \label{eq:time_int_step3}  \frac{\QQ^{n+1,k}\JJ_0 \, (\QQ^{n+1,k})^T \, {\vec{\omega}}^{n+1,k}-{\QQ^{n}}\JJ_0{(\QQ^n)^{T} {\vec{\omega}}^{n}}}{\dt}&= \int_{\Gamma^{\fs}_t} (\x - {\Y_{\com}}) \times  \tauf^{+}(\x)  \, {\mathrm d} a,\\
\label{eq:time_int_step4}  \frac{\QQ^{n+1,k} - {\QQ}^{n}}{\dt}&=\frac{1}{4}(\OO^{n+1,k} + \OO^{n})(\QQ^{n+1,k} + {\QQ}^{n}).
\end{align}
Eqs.~(\ref{eq:time_int_step3}) and (\ref{eq:time_int_step4}) can be solved through a few subiterations
to obtain the new rotation matrix $\QQ^{n+1}$ and 
the angular velocity ${\vec{\omega}}^{n+1}$; starting from $\QQ^{n+1,k=1}\equiv\QQ^{n}$, and $\vec{\omega}^{n+1,k=1}\equiv\vec{\omega}^{n}$ 
at each time step, subiterations are performed by looping over $k$ until $\norm{\QQ^{n+1,k}-\QQ^{n+1,k-1}}_{\infty}\leq \epsilon$ or
 $\norm{\vec{\omega}^{n+1,k}-\vec{\omega}^{n+1,k-1}}_{\infty}\leq \epsilon$ 
with $\epsilon=10^{-8}$. These iterations are inexpensive, and between one and three are typically needed to reach the convergence criteria. 
Seventh-order Gaussian quadrature rules are used to approximate integrals on the left side of Eqs.~(\ref{eq:time_int_step2}) and (\ref{eq:time_int_step4}).
The updated configuration of $\Y(\s,t)$ is then evaluated using Eq.~(\ref{eq:y_current_RBD}). 
Similarly, the moment of inertia is determined using a discretization of Eq.~(\ref{eq:QJoQ}),
\begin{equation}
\JJ^{n+1} =\QQ^{n+1}\JJ_0(\QQ^{n+1})^T.
\end{equation}

\subsubsection{IIM time integration scheme}
\label{subsec:time-integration}

Starting from $\X^{n}$ and $\vec{u}^{n}$ at time $t^n$ and $p^{n-\frac{1}{2}}$ at 
time $t^{n-\frac{1}{2}}$, we compute $\X^{n+1}$, $\u^{n+1}$, and $p^{n+\frac{1}{2}}$. 
Denote by $\cE_{\dt}$ the action of the IIM solver over a full time step that
acts on the solution vector $\vec{\Upsilon}$ that includes all of the Eulerian and Lagrangian variables but only advances $\u$,  $p$, and $\X$.
Briefly, using the discrete velocity restriction operator $\cJ^{n}=\cJ[\X^{n},\F^{n}]$, we first obtain 
initial approximations to the interface position at time $t^{n+\frac{1}{2}}$  via
\begin{align}
	\hat{\X}^{n+1} = \X^{n} + \frac{\dt}{2} \cJ^{n} \vec{u}^{n}, \\
	\X^{n+\frac{1}{2}}=\frac{\hat{\X}^{n+1}+\X^{n}}{2}.
\end{align}
Next, we solve for $\X^{n+1}$, $\u^{n+1}$, and $p^{n+\frac{1}{2}}$ via
\begin{align}
\label{eq:time-integ-momentum}\rho\left(\frac{\u^{n+1}-\u^{n}}{\Delta t}+\vec{A}^{n+\frac{1}{2}}\right) &=-\vec{G} \, p^{n+\frac{1}{2}}+\muf \, \vec{L}\left(\frac{\u^{n+1} + \u^{n}}{2}\right)+ \vec{\cS}\, \F^{n+\frac{1}{2}}, \\
\vec{D} \cdot \vec{u}^{n+1} &=0, \\
\label{eq:update-X-FSI} \frac{\X^{n+1}-\X^{n}}{\Delta t} &= \U^{n+\half} = \cJ^{n+\half} \left(\frac{\u^{n+1}+\u^{n}}{2}\right),
\end{align}
in which $\vec{A}^{n+\frac{1}{2}}=\frac{3}{2}\vec{A}^{n}-\frac{1}{2}\vec{A}^{n-1}$ is obtained from a high-order upwind spatial discretization
of the nonlinear convective term $\u \cdot \grad \u$ \cite{griffith2009}, and
$\cJ^{n+\half}=\cJ[\X^{n+\half},\F^{n+\half}]$ is the velocity restriction operator at the half time step configuration.
This time stepping scheme requires only linear solvers for the time-dependent incompressible Stokes equations.
We solve this system of equations by the flexible GMRES (FGMRES) algorithm with a preconditioner based on the projection method that uses inexact subdomain solvers \cite{griffith2009}.
In the initial time step, a two-step predictor-corrector method is used to determine the velocity, deformation, and pressure; see Griffith and Luo \cite{BEGriffith17-ibfe} for further details.

\subsubsection{Fluid-structure interaction time stepping scheme}

For the solution vector $\vec{\Upsilon}$ over one time step using the Strang splitting scheme we have,
\begin{equation}
\vec{\Upsilon}^{n+1} = \cL_{\Delta t/2} \cE_{\Delta t}\cL_{\Delta t/2}\vec{\Upsilon}^{n}.
\end{equation}
With this time-staggered approach, the overall algorithm to solve the fluid-structure problem 
is
	\\
\begin{steps}
	\item{ 
		Solve the rigid body dynamics in Eqs.~(\ref{eq:time_int_step1})--(\ref{eq:time_int_step4}) 
		over a half time step $\Delta t/2$ to advance from step $n$ to $n+\frac{1}{2}$ and obtain 
		${\disp_{\com}}^{n+\frac{1}{2}}$, ${\dispdot_{\com}}^{n+\frac{1}{2}}$, ${\vec{\omega}}^{n+\frac{1}{2}}$ and ${\QQ}^{n+\frac{1}{2}}$,
		and the new position of the volumetric Lagrangian mesh.
		}
	\item{ 
		Calculate the penalty force using the most recent position of the volumetric solid mesh and the surface mesh
		that moves with the fluid. 
	}
	\item{ 
		Solve for the IIM subsystem in Eqs.~(\ref{eq:time-integ-momentum})--(\ref{eq:update-X-FSI}) over a full time step 
		and obtain the updated Eulerian velocity field $\vec{u}^{n+1}$ and pressure $p^{n+\frac{1}{2}}$ as well as the 
		the Lagrangian velocity $\U^{n+\frac{1}{2}}$ and positions $\X^{n+1}$ of the surface mesh and the exterior 
		fluid traction forces $\tauf^{+}$.
	}
	\item{ 
		Using the exterior fluid traction force $\tauf^{+}$ from Step III, solve the rigid body dynamics
		 in Eqs.~(\ref{eq:time_int_step1})--(\ref{eq:time_int_step4}) 
		over a half time step $\Delta t/2$ to advance from step $n+\frac{1}{2}$ to $n+1$ and obtain 
		${{\disp_{\com}}}^{n+1}$, ${{\dispdot_{\com}}}^{n+1}$, ${\vec{\omega}}^{n+1}$ and ${\QQ}^{n+1}$.
	}
	\item{ 
			Move the Lagrangian mesh of the bulk solid and obtain the new positions $\Y^{n+1}$.	}
\end{steps}

%%%%%%%%%%%%%%%%%%%%%%%%%%%%%%%%%%%%%%%%%%%%%%%%%%%%%%%%%%%%%%%%%%%%%%%%%%%%%%%%%%%%%%%%%%%%%%%%%%%%%%%%%%%%%%%%%%%%%%%%%%%
\section{Numerical examples}
\label{sec:Results}
This section presents computational examples to characterize the performance of the present methodology in two and three spatial dimensions. 
As a demonstration of the method's ability to tackle more challenging problems,
applications to two biomedical models are also presented, 
including simulations of the dynamics of a bileaflet mechanical heart valve 
in a pulse duplicator system and the transport of blood clots in a patient-averaged anatomical model of the inferior vena cava.
Where possible, comparisons are made to available experimental or computational results. 
We begin by considering model problems involving a limited number of translational, but not rotational, degrees of freedom (DOF). 
We systematically increase the complexity of the tests by incorporating additional translation and rotational degrees of freedom. We also consider both smooth 
immersed structures as well as structures with sharp corners.
The fluid-solid interface representation that moves with the fluid is discretized 
either by piecewise-linear ($P^1$) elements for two-dimensional cases, or by piecewise linear ($P^1$) or piecewise bilinear ($Q^1$)
elements for three-dimensional cases.
Unless otherwise noted, the structural meshes are constructed so that the ratio of the Lagrangian element size to the Eulerian grid spacing,
denoted by $\Mfac$, is $\Mfac \approx 2$ at least along the fluid-structure interface.
In all cases, the Eulerian domain is discretized using an adaptively refined grid.
The Cartesian grid spacing on the finest level of the locally refined grid is $h_{\textrm{finest}}=r^{-(N-1)}h_{\textrm{coarsest}}$,
in which $h_{\textrm{coarsest}}$ is the grid spacing on the coarsest level, $r$ is the refinement ratio, and $N$ is the number of refinement levels.
The spring penalty parameters are computationally determined as approximately the largest values allowed by our explicit time stepping algorithm at the time step sizes 
used in the simulations. In each case, the penalty parameter values are tuned using  bisection.
Unless otherwise noted, centimeter-gram-second (CGS) units are used.
For models involving gravitational forcing, gravitational acceleration is set to $g = 981 \, \textrm{cm}\cdot\textrm{s}^{-2}$.
The large scale three-dimensional simulations in Secs.~\ref{subsec:falling_dense_sphere}, \ref{subsec:rising_sphere}, \ref{subsec:MHV}, 
and \ref{subsec:ivc_clot} were performed using the \textit{Dogwood} cluster provided by the Research Computing Division of University of North Carolina at
Chapel Hill Information Technology Services.
Each node is comprised of 2.4 GHz Intel Xeon E5-2699Av4 processors with the Broadwell-EP micro-architecture, 512 GB memory, and 44 cores per node, and nodes are connected by a high bandwidth Infini-band EDR switching fabric.
The smaller benchmark models are run on a workstation with two 32-core Intel Xeon E5-2680 v3 2.5~GHz processors and 32~GB of memory.

%%%%%%%%%%%%%%%%%%%%%%%%%%%%%%%%%%%%%%%%%%%%%%%%%%%%%%%%%%%%%%%%%%%%%%%%%%%%%%%%%%%%%%%%%%%%%%%%%%%%%%%%%%%%%%%%%%%%%%%%%%%
\subsection{Vortex-induced vibration of a cylinder}
\label{subsec:viv-cylinder}

 \begin{figure}[b!!!]
		\centering
			\includegraphics[width=0.9\textwidth]{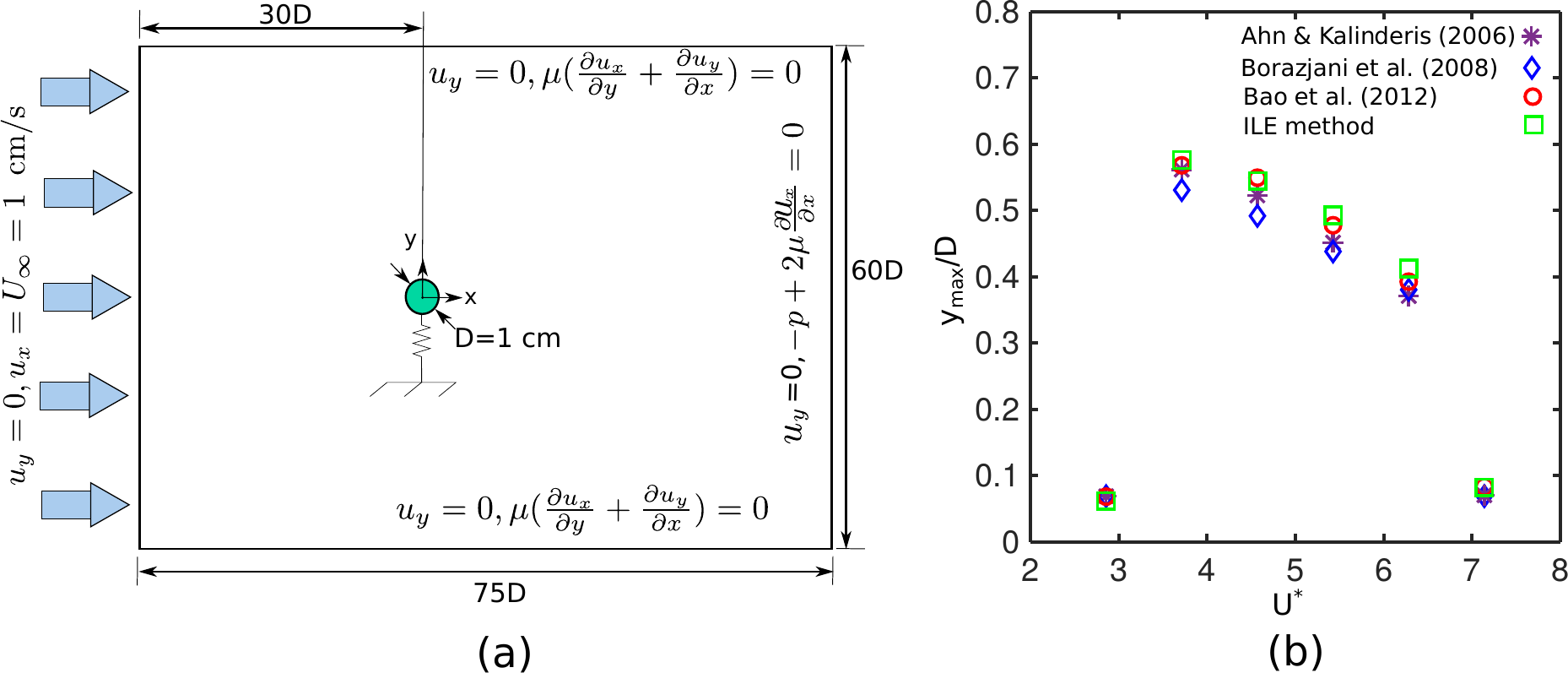}
\caption{(a) Schematic diagram of the computational domain and boundary conditions for flow around an elastically mounted rigid circular 
cylinder (Sec.~\ref{subsec:viv-1dof-cylinder}).
(b) Maximum transverse displacement of the oscillating cylinder (Sec.~\ref{subsec:viv-1dof-cylinder}) for different values of $\Ured$ 
with $m^{*}=8/\pi$, $\Re=150$, and $\zeta=0$.
         Results from the present ILE approach fall within the range of values reported in previous studies.}              
\label{fig:VIV-1DOF-cylinder-Ureds} 
\end{figure}
The problem of viscous flow past an elastically mounted two-dimensional cylinder undergoing vortex-induced vibration (VIV) 
has been widely studied both numerically and experimentally because of 
its broad range of engineering applications and its interesting vortex dynamics. 
This problem has also seen substantial use in benchmarking FSI algorithms \cite{ahn2006strongly,bao2012two, borazjani2008curvilinear,yang2008strongly,yang2015non,kim2018weak}.
Here we consider two separate cases in two spatial dimensions, one using a single vertical degree of freedom, and the second with two degrees of freedom (2-DOF). The governing equations for the 2-DOF cylinder motion are
\begin{align}
       \label{eq:2DOF-VIV-x} \Ms \ddot{d^{x}_{\com}}  + \Cs \dot{d^{x}_{\com}} + \Ks{d^{x}_{\com}} &= f^x, \\
       \label{eq:2DOF-VIV-y} \Ms  \ddot{d^{y}_{\com}} + \Cs \dot{d^{y}_{\com}} + \Ks{d^{y}_{\com}} &= f^y,
\end{align}
in which $d^{x}_{\com}$ and $d^{y}_{\com}$ are the horizontal and vertical displacements of the center of mass, $M_{\textrm{s}}$ is the mass
per unit length of the cylinder, $\Cs$ and $\Ks$ are the damping and stiffness constants characterizing the spring, and $f^x$ and $f^y$ are
the instantaneous drag and lift forces, respectively.
To facilitate comparisons to previous work, we define
the non-dimensional horizontal and vertical displacements of the cylinder center in the streamwise and transverse directions, respectively, 
as $\overline{d^{x}_{\com}}={d^{x}_{\com}}/{D}$ 
and $\overline{d^{y}_{\com}}={d^{y}_{\com}}/{D}$, in which $D$ is the diameter of the cylinder.
Taking $U_{\infty}$ as the free stream flow velocity, 
the mass ratio and reduced velocity are respectively defined as $\Mred=\rhos/\rhof$ and $\Ured=U_{\infty}/(f_{\textrm{n}} D)$, in which
$f_{\textrm{n}}=\sqrt{\Ks/\Ms}/(2\pi)$ is the natural frequency of the structure. The damping ratio is $\zeta = \Cs/(2\sqrt{\Ks \Ms})$.

%%%%%%%%%%%%%%%%%%%%%%%%%%%%%%%%%%%%%%%%%%%%%%%%%%%%%%%%%%%%%%%%%%%%%%%%%%%%%%%%%%%%%%%%%%%%%%%%%%%%%%%%%%%%%%%%%%%%%%%%%%%
 \begin{figure}[!!t!!]
		\centering
			\includegraphics[width=0.6\textwidth]{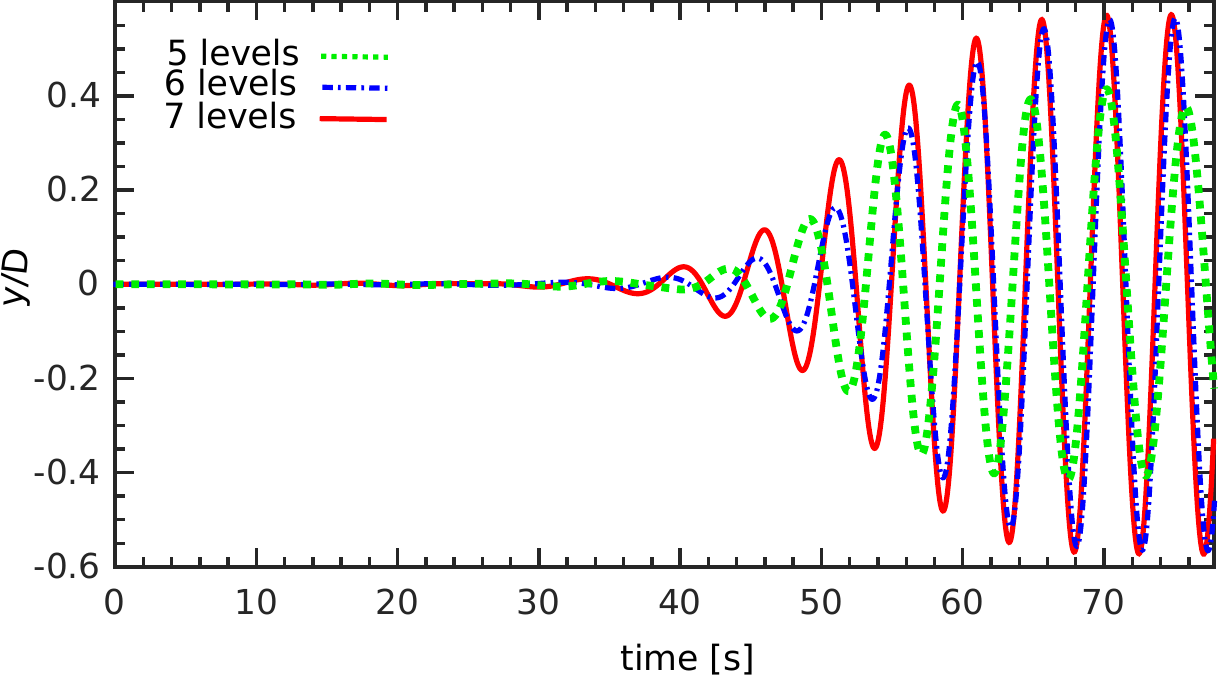}
\caption{Time history of the transverse displacement of the oscillating cylinder (Sec.~\ref{subsec:viv-1dof-cylinder}) under grid refinement. Simulation parameters
          include $\Ured=4$, $m^{*}=8/\pi$, $\Re=150$, and $\zeta=0$.}
\label{fig:VIV-1DOF-cylinder-convergence} 
\end{figure}
\subsubsection{1-DOF transverse oscillation}
\label{subsec:viv-1dof-cylinder}
We first consider
the benchmark problem of a circular cylinder undergoing VIV with a single vertical degree of freedom.
We are interested in capturing the well characterized vortex ``lock-in'' phenomenon observed in previous studies \cite{blackburn1993two, ahn2006strongly,borazjani2008curvilinear, bao2012two}.
Within the lock-in regime, the vortex shedding frequency is close to the natural frequency of the structure, which results
in large amplitude vibrations. Physical parameters are selected to match benchmark results in the literature.
A schematic of the simulation setup is shown in Fig.~\ref{fig:VIV-1DOF-cylinder-Ureds}(a).
The computational domain is $\Omega =[-30 \, \textrm{cm},45 \, \textrm{cm}]\times[-30 \, \textrm{cm},30 \, \textrm{cm}]$,
a rectangle of size $L_x \times L_y = 75 \, \textrm{cm} \times 60 \, \textrm{cm}$.
The cylinder has diameter $D=1$~cm, is initially at rest and centered at the origin.
A uniform inflow velocity $\vec{U}=(U_{\infty} = 1 \, \textrm{cm}\cdot\textrm{s}^{-1} ,0 \, \textrm{cm}\cdot\textrm{s}^{-1})$ is imposed on the left boundary ($x=-30$~cm), 
and zero normal traction and tangential velocity outflow conditions are imposed at the right boundary ($x=45$~cm).
Along the bottom ($y=-30$~cm) and top ($y=30$~cm) boundaries, zero normal velocity and tangential traction are imposed.
The domain is discretized using $N=6$ nested grid levels, with coarse grid spacing 
$h_{\textrm{coarsest}}=L_y/64=0.9375 \, \textrm{cm}$ and refinement ratio $r = 2$ between levels,
leading to $h_{\textrm{finest}} \approx 0.029 \, \textrm{cm}$.
With $\Mfac=2$, this results in approximately 54 linear elements around the perimeter of the disk.
The time step size is $\Delta t= (0.1 \, \textrm{s}\cdot\textrm{cm}^{-1}) \, h_{\textrm{finest}}$ with a penalty spring constant $\kappa=(0.00125 \, \textrm{g}\cdot\textrm{cm}^{-2})/\dt^2$.
The Reynolds number $\Re={\rhof U_{\infty} D}/{\muf}$ is fixed at 150, the damping is set to zero ($\zeta=0$), and the mass ratio is $m^{*} = 8/\pi$.

First, the effect of the reduced velocity on the maximum displacement of the cylinder is studied by systematically varying $\Ured$ within the range $ 3 \leq \Ured \leq 8$.
Vortex shedding occurs in all cases. Fig.~\ref{fig:VIV-1DOF-cylinder-Ureds}(b) shows the maximum displacement with respect to $\Ured$.
These results demonstrate that for $\Ured \in [4,7]$ there is a large increase in the vibration amplitude.
Fig.~\ref{fig:VIV-1DOF-cylinder-Ureds}(b) also compares results obtained by our method to previous numerical studies,
including
a geometrically conservative finite volume ALE method \cite{ahn2006strongly}, a curvilinear immersed boundary method 
\cite{borazjani2008curvilinear},
and a finite element based ALE approach \cite{bao2012two}.
Excellent agreement is obtained over the full range of $\Ured$ values considered here.

In addition, we perform a grid refinement study using the reduced velocity with the largest maximum displacement (i.e.~$\Ured=4$)
while fixing the previous values of all other parameters. To achieve this in our AMR framework, we vary the number of refinement levels $N$ between 5 and 7.
Fig.~\ref{fig:VIV-1DOF-cylinder-convergence} shows the time-history of the maximum displacement for $N=5$, 6, and 7 levels of refinement.
The displacement values obtained for $N = 6$ and 7 closely match each other, whereas the coarser case, using $N=5$ levels of refinement, under-predicts the maximum displacement 
in the vortex shedding region. 
Notice that these results indicate that using $N=6$ levels 
of refinement provides essentially grid-converged results for this benchmark problem.

%%%%%%%%%%%%%%%%%%%%%%%%%%%%%%%%%%%%%%%%%%%%%%%%%%%%%%%%%%%%%%%%%%%%%%%%%%%%%%%%%%%%%%%%%%%%%%%%%%%%%%%%%%%%%%%%%%%%%%%%%%%
\subsubsection{2-DOF oscillation}
\label{subsec:viv-2dof-cylinder}

The 2-DOF oscillation is studied using the same physical parameters as Blackburn and Karniadakis \cite{blackburn1993two}, who used 
a spectral element approach.
The size of the computational domain, the cylinder diameter and initial position, and the physical boundary conditions are all the same as the 1-DOF example in Sec.~\ref{subsec:viv-1dof-cylinder}.
The Reynolds number is $\Re=200$, the reduced velocity is $U^{*} = 0.5$, the damping ratio is $\zeta = 0.01$, and the mass ratio is $m^{*}=4/\pi$.

\begin{figure}[!!b!!]
		\centering
			\includegraphics[width=0.8\textwidth]{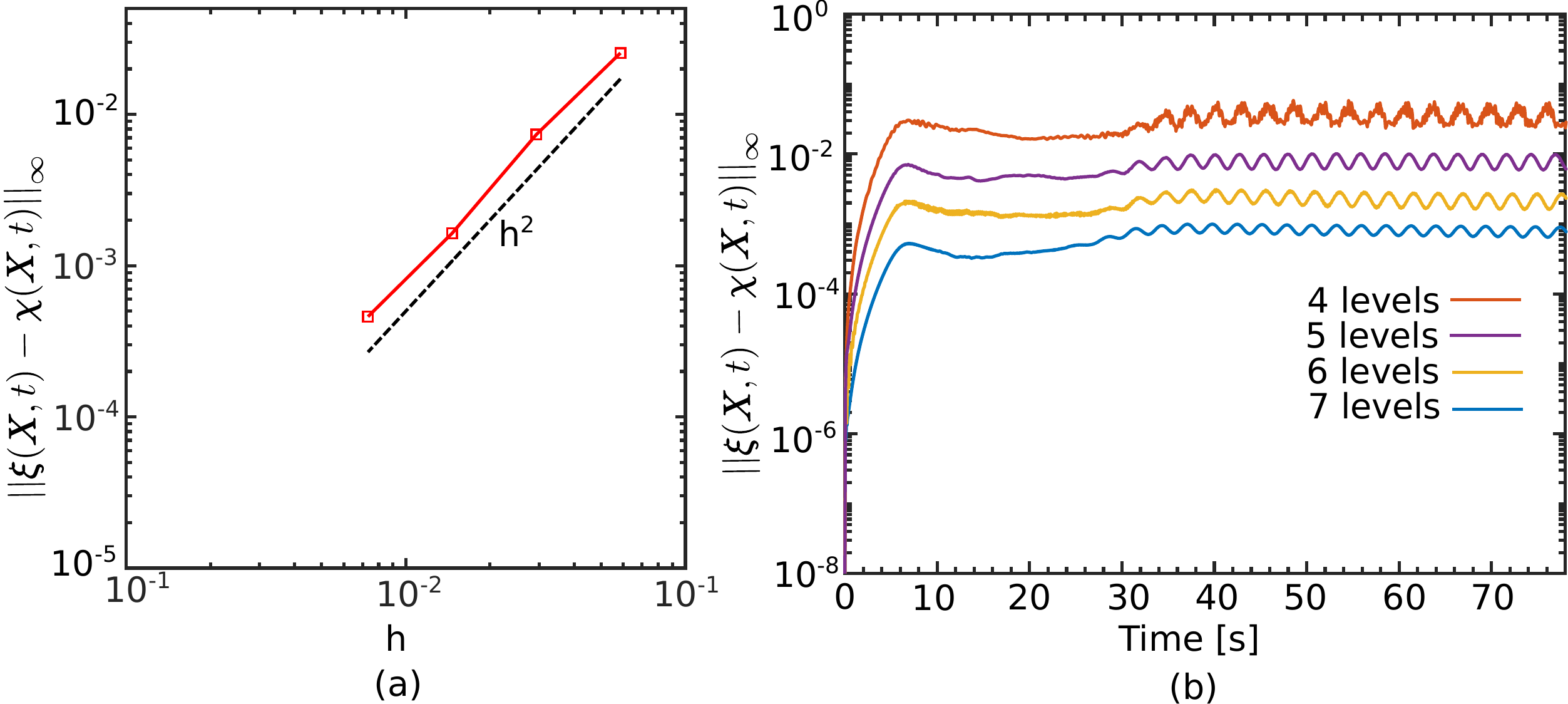}
\caption{Spatial convergence of the difference in the displacement between the two Lagrangian representations (Sec.~\ref{subsec:viv-2dof-cylinder}). (a) $\Linf$ norm of the difference between the positions of the two representations of the fluid-structure interface after 
the onset of the vortex shedding at $t=55$~s.
(b) Time history of the $\Linf$ difference, indicating a consistent behavior for all discretizations throughout the simulation.} 
\label{fig:Lmax_disp_VIV_2DOF_conv} 
\end{figure}
\begin{figure}[!!b!!]
		\centering
			\includegraphics[width=0.36\textwidth]{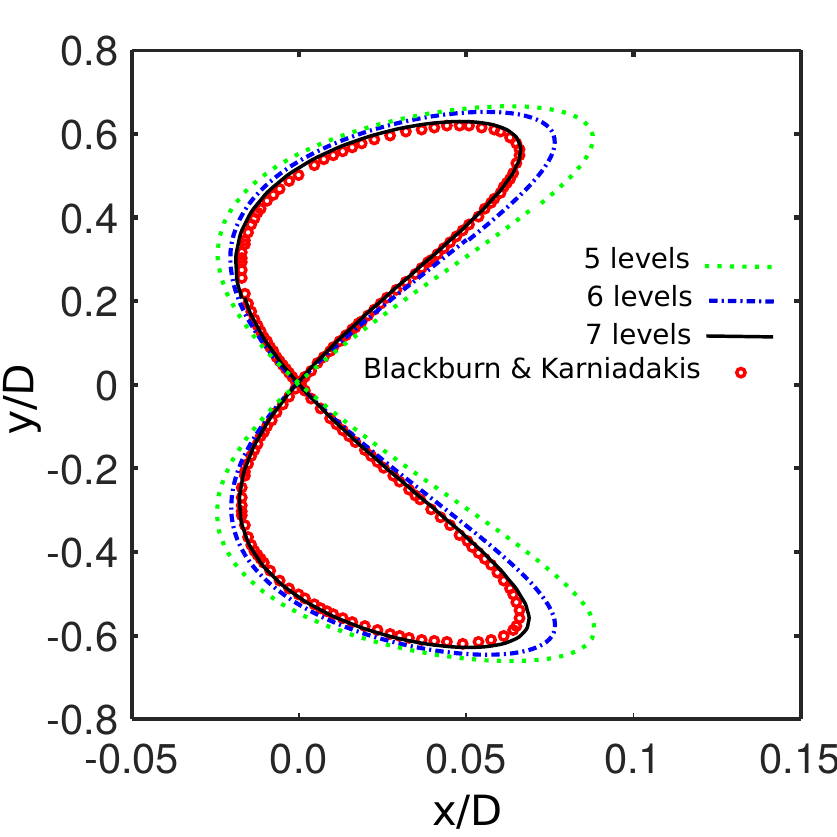}
\caption{Centerline displacement response of the 2-DOF elastically mounted rigid cylinder (Sec.~\ref{subsec:viv-2dof-cylinder}) under grid refinement. 
There is excellent agreement between the ILE result on the finest grid and 
the results of Blackburn and Karniadakis \cite{blackburn1993two} using a
high order spectral element method. Simulation parameters include $\Ured=5$, $m^{*}=4/\pi$, $\Re=200$, and $\zeta=0.01$.}
\label{fig:VIV-2DOF-cylinder-convergence}  
\end{figure}
We use this example to investigate the order of spatial convergence of the difference in the displacement
between the two Lagrangian 
representations (i.e.~the fluid-structure interface representation that moves with the fluid, and the surface of the volume mesh).
As discussed in Sec.~\ref{subsec:ILE_approach}, the two representations are tethered together by
spring-like forces that penalize relative motion between them. 
The Cartesian grid spacing on the coarsest level for all four cases is $h_{\textrm{coarsest}}=60/64$~cm. 
We consider four sets of discretizations, ranging from a very coarse composite grid with $N=4$ levels of refinement to a reasonably fine grid with $N=7$
levels of refinement, using refinement ratio $r=2$ between levels.
The penalty spring constant is $\kappa=(0.125 \, \textrm{g}\cdot\textrm{s}^{-2})/h^2$, and a small constant time step size $\dt = 0.001$~s is chosen for all cases.
In our recent work \cite{kolahdouz2020immersed}, it was shown that
the present IIM algorithm yields second-order convergence in the displacement of the interface
in both $\L2$ and $\Linf$ with suitable scalings for the penalty parameter $\kappa$.
As in our previous work, for the choice of the $\kappa$ used here, pointwise second order convergence is also expected for the displacement between the positions of 
the two interface representations. 
Fig.~\ref{fig:Lmax_disp_VIV_2DOF_conv} shows $\norm{\Y(\s,t)-\X(\s,t)}_{\infty}$,
the $\Linf$ norm of the discrepancy between the Lagrangian points at the surface mesh and the corresponding points on the boundary of the volume mesh. 
Fig.~\ref{fig:Lmax_disp_VIV_2DOF_conv}(a) indicates that the method converges at second order in the maximum norm at $t=55$~s.
To investigate the change of the $\Linf$ norm of the error over time, this value is plotted on a semi-log scale in Fig.~\ref{fig:Lmax_disp_VIV_2DOF_conv}(b) for
the four discretizations.
The overall behavior of the error norm remains consistent for all discretizations throughout the entire simulation.
Note that the difference shown in $\vec{\xi}(\vec{X},t)-\vec{\chi}(\vec{X},t)$ is proportional to the tethering penalty force that connects the two representations.
Because of the periodic ``figure-of-eight''  nature of the cylinder's dynamics, a undulatory pattern in the loading force is to be expected over time.
This pattern is expected to be present in $\|{\vec{\xi}(\s,t)-\vec{\chi}(\s,t)}\|$ difference as well due to its proportionality with the force.
Specifically if we wish to achieve $\|\vec{\xi}(\s,t)-\vec{\chi}(\s,t)\| = O(h^2)$, it is necessary 
that the penalty parameter $\kappa$ also satisfies $\kappa = O(1/h^2)$, 
so that an applied penalty force of the form $\F = \kappa (\vec{\xi}(\s,t)-\vec{\chi}(\s,t))$ satisfies $\|\F\| = O(1)$ under grid refinement \cite{kolahdouz2020immersed}.  
 \begin{table}[!!!t!!!]
	\centering		
	\caption{Dimensionless origin of oscillation ($x_\textrm{c}/D$) and the Strouhal number ($\St$) for the 2-DOF elastically mounted rigid cylinder (Sec.~\ref{subsec:viv-2dof-cylinder}).
	Simulation parameters include $m^{*}=4/\pi$, 
	$\Ured=5$, $\zeta=0.01$,  and $\Re=200$, which generate vortex-induced vibrations (VIV).}
	\label{table:2DOF_xc_st_heavy_comparison}
\begin{tabular}{l*{3}{c}r}
\hline
 & $x_{\textrm{c}}/D$ & $\St$  \\
\hline
Blackburn and Karniadakis \cite{blackburn1993two}          & 0.620 & -  \\
Yang et al.~\cite{yang2008strongly}      & 0.639 & -  \\
Yang \& Stern \cite{yang2015non}      & - & 0.187  \\
Kim et al.~\cite{kim2018weak}     & 0.622 & 0.186  \\
Qin et al.~\cite{jianhuaefficient}     & 0.626 & 0.187  \\
ILE method          &  0.619 & 0.187 
\end{tabular}
\end{table}
\begin{figure}[!!t!!]
		\centering
			\includegraphics[width=0.54\textwidth]{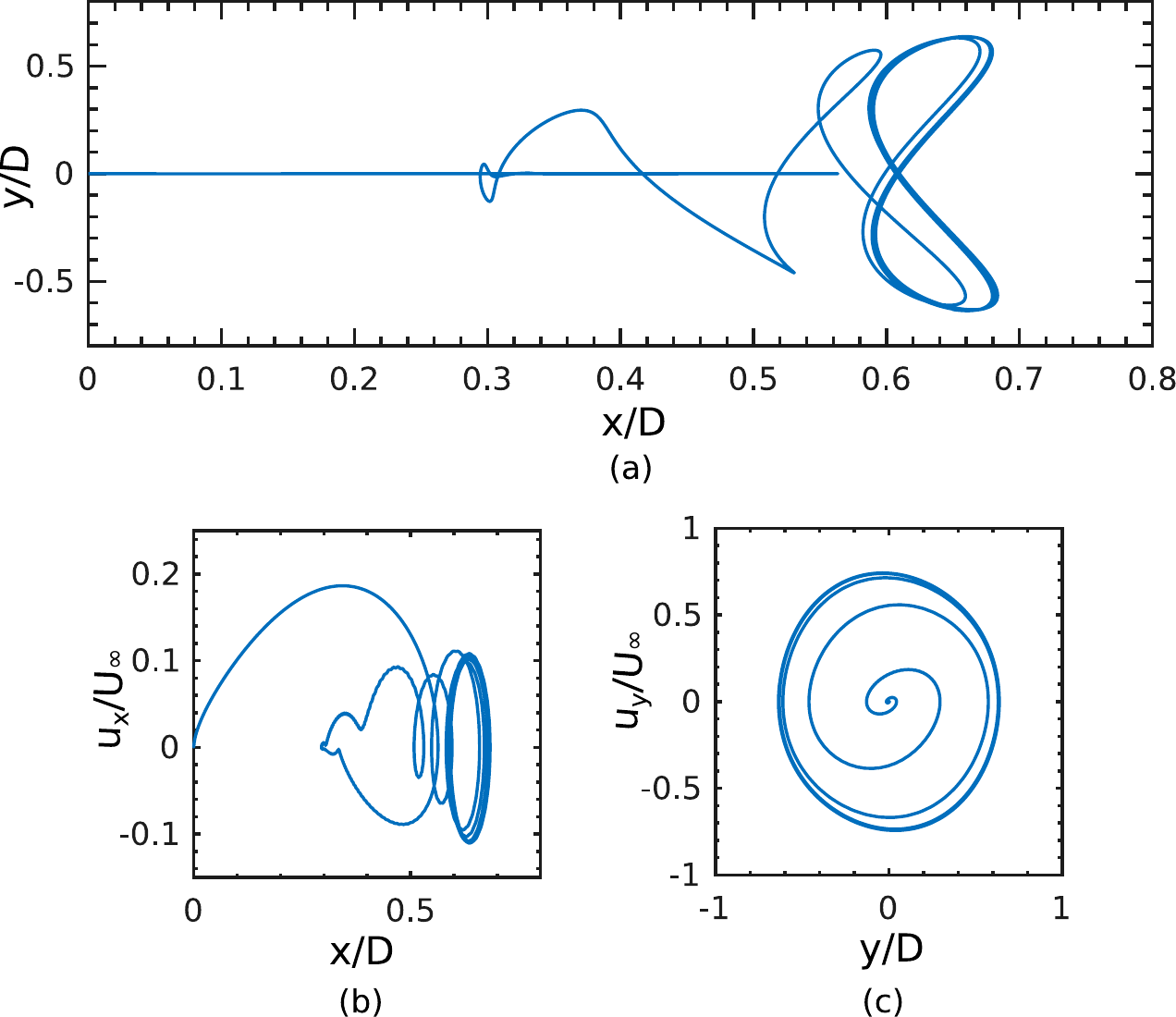}
			\caption{Phase plots of the center of mass displacement and velocity responses
			for an elastically mounted cylinder with mass ratio of  $m^{*}=4/\pi$ (Sec.~\ref{subsec:viv-2dof-cylinder}). Other simulation parameters include  $\Ured=5$, $\zeta=0.01$, and $\Re=200$. }
\label{fig:VIV-2DOF-heavy-cylinder}  
\end{figure}

As another verification, the centerline displacement response of the cylinder is compared to 
the numerical results of Blackburn and Karniadakis \cite{blackburn1993two}.
It is expected that the vortex shedding of the oscillating cylinder will lead to a periodic ``figure-of-eight''  behavior. 
This is shown in Fig.~\ref{fig:VIV-2DOF-cylinder-convergence} 
for three locally refined Cartesian grids with $N=$ 5, 6, and 7 levels along with the solution from 
Blackburn and Karniadakis \cite{blackburn1993two}.
The trajectory clearly converges under grid refinement. Moreover, the trajectory of the cylinder
for the finest discretization agrees extremely well with the spectral element solution \cite{blackburn1993two}.
The origin of oscillations $x_{\textrm{c}}$, which is defined as the intersection point in the periodic figure-of-eight trajectory of the cylinder's
center of mass, and the Strouhal number, which is calculated as $\St=f D/U_{\infty}$ with $f$ representing the oscillation frequency, are
reported in Table~\ref{table:2DOF_xc_st_heavy_comparison} for the finest ILE computation ($N=7$) along with the results 
of Blackburn and Karniadakis and several other studies.
The center of oscillations agrees very well with the original work of Blackburn and Karniadakis \cite{blackburn1993two}. Moreover, 
the Strouhal number matches the values reported
from studies by Yang et al.~\cite{yang2008strongly}, Yang and Stern \cite{yang2015non}, Kim et al.~\cite{kim2018weak}, and Qin et al.~\cite{jianhuaefficient}.
The centerline trajectory and the dimensionless displacement-velocity phases of $x/D-u_x/U_{\infty}$ and $y/D-u_y/U_{\infty}$ for this case are shown in
Fig.~\ref{fig:VIV-2DOF-heavy-cylinder}.
The phase response obtained from our solution agrees well with previous observations reported 
by Yang and Stern \cite{yang2012simple} and Liu and Hu \cite{liu2018block}.

\begin{figure}[t!!]
		\centering
			\includegraphics[width=0.53\textwidth]{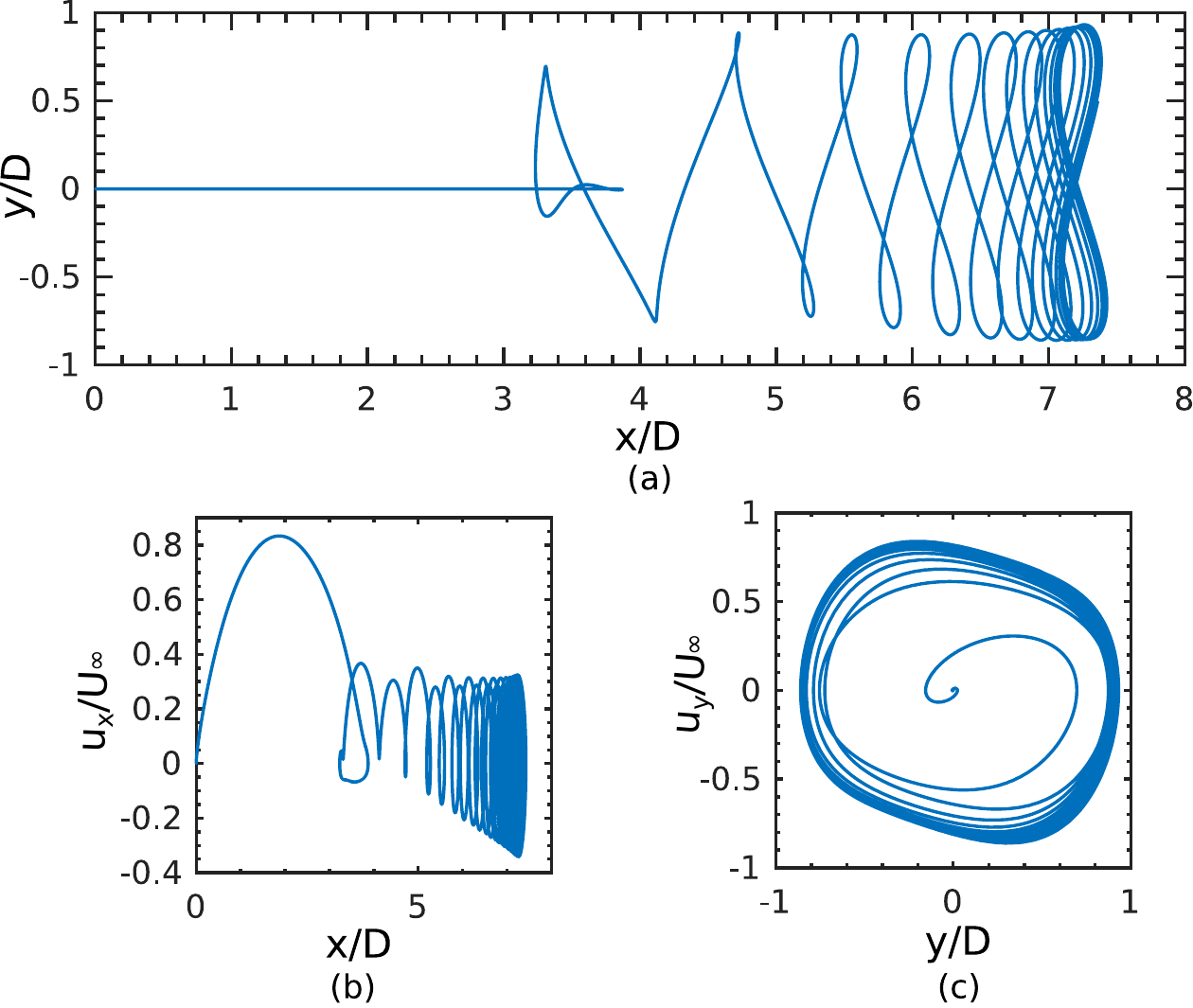}
			\caption{Phase plots of the center of mass displacement and velocity responses
			for an elastically-mounted cylinder (Sec.~\ref{subsec:viv-2dof-cylinder-low-mass-ratio}) with a low mass ratio $m^{*}=0.4/\pi$. 
			Other simulation parameters include  $\Ured=5$, $\zeta=0.01$, and $\Re=200$. }
\label{fig:VIV-2DOF-light-cylinder}  
\end{figure}
\begin{figure}[b!!!]
		\centering
			\includegraphics[width=0.8\textwidth]{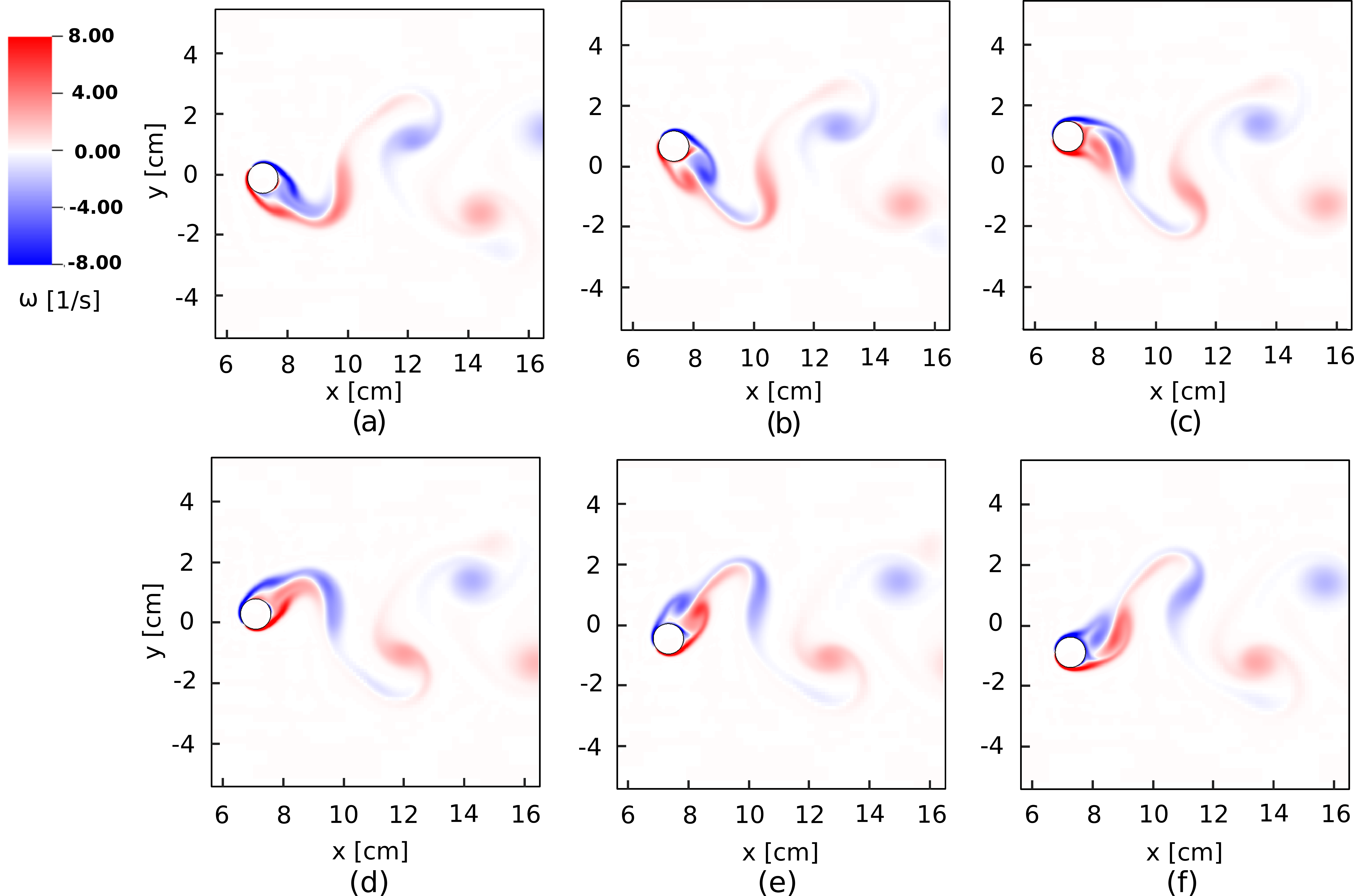}
\caption{Vorticity field of the 2-DOF oscillating cylinder (Sec.~\ref{subsec:viv-2dof-cylinder-low-mass-ratio}) with low mass ratio $m^{*}=0.4/\pi$ at times (a) $t=178.82$~s, (b) $t=179.82$~s, (c) $t=180.82$~s, (d) $t=181.82$~s, (e) $t=182.82$~s, and (f) $t=183.82$~s.} 
\label{fig:VIV-cylinder-2DOF-light-vorticity} 
\end{figure}  
%%%%%%%%%%%%%%%%%%%%%%%%%%%%%%%%%%%%%%%%%%%%%%%%%%%%%%%%%%%%%%

\subsubsection{2-DOF oscillation with low mass ratio}
\label{subsec:viv-2dof-cylinder-low-mass-ratio}
We now briefly consider significantly lower mass ratios than the one studied above.
We first consider a case with density ratio of $m^{*}=0.4/\pi$, which is 10 times smaller than the first example.
The remainder of the simulation
parameters are identical to the finest case in the previous example (Sec.~\ref{subsec:viv-2dof-cylinder}). Fig.~\ref{fig:VIV-2DOF-light-cylinder}(a) shows that the cylinder
travels more than 7 diameters downstream before undergoing
the same periodic figure-of-eight motion.
Fig.~\ref{fig:VIV-2DOF-light-cylinder} also reveals that the oscillations for this light cylinder occur
at a higher amplitude and with larger horizontal and, to a lesser extent, vertical velocity magnitudes.    
Fig.~\ref{fig:VIV-cylinder-2DOF-light-vorticity} shows the 
instantaneous vorticity contours around this light structure at six time points.
Vortex shedding is observed, with two vortices
shed during one cycle of the quasi-periodic oscillation. The wake footprint and the
associated vorticity patterns
also indicate that the cylinder is undergoing periodic motion.
Notice that previous work has reported severe instabilities
in computing such cases using both weak and strong coupling approaches \cite{kim2018weak}. 
The present method appears to remain stable even for extremely light structures. 
Table~\ref{table:2DOF_xc_st_light_comparison} reports
the center of oscillations along the $x$-axis 
for a wide range of mass ratios along with the results of Kim et al.~\cite{kim2018weak}, the only other study that we are aware of to
also consider such low mass ratios.
To the authors' knowledge, results obtained using mass ratios smaller than 0.3 have not been previously reported for this benchmark case.
The present method yields a slightly larger distance than the work of Kim et al.~\cite{kim2018weak} for the origin of the oscillations. 
In the work of Kim et al.~\cite{kim2018weak}, the mass ratio of $m^{*}=0.3$ was reported as the lowest ratio for a stable solution using
a strong coupling approach.
The present method is able to predict the dynamics at substantially smaller mass ratios (nearly two orders of magnitudes smaller) as shown for $m^{*}=0.005$ in Table~\ref{table:2DOF_xc_st_light_comparison}.
Although stable results were obtained for lower density ratios, the accuracy of the dynamics requires
further validation that can be performed once additional experimental data are available.
Note that to accommodate the extended horizontal movement of the cylinder for $m^{*} < 0.1$, the computational domain for these cases was changed to $\Omega = [-120 \, \textrm{cm}, 180 \, \textrm{cm}]\times[-120 \, \textrm{cm},120 \, \textrm{cm}]$,
a rectangle of size $L_x \times L_y  = 300 \, \textrm{cm} \times  240 \, \textrm{cm}$.
To keep the grid spacing the same as before, $N=8$ nested grid levels were used, and the Cartesian grid spacing on the coarsest level was set to
$h_{\textrm{coarsest}}=L_y/64$.
\begin{table}[t!!]
	\centering	
	\caption{Dimensionless origin of oscillation ($x_\textrm{c}/D$) for low and very low mass ratios of the 2-DOF elastically mounted rigid cylinder 
	(Sec.~\ref{subsec:viv-2dof-cylinder-low-mass-ratio}).}
	\label{table:2DOF_xc_st_light_comparison}	
\begin{tabular}{l*{7}{c}r}
\hline
 & $m^{*} = 0.3$ & $m^{*} = 0.2$ & $m^{*} = 0.1$ & $m^{*} = 0.05$ & $m^{*} = 0.01$ & $m^{*} = 0.005$   \\
\hline
Kim et al.~\cite{kim2018weak}     & 2.27 & -  & -  & - & - & - \\
ILE method            & 2.83  & 4.41 & 9.39 & 19.01 & $91.50$ & $160.10$ 
\end{tabular}
\end{table}

%%%%%%%%%%%%%%%%%%%%%%%%%%%%%%%%%%%%%%%%%%%%%%%%%%%%%%%%%%%%%%%%%%%%%%%%%%%%%%%%%%%%%%%%%%%%%%%%%%%%%%%%%%%%%%%%%%%%%%%%%%%

\subsubsection{Stability and sensitivity studies for $\kappa$}
\label{subsec:sensitivity-kappa}

\begin{figure}[b!!!]
		\centering
			\includegraphics[width=\textwidth]{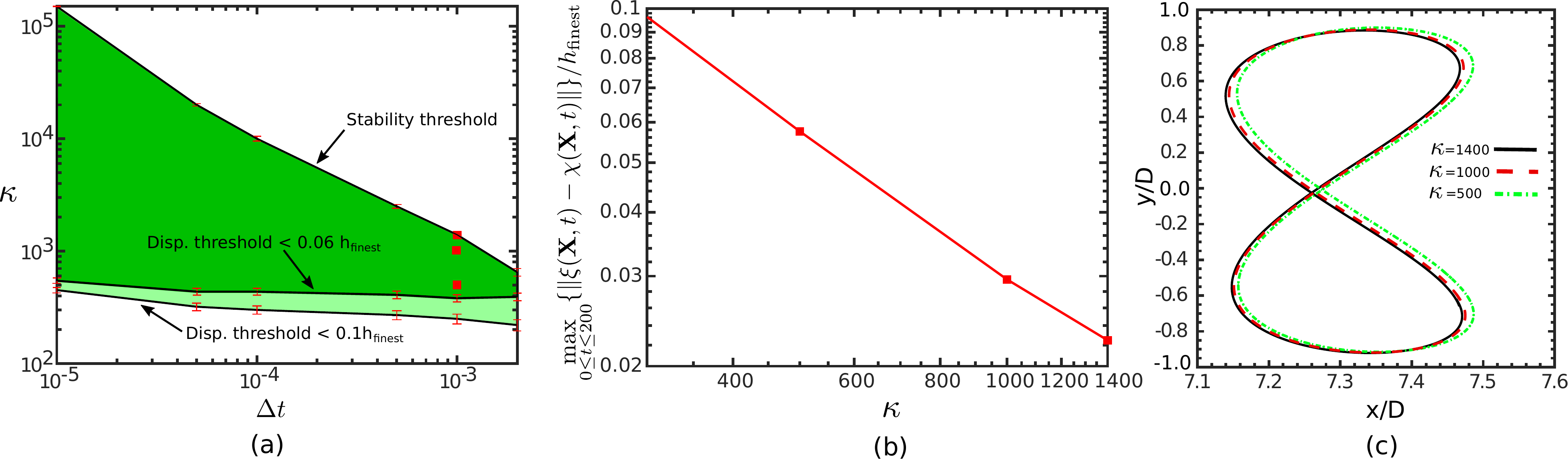}
\caption{Impact of the penalty parameter $\kappa$ on the computed dynamics for the 2-DOF elastically mounted light cylinder (Sec.~\ref{subsec:sensitivity-kappa}).
Other simulation parameters in these computations are $m^{*}=0.4/\pi$, $\Ured=5$, $\zeta=0.01$, and $\Re=200$
(a) Accuracy and stability thresholds of $\kappa$ for a range of time step sizes. 
(b) The ratio between the maximum norm of the difference in the positions of the two Lagrangian representations to the finest Cartesian grid spacing 
over the entire course of the simulation, determined using $\dt=0.001$ and $\kappa=500, 1000, \, \textrm{and} \, 1400$.
(c) Centerline displacements of the cylinder for the three values of $\kappa$ highlighted in panel (b).}
\label{fig:kappa-sensitivity} 
\end{figure}
Because we use an explicit time stepping approach to couple the fluid and solid degrees of freedom, 
if we fix the spatial and temporal discretization parameters, there will be a largest stable value of 
the penalty spring parameter $\kappa$. As mentioned at the beginning of Sec.~\ref{sec:Results}, 
in practice, we compute using values of $\kappa$ close to the stability limit, which we approximately determine using the method of bisection.
%~ As described at the beginning of Sec.~\ref{sec:Results}, for a given time step size the spring penalty parameter $\kappa$ is determined to be approximately the largest value that satisfies the CFL condition in
%~ the explicit time stepping algorithm. This value can be tuned using a simple bisection approach. One question that arises is 
%~ to what extent bisecting of the penalty parameter has to continue in order to obtain a reliably accurate solution. 
To study the sensitivity of the numerical algorithm to the penalty parameter $\kappa$ and the influence of this parameter on the accuracy of the 
computed solution, we
 reconsider the 2-DOF vortex-induced vibration of a light cylinder with density ratio of $m^{*}=0.4/\pi$ (Sec.~\ref{subsec:viv-2dof-cylinder-low-mass-ratio}). 
 %~ First, we aim to find the maximum and minimum $\kappa$ values for
  %~ which the solution 1) stays stable (stability threshold) and 2) the difference in the displacement between the two Lagrangian representations stays 
  %~ smaller than $0.1\,h_{\textrm{finest}}$ (accuracy threshold).
Except for the values of $\kappa$ and the time step size $\dt$, all other simulation parameters are the same as Sec.~\ref{subsec:viv-2dof-cylinder-low-mass-ratio}. 
We use a locally refined grid with $N=7$ levels of refinement, and the final simulation time is $t=200$~s, which is long after vortex induced vibration has been established. 
  We computationally determine minimum and maximum values of $\kappa$ that 
  1) satisfy a minimum threshold in displacement between the two Lagrangian representations, 
  and 2) remain stable. 
  Accuracy thresholds are determined to be the minimum values of $\kappa$ required to achieve tolerances of $0.1 \, h_\text{finest}$, 
  which is commonly used in IB models of rigid structures, along with a tighter tolerance of $0.06 \, h_\text{finest}$. 
  The maximum stable value of $\kappa$ determines a stability threshold.
  %~ Stable $\kappa$ values satisfying the above-mentioned stability and accuracy thresholds are experimentally determined over a range of time step sizes.
Fig.~\ref{fig:kappa-sensitivity}(a) shows 
the computationally determined accuracy and stability thresholds for this case.
%~ the effect of $\kappa$ on the numerical instability. 
%~ The Stable solutions are obtained for a large range of $\kappa$ values with the displacement between the two Lagrangian representations remaining 
%~ smaller than $0.1\,h_{\textrm{finest}}$.
Because it is not practical to determine the precise values of these thresholds across all time step sizes,
error bars indicate the small uncertainty in the lower and upper bounds.
For the smallest time step size considered here, 
the range of values of $\kappa$ that satisfy both criteria spans nearly three orders of magnitude, 
from the minimum value that achieves the minimum accuracy 
criterion to the maximum value that results in a stable computation with $0.06 \, h_{\textrm{finest}}$ displacement threshold.
%~ there is a wide interval of stable $\kappa$ values with roughly three orders of magnitude difference between the
%~ highest and lowest bounds. 
For the largest time step, there is at least about three-fold difference between the minimum and maximum values of $\kappa$.
Fig.~\ref{fig:kappa-sensitivity}(b) shows the effect of $\kappa$ on the ratio between the maximum norm of the two Lagrangian representations
and the finest grid spacing for three stable spring constants of $\kappa = 500, 1000, \, \text{and} \, 1400$ using the time step $\dt=0.001$.
Notice that the maximum displacement between the representations scales like $\kappa^{-1}$.
Fig.~\ref{fig:kappa-sensitivity}(c) shows the trajectory of the center of mass of the cylinder for the same three $\kappa$ values.
Deviations in the trajectories are small and converge as $\kappa$ increases.
It is worth noting that the maximum discrepancies for the highest and lowest choices of
$\kappa$ are less than the spacing on the finest grid level, $h_{\textrm{finest}}$.
 %~ Except for the smallest spring constant of $\kappa=300$, the deviation in trajectories is negligible.

%%%%%%%%%%%%%%%%%%%%%%%%%%%%%%%%%%%%%%%%%%%%%%%%%%%%%%%%%%%%%%%%%%%%%%%%%%%%%%%%%%%%%%%%%%%%%%%%%%%%%%%%%%%%%%%%%%%%%%%%%%%
\subsection{Galloping rectangular structure}
\label{subsec:galloping}

This example uses a rectangular plate undergoing galloping motion to test the accuracy of the method for models involving only rotational degrees of freedom.
Flow-induced rotational galloping oscillations occur
in many areas of structural \cite{mannini2016aeroelastic}, aeronautical \cite{alonso2010galloping}, and mechanical \cite{yang2013comparative} engineering applications. 
This problem has also been widely used as a benchmark problem to validate numerical 
algorithms \cite{robertson2003numerical,yang2012simple,yang2015non}.
In this section, only a single rotational degree of freedom is used, and
the translational heave (horizontal) and surge (longitudinal) motions are eliminated.
The governing equation for the mass-spring-damper model with one rotational degree of freedom is
\begin{equation}
       \Is^{\theta} \ddot{\theta}  + \Cs^{\theta} \dot{\theta} + \Ks^{\theta}\theta = M^{\theta}, \\
       \label{eq:VIV-galloping} 
\end{equation}
in which $\theta$ is the rotational angle of the body, $\Is^{\theta}$ is the rotational moment of inertia, $\Cs^{\theta}$ is the torsional damping constant,
$K_{\textrm{s}}^{\theta}$ is the torsional spring constant, and $M^{\theta}$ is the moment acting on the rigid structure from 
exterior fluid forces.
In our simulations, we consider a rectangular structure with a width-to-thickness ratio of $\Lambda=L/H=4$.
\begin{figure}[b!!!]
		\centering
			\includegraphics[width=0.95\textwidth]{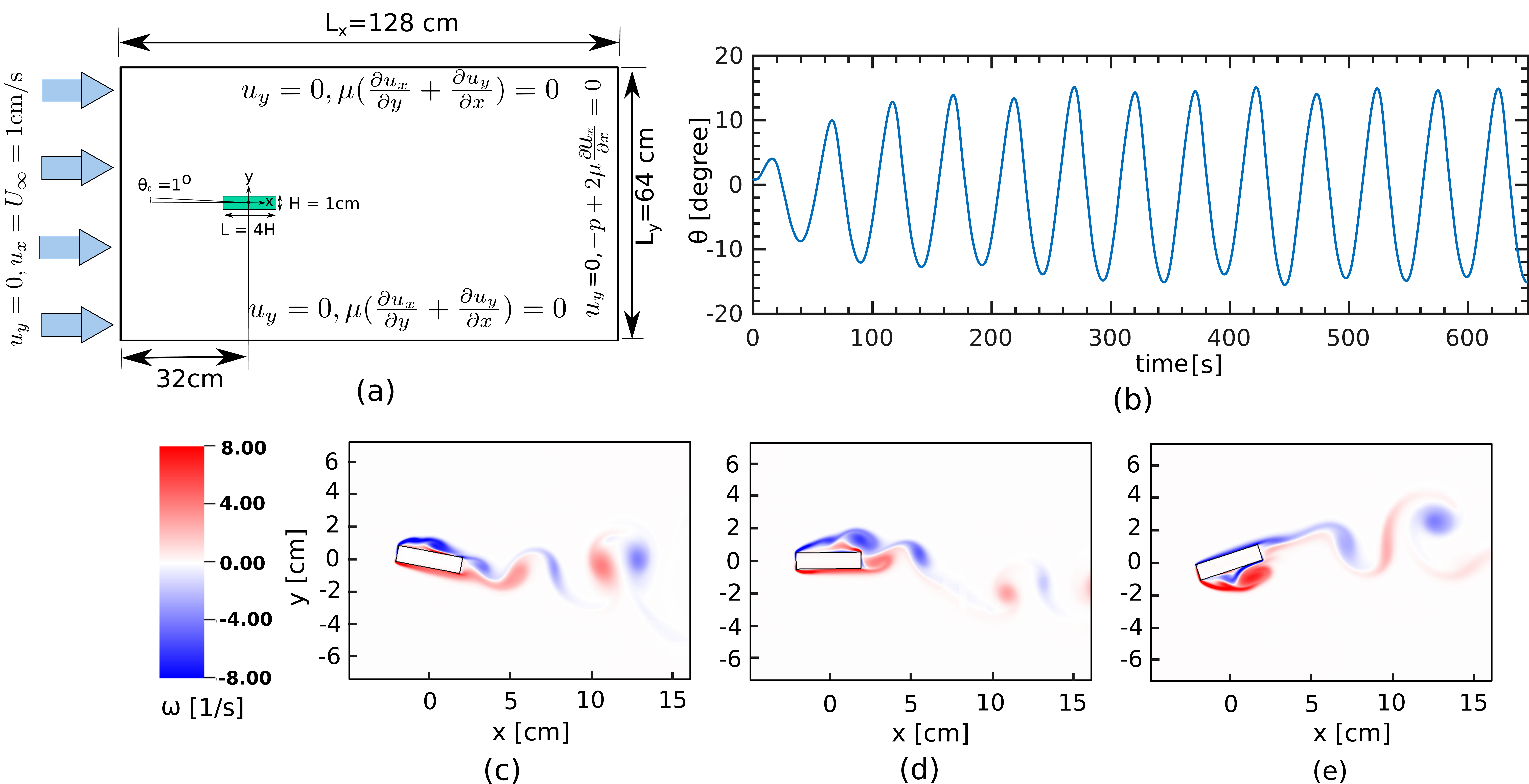}
\caption{(a) Schematic diagram of the computational domain and boundary conditions for the galloping rectangle (Sec.~\ref{subsec:galloping}) with $\zeta_{\textrm{s}}^{\theta}=0.25$. (b) Pitch angle ($\theta$) as a function of time. (c)--(e) Vorticity
 field at times $t=293.9$~s, $t=305.9$~s, and $t=321.9$~s, respectively.} 
\label{fig:galloping-angles-damping} 
\end{figure}
\begin{figure}[b!!]
		\centering
			\includegraphics[width=0.9\textwidth]{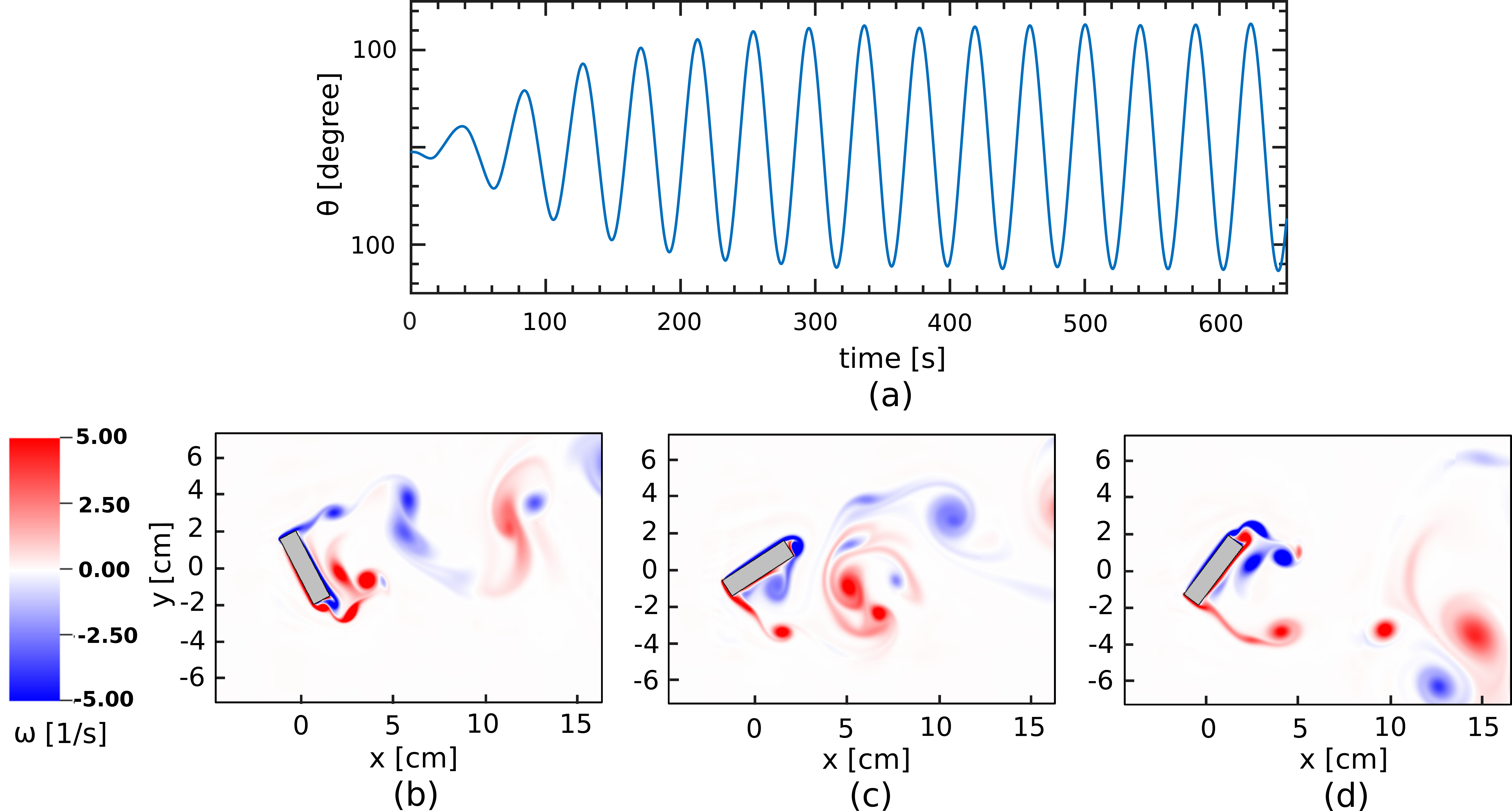}
\caption{Representative results for the galloping rectangle with zero damping ($\zeta_{\textrm{s}}^{\theta}=0$). (a) Pitch angle ($\theta$) as a function of time.
(b)--(d) Vorticity field at times $t=377.9$~s, $t=385.9$~s, and $t=397.8$~s, respectively. 
Here both the bulk region (shown in gray) and the surface mesh of the structure (shown in black) are illustrated
to show the effectiveness of the coupling approach and to confirm that the two representations move together.} 
\label{fig:galloping-angles-zero-damping} 
\end{figure}
To enable comparisons with prior studies, we define the non-dimensional moment of inertia by $\Is^{*}=\Is^{\theta}/(\rhos H^4)$, 
the non-dimensional damping ratio by $\zeta_{\textrm{s}}^{\theta}=\Cs^{\theta} /(2\sqrt{\Ks^{\theta}\Is^{\theta}})$, 
and the reduced velocity by $U^{*}=U_{\infty}/(f^{\theta}H)$.
In the latter formula, $U_{\infty}$ is the free stream velocity and ${f_{\textrm{s}}}^{\theta}=\sqrt{\Ks^{\theta}/\Is^{\theta}}/{2\pi}$ is the natural frequency of the body. 
Following the work of Robertson et al.~\cite{robertson2003numerical}, the 
non-dimensional parameters are taken to be $I_{\textrm{s}}^{*}=400$, $U^{*}=40$, and $\zeta_{\textrm{s}}^{\theta}=0.25$.
A schematic of the computational setup and the boundary conditions is given in Fig.~\ref{fig:galloping-angles-damping}(a).
The rectangular structure is centered at the origin with an initial zero angular velocity and a non-zero angle of $\theta_0=1^{\circ}$.
The computational domain is $\Omega = [-32 \, \textrm{cm},96 \, \textrm{cm}]\times[-32 \, \textrm{cm},32 \, \textrm{cm}]$,
a rectangle of size $L_x \times L_y = 128 \, \textrm{cm} \times 64 \, \textrm{cm}$.
The domain is discretized using $N=4$ nested grid levels, with coarse grid spacing $h_{\textrm{coarsest}}=L_y/32 = 2.0 \, \textrm{cm}$ and 
refinement ratio $r = 4$ between levels, leading to $h_{\textrm{finest}} = 0.0625 \, \textrm{cm}$.
A uniform inflow velocity of $\vec{U}=(U_{\infty}=1\,\textrm{cm}\cdot\textrm{s}^{-1},0\,\textrm{cm}\cdot\textrm{s}^{-1})$ is imposed on the left boundary ($x=-32$~cm).
Using the free stream velocity, the Reynolds number $\Re=\rhof U_{\infty} H/\muf$
is set to 250.
A penalty spring constant of $\kappa=400$~$\textrm{g}\cdot\textrm{cm}^{-2} \cdot \textrm{s}^{-2}$ is used, and the time step size
 is $\Delta t=(0.02 \, \textrm{s}\cdot\textrm{cm}^{-1}) \, h_{\textrm{finest}}$.
Zero normal traction and tangential velocity are imposed at the right boundary ($x=96$~cm).
Along the bottom ($y=-32$~cm) and top ($y=32$~cm) boundaries, the normal velocity and tangential traction are set to zero.
Fig.~\ref{fig:galloping-angles-damping}(b) shows the time history of the pitch angle of the galloping rectangle
for the damped oscillation of the structure.
Once the vortex shedding state is established, the structure starts to undergo a
periodic rotation with well-characterized frequency and upper bound of the maximum angle.
This behavior is clearly captured by the present method.
Fig.~\ref{fig:galloping-angles-damping} panels (c)--(e) show
snapshots of the structural rotation along with the vortex structure of the flow at three different times.
\begin{table}[t!!!]
	\centering	
	\caption{Comparison of the maximum pitch angle ($\theta_{\textrm{max}}$) and galloping frequency ($f^{\theta}$) for the galloping rectangular structure (Sec.~\ref{subsec:galloping})
	with zero damping ($\zeta_{\textrm{s}}^{\theta}=0$) and non-zero damping ($\zeta_{\textrm{s}}^{\theta}=0.25$).    
	}
	\label{table:Galloping_table}	
\begin{tabular}{l*{7}{c}r}
             & \hspace{1.8cm} $\zeta_{\textrm{s}}^{\theta}=0$& & \hspace{1.9cm} $\zeta_{\textrm{s}}^{\theta}=0.25$\\
\cmidrule(lr){2-3}\cmidrule(lr){4-5}
 & $\theta_{\textrm{max}}$ & $f^{\theta}$ & $\theta_{\textrm{max}}$ & $f^{\theta}$  \\
\hline
Robertson et al. \cite{robertson2003numerical}        & - &  -  &  $15^{\circ}$ & $0.0191$   \\
Yang \& Stern  \cite{yang2012simple}      & $123^{\circ}$ & $0.0244$ & $15.7^{\circ}$   & $0.0198$  \\
Yang et al.  \cite{yang2015non}      & $125^{\circ}$  & $0.0243$ & $16.1^{\circ}$ & $0.0197$ \\
ILE method    & $124^{\circ}$ & $0.0243$ &  $15^{\circ}$ & $0.0198$  \\
\end{tabular}
\end{table}

To demonstrate the ability of the algorithm in modeling larger rotational angles, an additional case is
considered with $\zeta_{\textrm{s}}^{\theta}=0$. The initial angle is set to $\theta_0=-5^{\circ}$. All the other parameters, 
including the time step size and penalty spring constant, are the same as before.
As seen in Fig.~\ref{fig:galloping-angles-zero-damping}(a), the method generates periodic behavior 
with maximum galloping amplitude of approximately $\theta_{\textrm{max}}=124^{\circ}$.
Fig.~\ref{fig:galloping-angles-zero-damping} panels (b)--(d) show the flow patterns around the structure at three different times.

To compare the rotational response of the structure with previous work,
the maximum pitch angle and frequency for the above two cases are reported in Table~\ref{table:Galloping_table} along with values from other studies. 
The work of
Robertson et al.~\cite{robertson2003numerical} uses a body-fitted spectral element method 
in a non-inertial reference frame.
The methods used by Yang and colleagues \cite{yang2012simple,yang2015non}
are different variations of a strongly-coupled direct forcing approach with 
a field extension strategy for the pressure/velocity derivatives.
Table~\ref{table:Galloping_table} demonstrates 
excellent quantitative agreement of both vortex shedding characteristics of $\theta_{\textrm{max}}$  and the galloping frequency $f^{\theta}$
for the two cases in comparison with other numerical studies.

%%%%%%%%%%%%%%%%%%%%%%%%%%%%%%%%%%%%%%%%%%%%%%%%%%%%%%%%%%%%%%%%%%%%%%%%%%%%%%%%%%%%%%%%%%%%%%%%%%%%%%%%%%%%%%%%%%%%%%%%%%%
\subsection{Freely falling rectangular plate} 
\label{subsec:falling_plate}

We now consider a model of a falling rigid plate in a water tank that is based on the experiments by Andersen et al.~\cite{andersen2005unsteady}.
This example tests the action of the instantaneous fluid forces on an object with sharp corners that leads to extremely complex trajectories.
This benchmark case has also been investigated in the context of fluid-structure interaction
algorithms in prior studies  \cite{vanella2010direct,yang2012simple, liu2018block}.
We consider two cases from Anderson et al.~\cite{andersen2005unsteady}, one undergoing 
fluttering motion at $\Re=1147$ and the other undergoing tumbling motion at $\Re=837$.
Different modes of fluttering and tumbling motions 
are captured in both the experiments \cite{andersen2005unsteady} and the simulations.

%%%%%%%%%%%%%%%%%%%%%%%%%%%%%%%%%%%%%%%%%%%%%%%%%%%%%%%%%%%%%%%%%%%%%%%%%%%%%%%%%%%%%%%%%%%%%%%%%%%%%%%%%%%%%%%%%%%%%%%%%%%
\subsubsection{Fluttering motion} 
\label{subsubsec:fluttering_motion}

For the case of a freely fluttering plate, the plate thickness is taken to be $ H = 8.1 \, \times 10^{-2}$~cm with width-to-thickness ratio of $\Lambda = L/H = 14$.
 \begin{figure}[b!!!]
		\centering
			\includegraphics[width=0.9\textwidth]{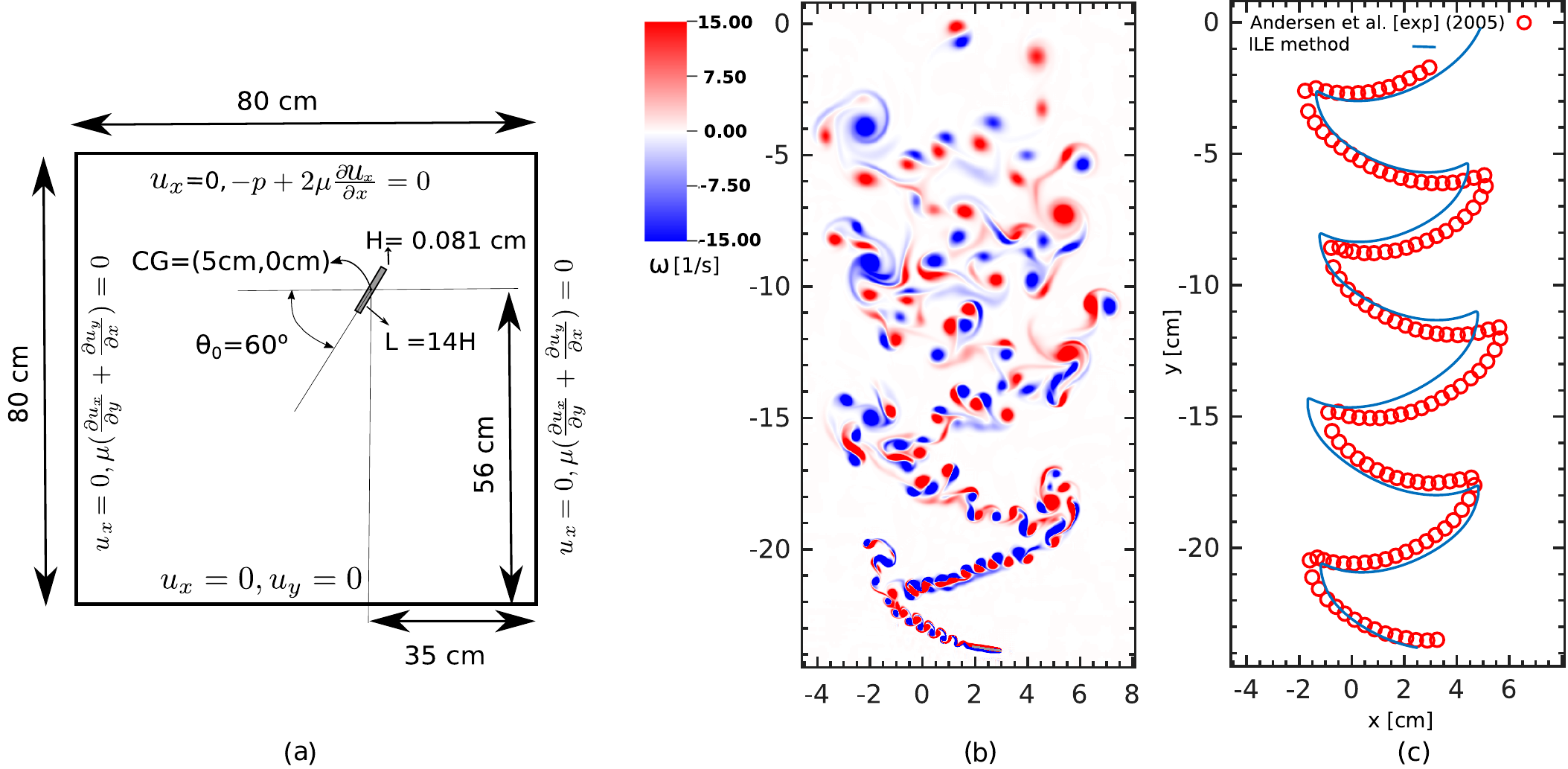}
\caption{
(a) Schematic diagram of the computational domain and boundary conditions for freely falling fluttering plate (Sec.~\ref{subsubsec:fluttering_motion}). (b) Vorticity field at $t=2.68$~s. (c) Comparison of the trajectory of the center of mass using the present ILE method (solid blue line) to the experimental data of Andersen 
et al.~\cite{andersen2005unsteady} (red circles).} 
\label{fig:trajectory-vorticity-fluttering-plate}
\end{figure} 
 \begin{figure}[b!!!]
		\centering
			\includegraphics[width=0.9\textwidth]{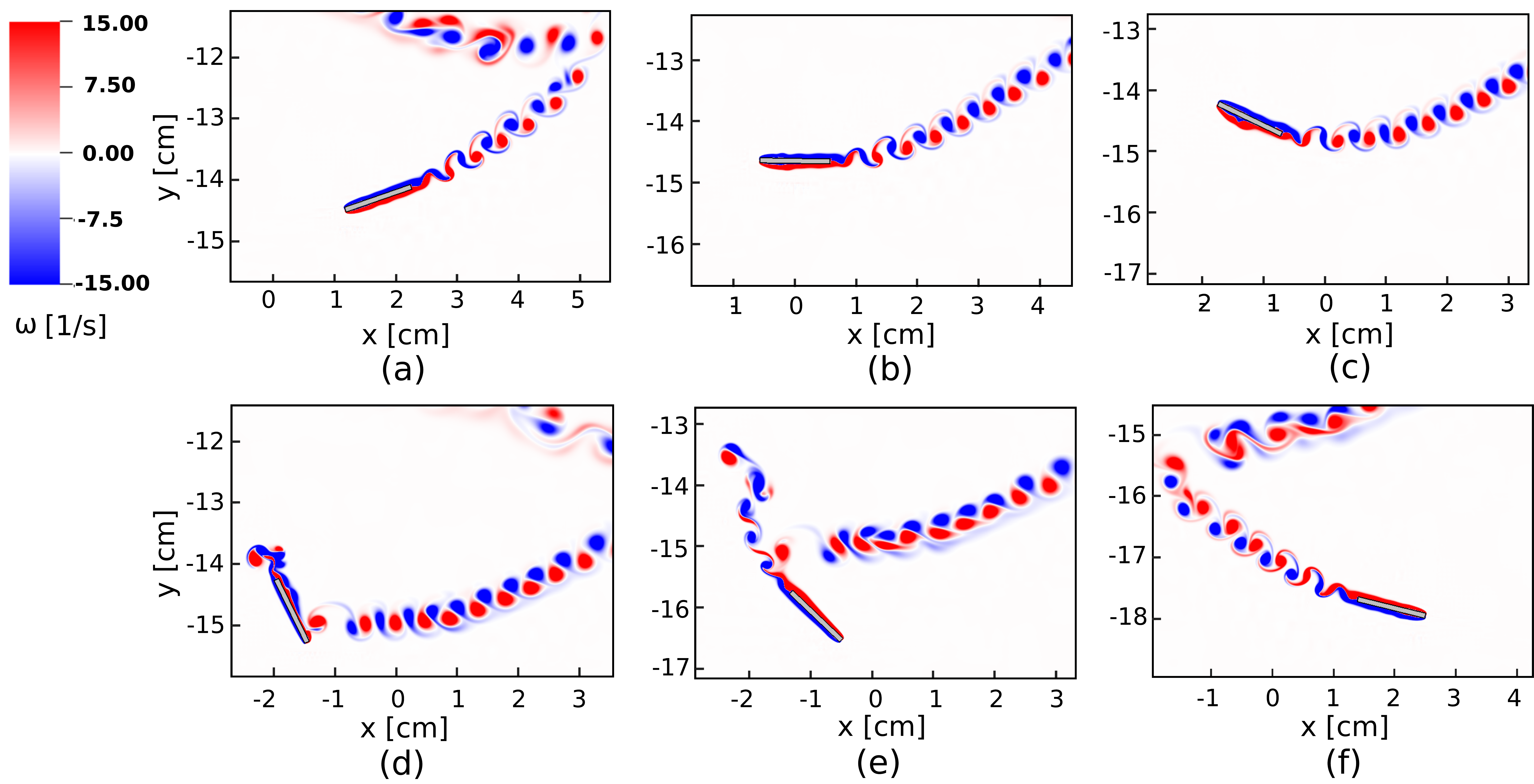}
\caption{Enlarged views of the vorticity field for the fluttering plate (see Fig.~\ref{fig:trajectory-vorticity-fluttering-plate}(b)) at times (a) $t=1.40$~s, (b) $t=1.45$~s, (c) $t=1.50$~s, (d) $t=1.59$~s, (e) $t=1.65$~s, and (f) $t=1.74$~s.} 
\label{fig:vorticity-fluttering-plate-closeup}
\end{figure}
The density of the plate and fluid are $\rhos = 2.7$~$\textrm{g}\cdot\textrm{cm}^{-3}$ and $\rhof = 0.997$~$\textrm{g}\cdot\textrm{cm}^{-3}$, respectively. To achieve a Reynolds number
of $\Re=1147$, the fluid viscosity is set to $\muf = 8.87 \times 10^{-3}$~$\textrm{g}\cdot\textrm{cm}^{-1}\cdot\textrm{s}^{-1}$.
A schematic of the problem setup is shown in Fig.~\ref{fig:trajectory-vorticity-fluttering-plate}(a).
The computational domain is $\Omega = [-40 \, \textrm{cm} , 40 \, \textrm{cm}] \times [-56 \, \textrm{cm}, 24 \, \textrm{cm}]$, 
a square of size $L_x \times L_y = 80 \, \textrm{cm} \times 80 \, \textrm{cm}$.
The center of the plate is initially located at $(x_0, y_0)=(5 \, \textrm{cm},0 \, \textrm{cm})$ 
with an initial angle of $\theta_0 = 60^{\circ}$ with respect to the $x$-axis.
Zero normal traction and tangential velocity conditions are imposed at the top boundary. 
Along the left and right boundaries zero normal velocity and tangential traction are imposed.
The no-slip condition is imposed at the bottom boundary.
\begin{figure}[t!!!]
		\centering
			\includegraphics[width=0.92\textwidth]{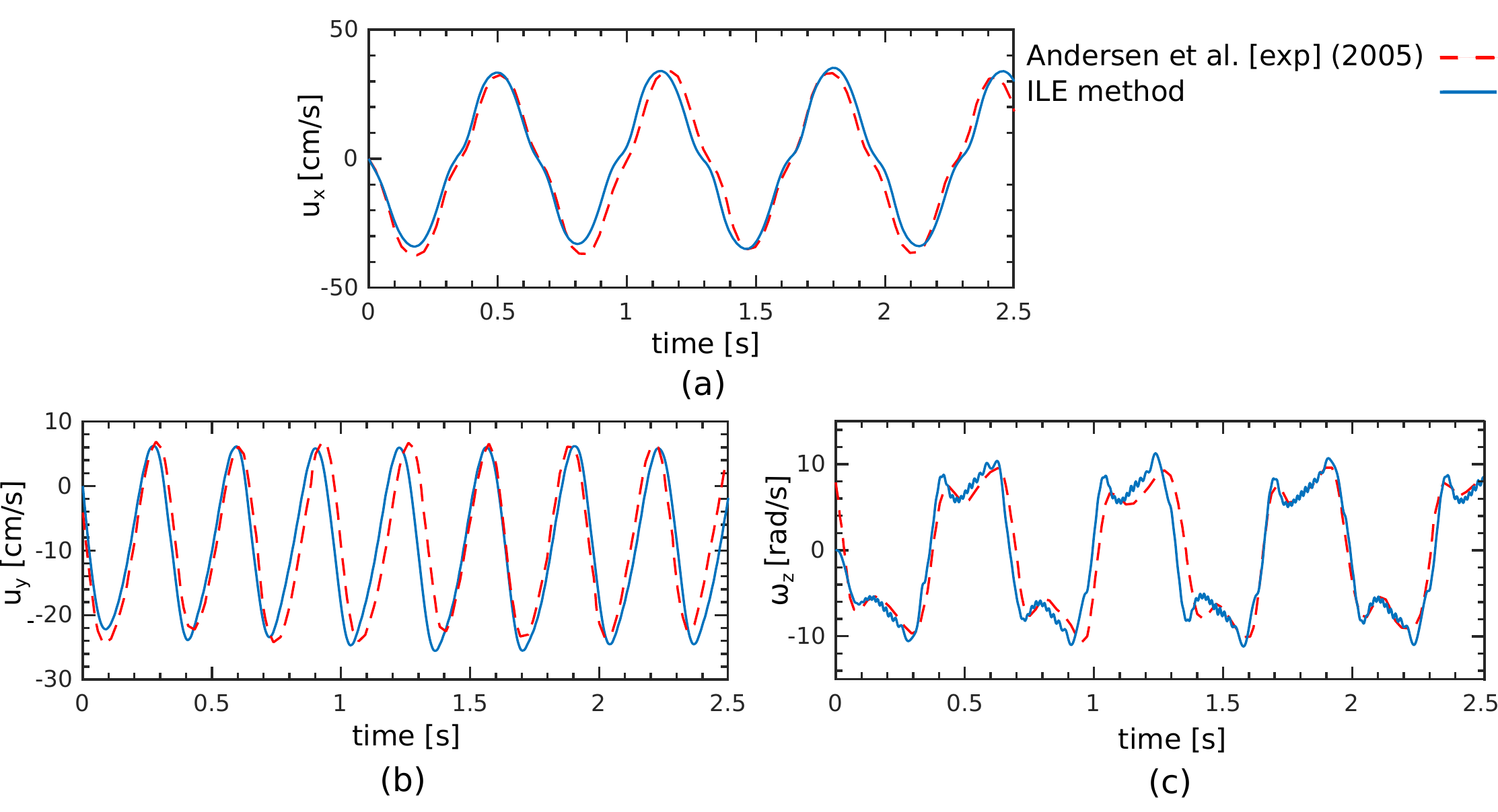}
\caption{Comparison of the (a) horizontal velocity, (b) vertical velocity, and (c) angular velocity between the 
ILE method (solid blue lines) and the experimental data of Andersen et al.~\cite{andersen2005unsteady} (dashed red lines) for the fluttering plate (Sec.~\ref{subsubsec:fluttering_motion}).} 
\label{fig:velocity-vorticity-fluttering-plate-comparison}
\end{figure} 

The domain is discretized using $N=6$ nested grid levels, with coarse grid spacing $h_{\textrm{coarsest}}=L_y/16 = 5 \, \textrm{cm}$ and refinement ratio $r = 4$ between levels,
leading to $h_{\textrm{finest}} = 0.00488 \, \textrm{cm}$.
Using $\mfac=1.8$, this leads to the thickness of the plate being discretized by approximately 10 linear elements. 
A constant time step size of $\Delta t=0.01$~ms is used, and the penalty spring constant is $\kappa=7.45 \times 10^{5}$~$\textrm{g}\cdot\textrm{cm}^{-2}\cdot\textrm{s}^{-2}$.
Figs.~\ref{fig:trajectory-vorticity-fluttering-plate}(b) and ~\ref{fig:trajectory-vorticity-fluttering-plate}(c) shows the overall dynamics of the plate during its fluttering free fall. 
Fig.~\ref{fig:trajectory-vorticity-fluttering-plate}(b) shows
the vortex structure over the course of the simulation.
Periodic fluttering motion is clearly observed. 
Fig.~\ref{fig:trajectory-vorticity-fluttering-plate}(c) shows the trajectory of the center of mass
along with the experimental data of Andersen et al.~\cite{andersen2005unsteady}.
Overall, our numerical results are in good agreement with the experimental data for the gliding motion of the plate from side to side as it flutters in its free fall.
Discrepancies around turning points can be, in part, attributed to the interaction of the complex vortex dynamics with sharp corners of the rectangular object.
There are also uncertain differences in the simulation and experimental operating conditions that could explain these differences.
Future studies should further investigate the accuracy of the present IIM coupling strategy for objects with sharp corners, including assessing grid sensitivity.

Fig.~\ref{fig:vorticity-fluttering-plate-closeup} shows close-up views of the fluttering plate at six different time points.
With a low angle of attack, the plate glides from a turning point on one side to reach a new turning point on the opposite side.
Towards the end of its glide and before it changes direction, there is a sharp increase in the magnitude of the angular velocity.
Because of the sharp edge of the structure,
the flow separates on the lower surface (Fig.~\ref{fig:vorticity-fluttering-plate-closeup}(c)) and shortly thereafter flow separation also 
occurs on the upper surface as well (Fig.~\ref{fig:vorticity-fluttering-plate-closeup}(d)). As the plate pitches upward,
it begins to glide in the opposite direction, and this process repeats itself in a periodic manner.
The interaction of vortices in the locations where the plate reverses direction creates a complex pattern of vortices, as shown in the figure.
Fig.~\ref{fig:velocity-vorticity-fluttering-plate-comparison} shows
the time history of the horizontal and vertical velocity components as well as the angular velocity
along with a comparison to the experimental results \cite{andersen2005unsteady}. Overall, there is very good agreement between the two results despite the complex dynamics.
Although the algorithm is able to correctly predict the periodic dynamics, 
	small deviations from measured data are observed particularly around troughs and crests of both translational and angular velocity plots. 
	We speculate that these deviations could potentially be due to sensitivity of the results to the complex interaction of the shed vortices 
	with the sharp corners, or possibly other uncertainties that are not accounted for in the comparison.

 \begin{figure}[b!!!]
		\centering
			\includegraphics[width=0.4\textwidth]{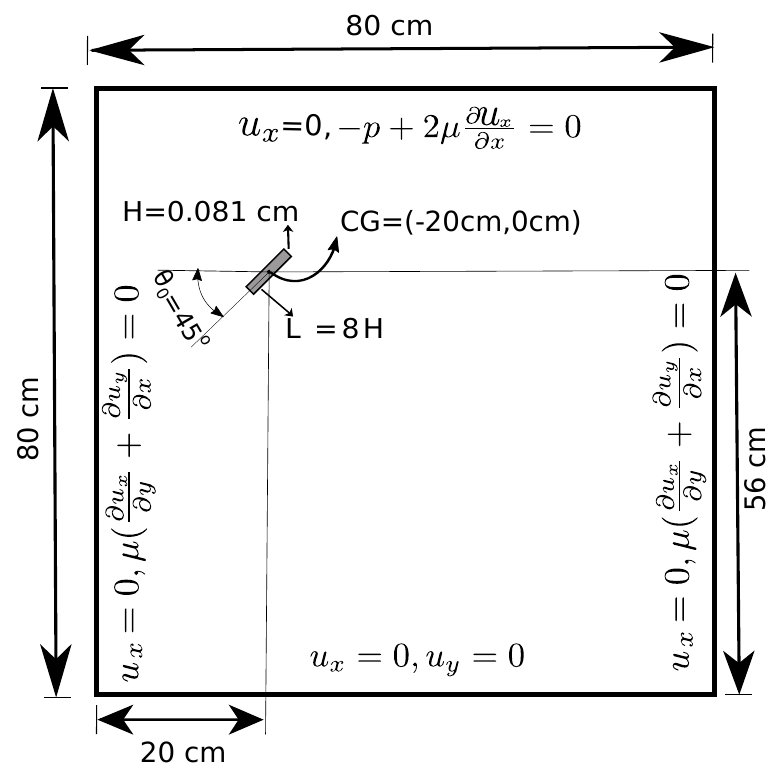}
\caption{Schematic diagram of the computational domain and boundary conditions for freely falling tumbling plate (Sec.~\ref{subsubsec:tumbling_motion}).} 
\label{fig:schematic-tumbling_motion}
\end{figure} 
\subsubsection{Tumbling motion} 
\label{subsubsec:tumbling_motion}
For the tumbling case, the plate thickness is kept fixed at $ H = 8.1 \times 10^{-2}$~cm, but the width-to-thickness ratio is changed to $\Lambda = L/H = 8$.
The initial location and angle of attack are  $(x_0, y_0)=(-20 \, \textrm{cm},0 \, \textrm{cm})$ and $\theta_0 = -45^{\circ}$, respectively. 
Note that because of the different width of the plate in this example, the Reynolds number is $\Re=837$.
%%%%%%%%%%%%%%%%%%%%%%%%%%%%%%%%%%%%%%%%%%%%%%%%%%%%%%%%%%%%%%%%%%%%%%%%%%%%%%%%%%%%%%%%%%%%%%%%%%%%%%%%%%%%%%%%%%%%%%%%%%% 
The penalty spring constant is $\kappa=5.5 \times 10^{5}$~$\textrm{g}\cdot\textrm{cm}^{-2}\cdot\textrm{s}^{-2}$. 
The remaining simulation parameters, including
the computational domain extent and size, spatial resolution, time step size, boundary conditions, and fluid properties, are identical to the fluttering case;
see also the schematic in Fig.~\ref{fig:schematic-tumbling_motion}.  
\begin{figure}[t!!!]
		\centering
			\includegraphics[width=0.95\textwidth]{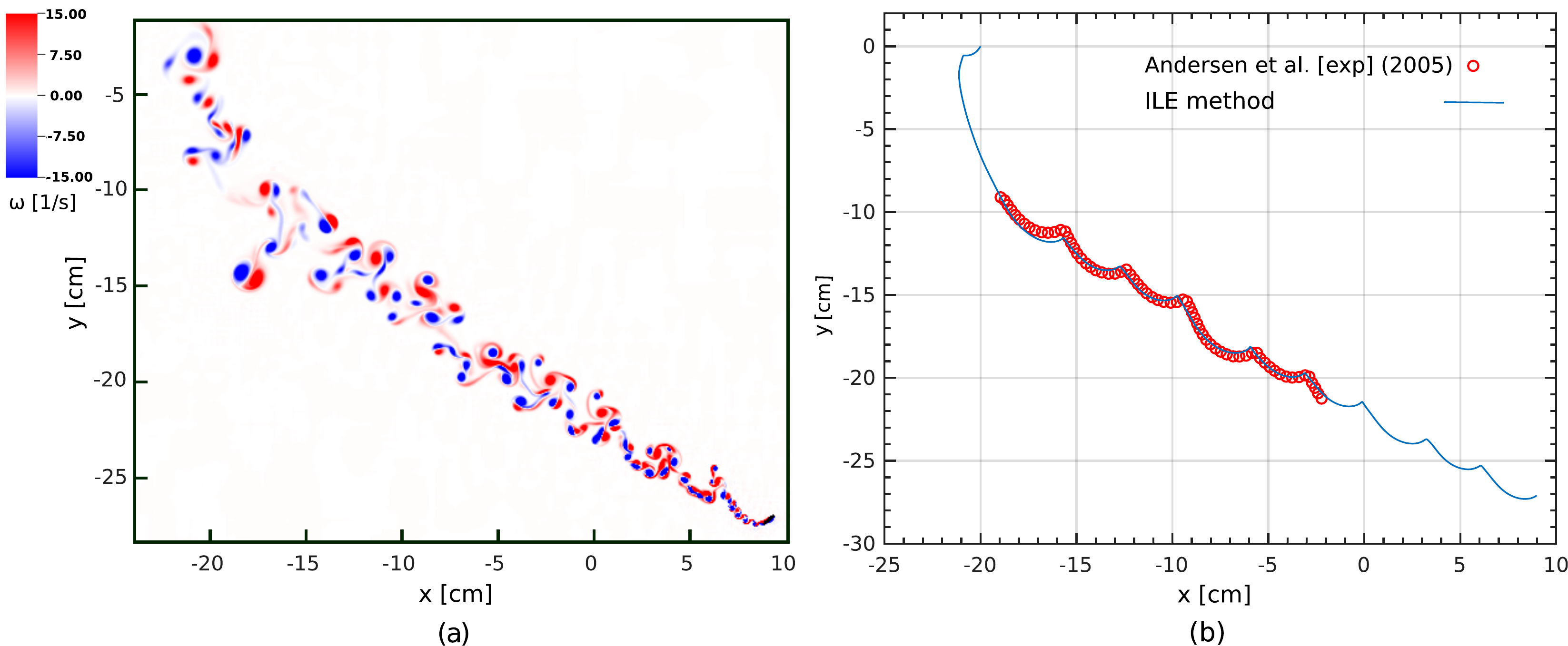}
\caption{Freely falling rectangular plate undergoing tumbling motion at $\Re=837$. (a) Vorticity field at time $t=2.06$~s. 
         (b) Comparison of the trajectory of the center of mass using the present ILE method (solid blue line) to the experimental data of Andersen et al.~\cite{andersen2005unsteady}.}
\label{fig:trajectory-vorticity-tumbling-plate}
\end{figure}
\begin{figure}[b!!!]
		\centering
			\includegraphics[width=0.9\textwidth]{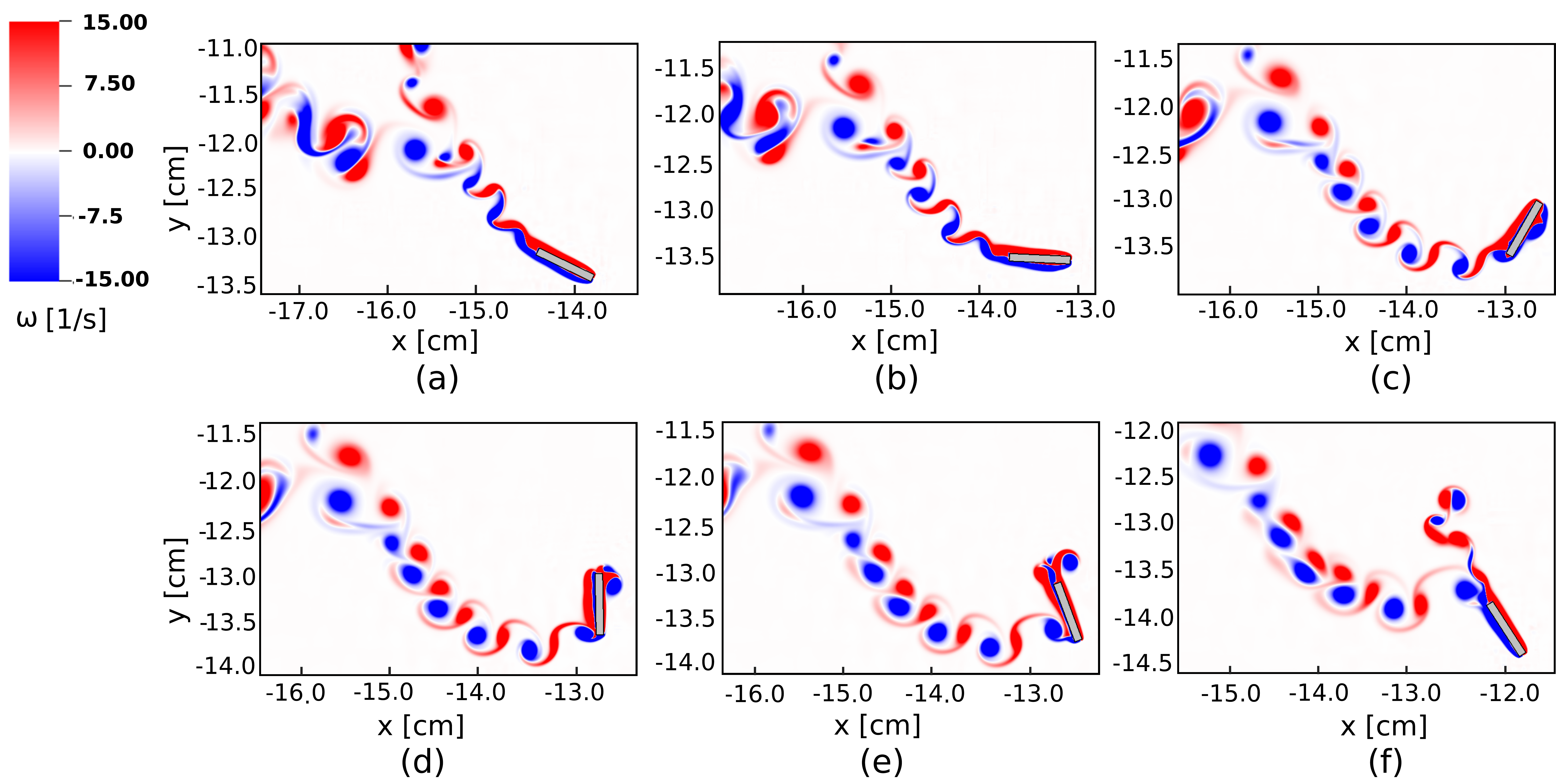}
\caption{Enlarged views of the vorticity field for the tumbling plate (see Fig.~\ref{fig:trajectory-vorticity-tumbling-plate}(a)) at times (a) $t=0.64$~s, (b) $t=0.66$~s, (c) $t=0.70$~s, (d) $t=0.72$~s, (e) $t=0.74$~s, and (f) $t=0.78$~s. } 
\label{fig:vorticity-closeup-tumbling-plate} 
\end{figure}
\begin{figure}[t!!!]
		\centering
			\includegraphics[width=0.7\textwidth]{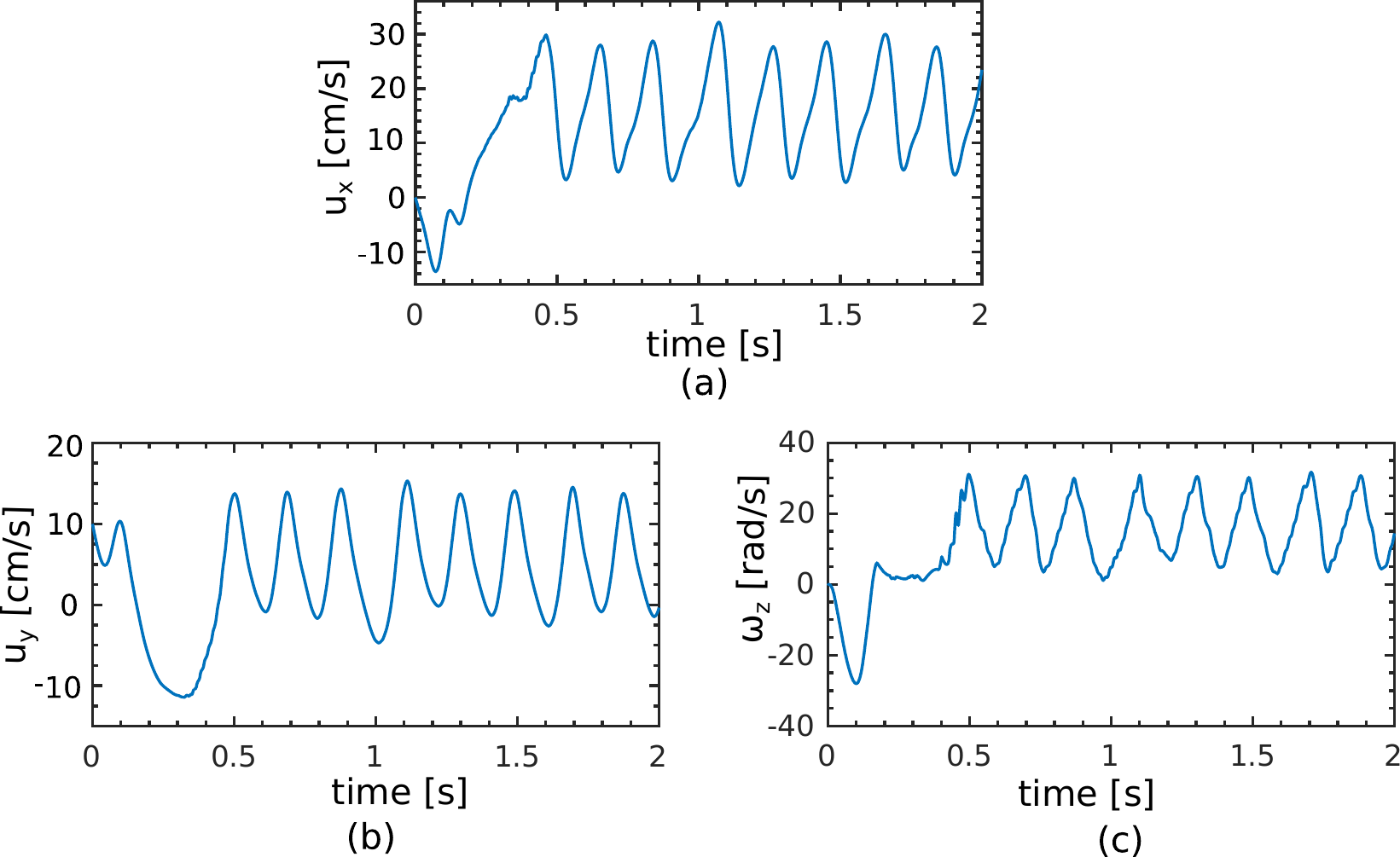}
\caption{Time history of the (a) horizontal velocity, (b) vertical velocity, and (c) angular velocity for the tumbling plate (Sec.~\ref{subsubsec:tumbling_motion}).}
\label{fig:uvw-tumbling-plate} 
\end{figure}
Fig.~\ref{fig:trajectory-vorticity-tumbling-plate} 
shows the overall dynamics of the tumbling plate.
As shown in Fig.~\ref{fig:trajectory-vorticity-tumbling-plate}(a), the complex vortex structure of the tumbling plate is well resolved by the simulation.
After the plate is released, it begins a gliding motion. 
Shortly after, it pitches upward, similar to the fluttering case. 
Because of the large angular momentum, however,
at the turning point the plate rotates more than $90^{\circ}$ clockwise. 
This large rotation, creates a large restoring moment that
causes the plate to rotate slightly counter-clockwise and then continue falling, with an inclination to the 
right side. 
The plate travels a path towards the bottom-right corner of the computational domain by a sequence of descending and accelerating motions.
At the turning points, the plate undergoes a full $360^{\circ}$ tumbling rotation, such that the lower surface during the gliding re-configures as upward facing.
Fig.~\ref{fig:trajectory-vorticity-tumbling-plate}(b) compares the vortex structure and trajectory of the tumbling plate 
to experimental measurements \cite{andersen2005unsteady}.
As in the fluttering case, the results are in very good agreement with the experimental data, and the trajectory of the plate agrees very well for 
the portion of the trajectory where experimental data are available. 
Fig.~\ref{fig:vorticity-closeup-tumbling-plate} shows enlarged views of the vorticity field at different times.
Flow separation on the lower side of the plate is clearly observed in Fig.~\ref{fig:vorticity-closeup-tumbling-plate}(c). This is followed
by separation of the flow on the opposite side as the plate rotates (Fig.~\ref{fig:vorticity-closeup-tumbling-plate}(d)). The gliding and diving 
towards the bottom left resumes at the end of each tumbling motion, as seen in Fig.~\ref{fig:vorticity-closeup-tumbling-plate}(f).

Fig.~\ref{fig:uvw-tumbling-plate} shows the time history of translational and rotational velocities.
Unlike the fluttering case, in which the vertical velocity appeared to have twice as larger frequency than the other two velocities,
here the number of periods demonstrates that approximately the same frequency is observed among all velocities. In addition, the significantly larger angular velocity in this case compared to the fluttering case in Fig.~\ref{fig:velocity-vorticity-fluttering-plate-comparison}
indicates faster rotation near the turning points.
The average horizontal and vertical velocity components and the average angular velocity obtained from three full cycles
of the present simulation are found to be $\overline{u_x}=15.87$~$\textrm{cm}\cdot\textrm{s}^{-1}$, $\overline{u_y}=-11.32$~$\textrm{cm}\cdot\textrm{s}^{-1}$, 
and $\overline{\omega_z}=15.95$~$\textrm{rad}\cdot\textrm{s}$.
For comparison, the experimental measurements of Andersen et al.~\cite{andersen2005unsteady} 
for the same average velocities are 
$\overline{u_x}=15.94$~$\textrm{cm}\cdot\textrm{s}^{-1}$, $\overline{u_y}=-11.5$~$\textrm{cm}\cdot\textrm{s}^{-1}$, and 
$\overline{\omega_z}=14.5$~$\textrm{rad}\cdot\textrm{s}^{-1}$.
The relative discrepancies in these quantities are $0.19$\%, $1.57$\%, and $9.93$\%, respectively.

%%%%%%%%%%%%%%%%%%%%%%%%%%%%%%%%%%%%%%%%%%%%%%%%%%%%%%%%%%%%%%%%%%%%%%%%%%%%%%%%%%%%%%%%%%%%%%%%%%%%%%%%%%%%%%%%%%%%%%%%%%%
\begin{figure}[t!!]
		\centering
			\includegraphics[width=0.7\textwidth]{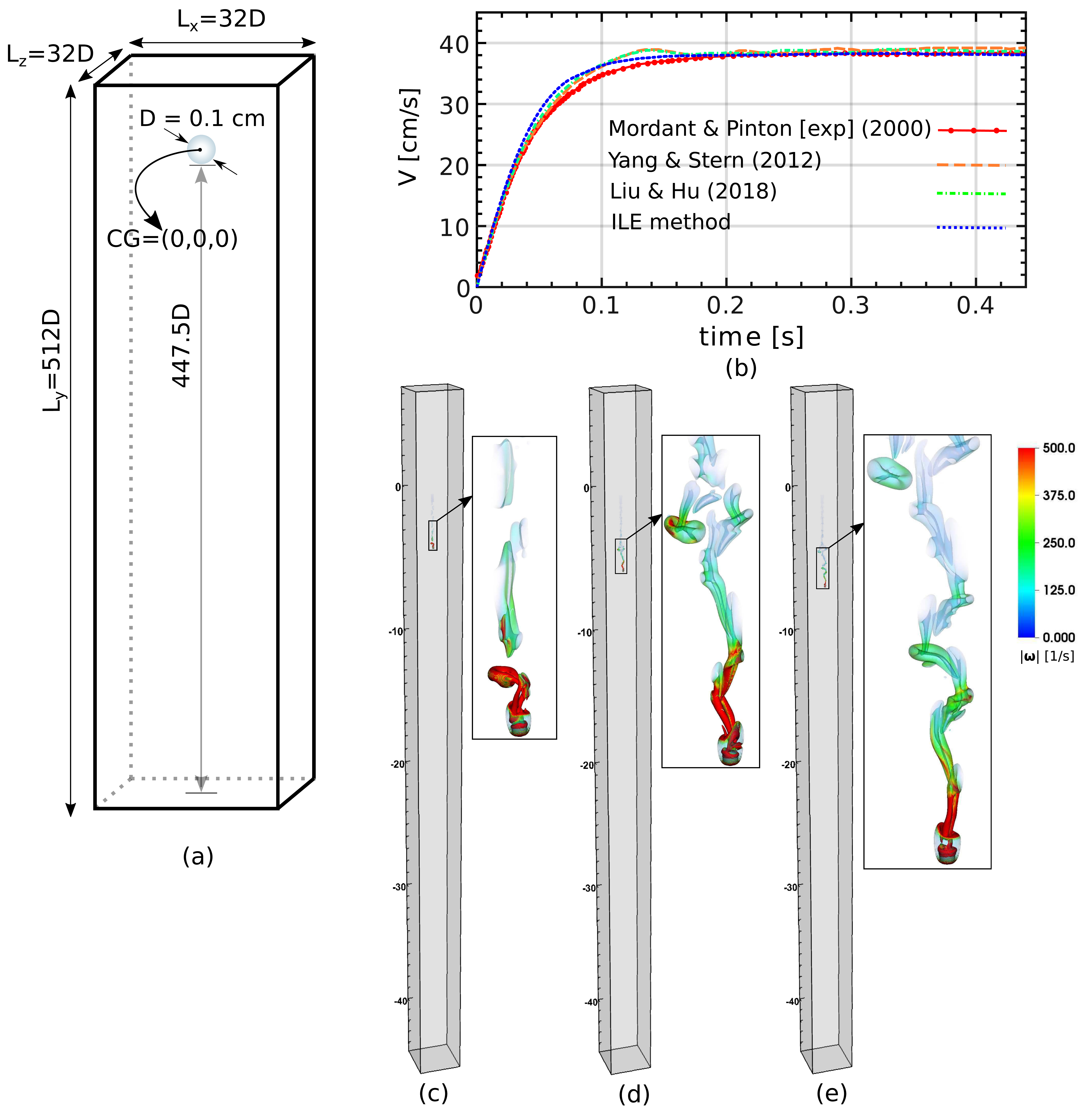}
\caption{ (a) A schematic of the initial setup for the freely falling steel bead (Sec.~\ref{subsec:falling_dense_sphere}) 
with $\rhos/\rhof=7.8$. Note that lengths in the schematic are not to scale. (b) Vertical velocity of the freely falling steel bead, 
including a comparison to the experimental data of Mordant and Pinton \cite{mordant2000velocity}
and the numerical results of Yang and Stern \cite{yang2012simple} and Liu and Hu \cite{liu2018block}.
(c)--(e) Instantaneous wake vortex patterns illustrated as isosurfaces of the Q criterion 
at the time $t=0.132$~s, $t=0.146$~s, and $t=0.179$~s, respectively.
} 
\label{fig:free-falling-steel-voritcity} 
\end{figure}
 
%%%%%%%%%%%%%%%%%%%%%%%%%%%%%%%%%%%%%%%%%%%%%%%%%%%%%%%%%%%%%%%%%%%%%%%%%%%%%%%%%%%%%%%%%%%%%%%%%%%%%%%%%%%%%%%%%%%%%%%%%%%
 
\subsection{Free falling of a dense sphere}
\label{subsec:falling_dense_sphere}
This section investigates the dynamics of a freely falling steel bead in 
water using an unconstrained rigid-body structure model.
This problem follows the experimental setup of Mordant and Pinton \cite{mordant2000velocity}.
The density of the steel bead is $\rhos=7.85$~$\textrm{g}\cdot\textrm{cm}^{-3}$, and the diameter is $D=0.1$~cm.
The Reynolds number is $\Re=(\rhof \bar{V} D)/\muf =430$, and a terminal velocity of $\bar{V}=38.3$~$\textrm{cm}\cdot\textrm{s}^{-1}$ is reported in the original work.
With a fluid density of $\rhof=997$~$\textrm{g}\cdot\textrm{cm}^{-3}$, this results in a dynamic viscosity of $\muf=8.88\times 10^{-3}$~$\textrm{g}\cdot\textrm{cm}^{-1}\cdot\textrm{s}^{-1}$. 
The computational domain is $\Omega = [-16 D,16 D]\times[-448 \, D,64 \, D]\times[-16 \, D,16 \, D]$,
a rectangular cuboid of size $L_x \times L_y \times L_z = 32 D \, \times 512 D \, \times 32 D $; see Fig.~\ref{fig:free-falling-steel-voritcity}(a) for a schematic of the problem setup.
Zero normal traction and tangential velocity conditions are imposed at the top boundary.
Along the peripheral walls zero normal velocity and tangential traction are imposed.
The no-slip condition is imposed at the bottom boundary.
The domain is discretized using $N=3$ nested grid levels, with coarse grid spacing 
$h_{\textrm{coarsest}}=L_x/64 = 0.05 \, \textrm{cm}$ and refinement ratio $r = 4$ between levels,
leading to $h_{\textrm{finest}}=0.003125 \, \textrm{cm}$.
The volumetric mesh of the sphere consists of hexahedral elements leading to a surface mesh composed of  bilinear quadrilateral elements with $\Mfac=2$. 
A fixed time step size of $\dt=0.01$~ms is used, and the penalty spring constant is $\kappa=5 \times 10^{5}$~$\textrm{g}\cdot\textrm{cm}^{-2}\cdot \textrm{s}^{-2}$. 
The center of the sphere is initially positioned at the origin and released with zero initial translational and angular velocities. 
Fig.~\ref{fig:free-falling-steel-voritcity}(b) shows the time history of the vertical velocity. 
After being released, the steel bead starts to accelerate until reaching its terminal velocity. 
The acceleration of the bead is in good agreement
 with the experimental data of Mordant and Pinton \cite{mordant2000velocity}, and the terminal velocity
 is also in excellent agreement with the experimental measurement.
 Fig.~\ref{fig:free-falling-steel-voritcity}(c)--(e) show isosurfaces of the Q-criterion \cite{hunt1988eddies} to visualize the vortex dynamics.  
Note that there is no constraint in the motion of the sphere, and that all of the degrees of freedom including the rotational ones are included in the solution.
Shortly after the sphere reaches its terminal velocity $\bar{V}=38.3$~$\text{cm}\cdot\textrm{s}^{-1}$, the flow behind the sphere becomes
irregular due to the moderately large Reynolds number. Although small lateral movements
are observed in the simulation with the deviation of center of mass reaching to about $0.25D = 8 h_{\textrm{finest}}$, no consistent pattern of ``zig-zagging'' motion is apparent.
Fig.~\ref{fig:free-falling-steel-voritcity} also compares our results to two other numerical studies \cite{yang2012simple,liu2018block}. 

%%%%%%%%%%%%%%%%%%%%%%%%%%%%%%%%%%%%%%%%%%%%%%%%%%%%%%%%%%%%%%%%%%%%%%%%%%%%%%%%%%%%%%%%%%%%%%%%%%%%%%%%%%%%%%%%%%%%%%%%%%%

\begin{table}[t!!!!!!!!!!!!]
	\centering	
	\caption{Parameters for freely falling sphere with near-unity density ratio (Sec.~\ref{subsec:falling_light_sphere})}
	\label{table:Cate_exp_info}	
\begin{tabular}{l*{6}{c}r}
\hline
Case  & $\rhof \, [\textrm{g}\cdot\textrm{cm}^{-3}]$ & $\rhos/\rhof$ & $\muf \, [\textrm{g}\cdot\textrm{cm}^{-1} \cdot \textrm{s}^{-1}]$ & $\bar{V} \, [\textrm{cm}\cdot\textrm{s}^{-1}]$ & $\Re$  \\
\hline
1 (case 1 of ten Cate et al.~\cite{ten2002particle})    & 0.970 & 1.155  & 0.0373  & 3.8 & 1.5  \\
2 (case 4 of ten Cate et al.~\cite{ten2002particle})         & 0.960  & 1.167 & 0.058 & 12.8 & 31.9  
\end{tabular}
\end{table}
\subsection{Free falling of a sphere with near-unity density ratio}
\label{subsec:falling_light_sphere}
\begin{figure}[b!!]
		\centering
			\includegraphics[width=0.85\textwidth]{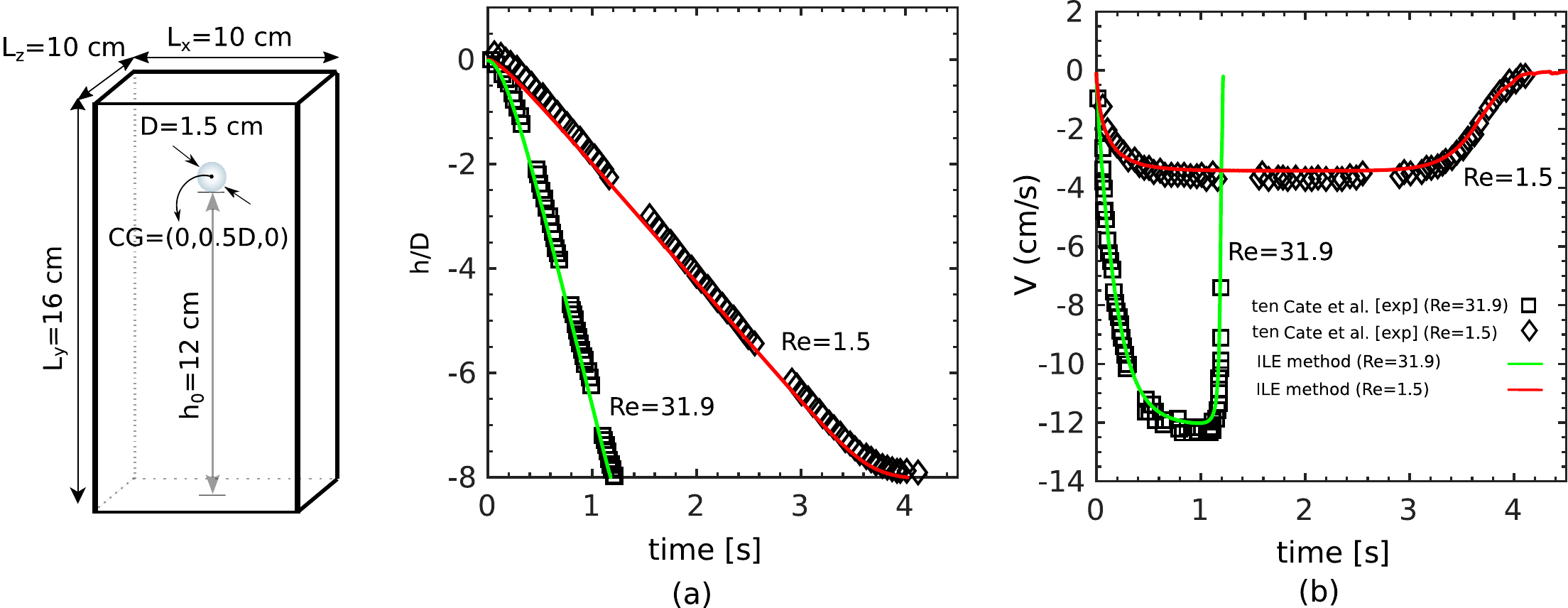}
\caption{(a) A schematic of the initial setup of a single sphere in a small flow chamber with near unity solid-fluid density ratio (Sec.~\ref{subsec:falling_light_sphere}).
(b) Time history of the vertical position and (c) vertical sedimentation velocity.}
\label{fig:cate-experiment}  
\end{figure}

This section explores the performance of the ILE methodology with respect to two challenging aspects of FSI: the influence of the wall and the ability of the method to handle near-unity 
density ratios in a fully unconstrained rigid-body motion.
This case is based on the experimental data reported by ten Cate et al.~\cite{ten2002particle} of a sphere
settling in a confined flow chamber.
The sphere has diameter $D=1.5$~cm and density $\rhos=1.120$~$\textrm{g}\cdot\textrm{cm}^{-3}$, which is close to that of the  surrounding fluid; see Table~\ref{table:Cate_exp_info}.
The Reynolds number is $\Re=\rhos \bar{V} D /\muf$, in which the terminal velocity $\overline{V}$ is taken to be the settling velocity in an infinite domain.
Two different cases with Reynolds numbers of $\Re=1.5$ and $\Re=31.9$ are considered, which are achieved in the experiments by varying the density and viscosity of the fluid.
The computational domain is $\Omega = [-5 \, \textrm{cm}, 5 \, \textrm{cm}]\times[-12 \, \textrm{cm}, 4 \, \textrm{cm}]\times[-5 \, \textrm{cm}, 5 \, \textrm{cm}]$,
a rectangular cuboid of size $L_x \times L_y \times L_z = 10 \, \textrm{cm} \times 16 \, \textrm{cm} \times 10 \, \textrm{cm}$.
The sphere is initially positioned in the $x$-$z$ mid-plane 
and towards the top of the box at a height of $12$ cm, measured from the bottom of the sphere to the bottom of the box; see Fig.~\ref{fig:cate-experiment}(a).
Zero normal traction and tangential velocity conditions are applied at the top boundary while no-slip condition is imposed at all other boundaries.

The domain is discretized using $N=3$ nested grid levels, with coarse grid spacing $h_{\textrm{coarsest}}=L_x/32 =0.3125 \, \textrm{cm}$ and refinement ratio $r = 2$ between levels,
leading to $h_{\textrm{finest}}=0.078125 \, \textrm{cm}$ .
The volumetric mesh of the sphere consists of hexahedral elements leading to a surface mesh composed of bilinear quadrilateral elements with $\Mfac=2$. 
A fixed time step size of $\dt=0.5$~ms is used, and the penalty spring constant is $\kappa=5 \times 10^{4}$~$\textrm{g}\cdot\textrm{cm}^{-2}\cdot \textrm{s}^{-2}$.
Because of the confinement effect of the walls, we expect the settling of the sphere to be different from a free falling sphere in an infinite medium.
Upon release, the particle accelerates until reaching its terminal velocity, which occurs quickly in the lower Reynolds number case.
The sphere then starts to decelerate as it approaches to the bottom the flow chamber. Fig.~\ref{fig:cate-experiment} panels (b)--(c)
show the time histories of the non-dimensional vertical gap between the sphere and bottom wall (h/D) and the 
vertical velocity of the sphere ($\overline{V}$) for both cases. 
The method yields excellent agreement with the 
experimental measurements of ten Cate et al.~\cite{ten2002particle} for both the gap size over time and the terminal velocity. Specifically, the present method 
captures the dynamics of the sphere while also remaining stable even as the sphere reaches the bottom of the box. For the lower Reynolds number case, the sphere has 
more time to travel through the box because of its lower speed, and this results in a larger flat region, when
the sphere experiences a nearly constant downward velocity, before it starts to decelerate as it approaches the bottom wall. Indeed, at $\Re=1.5$, we observe stable 
sedimentary motion of the sphere in the simulation even past the latest time reported in the experiment.
At $\Re=31.9$, the sphere does not experience an extended period of time at terminal-like velocity as it quickly reaches the bottom of the box. 
The simulation is also stable for this case even after the sphere impacts the bottom; see Fig.~\ref{fig:cate-experiment}(b).

%%%%%%%%%%%%%%%%%%%%%%%%%%%%%%%%%%%%%%%%%%%%%%%%%%%%%%%%%%%%%%%%%%%%%%%%%%%%%%%%%%%%%%%%%%%%%%%%%%%%%%%%%%%%%%%%%%%%%%%%%%%
\subsection{Freely rising sphere} 
\label{subsec:rising_sphere}

 \begin{figure}[b!!!]
		\centering
			\includegraphics[width=0.15\textwidth]{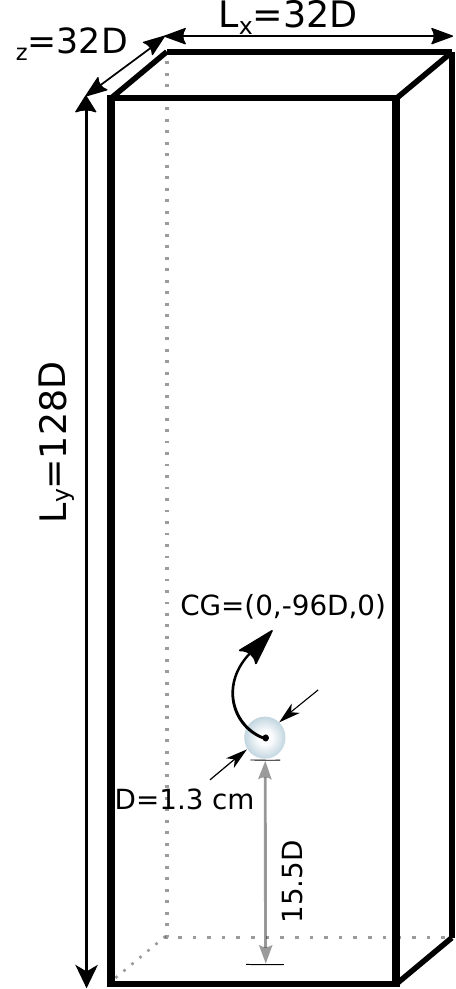}
\caption{Schematic diagram of the computational domain for freely rising sphere (Sec.~\ref{subsec:rising_sphere}).} 
\label{fig:schematic-rising-sphere}
\end{figure}
\begin{figure}[t!!]
		\centering
			\includegraphics[width=0.85\textwidth]{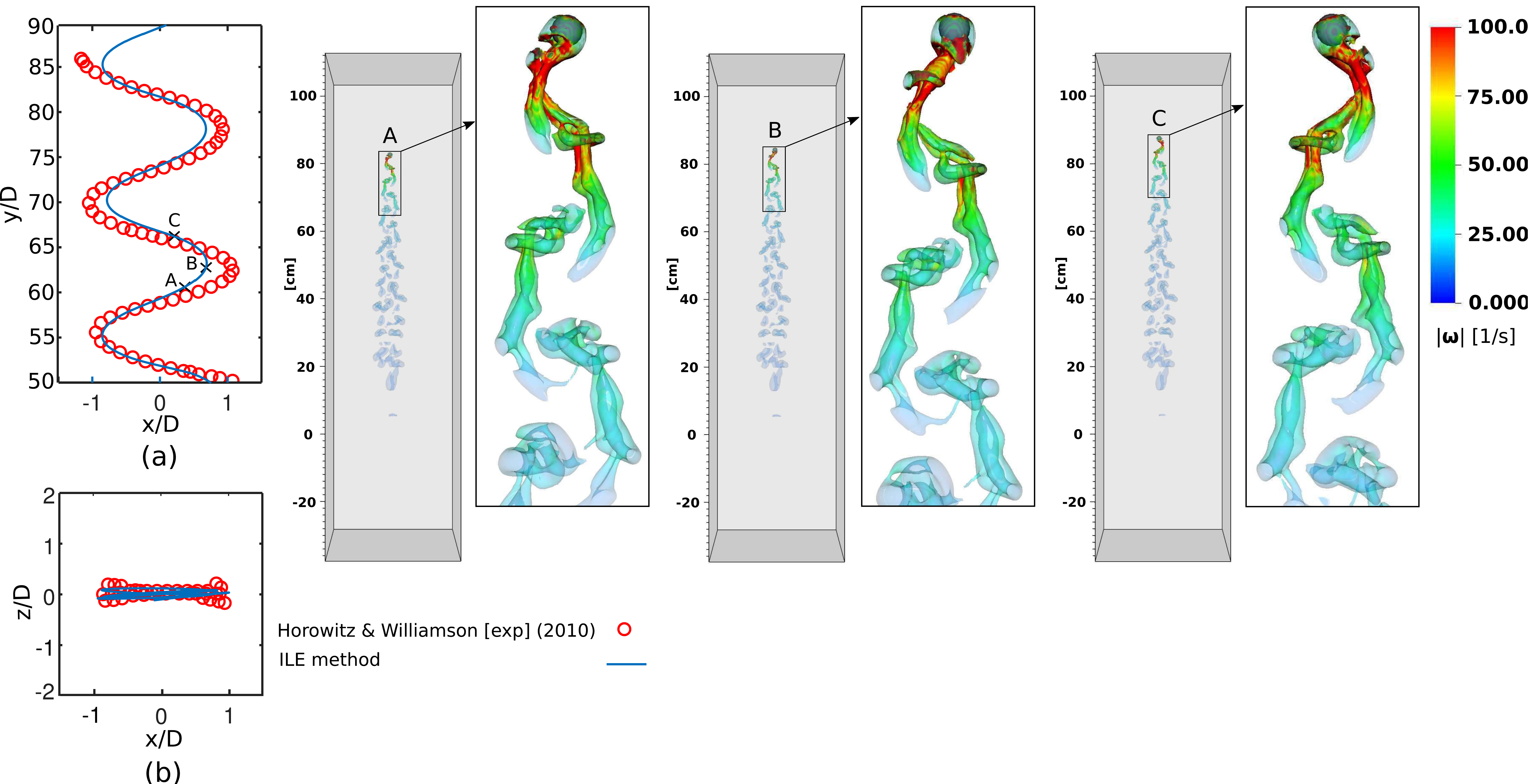}
\caption{Trajectory of the center of mass for the rising sphere (Sec.~\ref{subsec:rising_sphere}) in the (a) $x$-$y$ and (b) $x$-$z$ planes along with the experimental data of Horowitz and Williamson \cite{horowitz2010effect}.
		The vortex wake structure for the rising sphere at locations A, B, and C shown on the $x$-$y$ plane trajectory. Simulation parameters 
		include $\Re=450$, $m^{*}=0.11$, and $D=1.3$~cm.} 
\label{fig:rising-sphere}  
\end{figure}  
This section considers a rising sphere using a fully unconstrained rigid-body structural model.
Such problems can pose substantial challenges to FSI algorithms using weak coupling schemes, which can become unstable \cite{liu2018block}.
Experimental results of Horowitz and Williamson \cite{horowitz2010effect} are used as a benchmark 
to validate the dynamics generated by our numerical method. We consider a challenging case in which the parameters
are chosen such that the sphere oscillates periodically and vigorously in a  ``zig-zag'' trajectory within a tight vertical plane.
The sphere has diameter $D=1.3$~cm and mass ratio $m^{*}=0.11$.
The experiments are at $\Re = \rhof U D/\muf = 450$. To match the experimental Reynolds number,
the dynamic viscosity is $\muf = 0.125$~$\textrm{g}\cdot(\textrm{cm}\cdot \textrm{s})^{-1}$.
The computational domain is $\Omega = [-16 \, D,16 \, D]\times[-112 \, D,16 \, D]\times[-16 \, D,16 \, D]$; see the schematic in Fig.~\ref{fig:schematic-rising-sphere}.
The domain is discretized using $N=3$ nested grid levels, with coarse grid 
spacing $h_{\textrm{coarsest}}=L_x/64=0.65 \, \textrm{cm}$ and refinement ratio $r = 4$ between levels,
leading to $h_{\textrm{finest}} \approx 0.04 \, \textrm{cm}$.
The volumetric mesh of the sphere consists of hexahedral elements leading to a surface mesh composed of bilinear quadrilateral elements with $\Mfac=2$. 
The center of the sphere is initially positioned at $(0, -96 \, D, 0)$ and is released with zero initial translational and angular velocity.
Zero normal traction and tangential velocity conditions are imposed at the top boundary.
The no-slip condition is imposed at the bottom boundary.
Along the peripheral walls zero normal velocity and tangential traction are imposed.
The no-slip condition is imposed at the bottom boundary.
A fixed time step size of $\dt=0.01$~ms is used, and the penalty spring constant is $\kappa=5.2 \times 10^{5}$~$\textrm{g}\cdot(\textrm{cm}\cdot \textrm{s})^{-2}$.

Fig.~\ref{fig:rising-sphere} demonstrates that the sphere motion is approximately planar, and captures the periodic zig-zag trajectory observed experimentally. We emphasize that in our simulation, 
no constraints are imposed on the motion of the sphere. Specifically we do not impose either planar motion or the zig-zag trajectory. 
According to Horowitz and Williamson \cite{horowitz2010effect}, the periodic motion of the sphere at this particular density ratio
resembles the dynamics at much higher Reynolds number and is always confined to a plane. 
The simulation results 
predict the same kind of in-plane motion as observed in the experiment. As the sphere goes through the zig-zag motion it creates a complex
but organized wake pattern with vortex rings forming at the turning points. The numerical simulation clearly captures the vortices behind the sphere.
Further, there is generally excellent agreement with the trajectory obtained from the experiment.

%%%%%%%%%%%%%%%%%%%%%%%%%%%%%%%%%%%%%%%%%%%%%%%%%%%%%%%%%%%%%%%%%%%%%%%%%%%%%%%%%%%%%%%%%%%%%%%%%%%%%%%%%%%%%%%%%%%%%%%%%%% 
\subsection{Bileaflet mechanical heart valve at physiological conditions}
\label{subsec:MHV}

 \begin{figure}[b!!]
		\centering
			\includegraphics[width=0.98\textwidth]{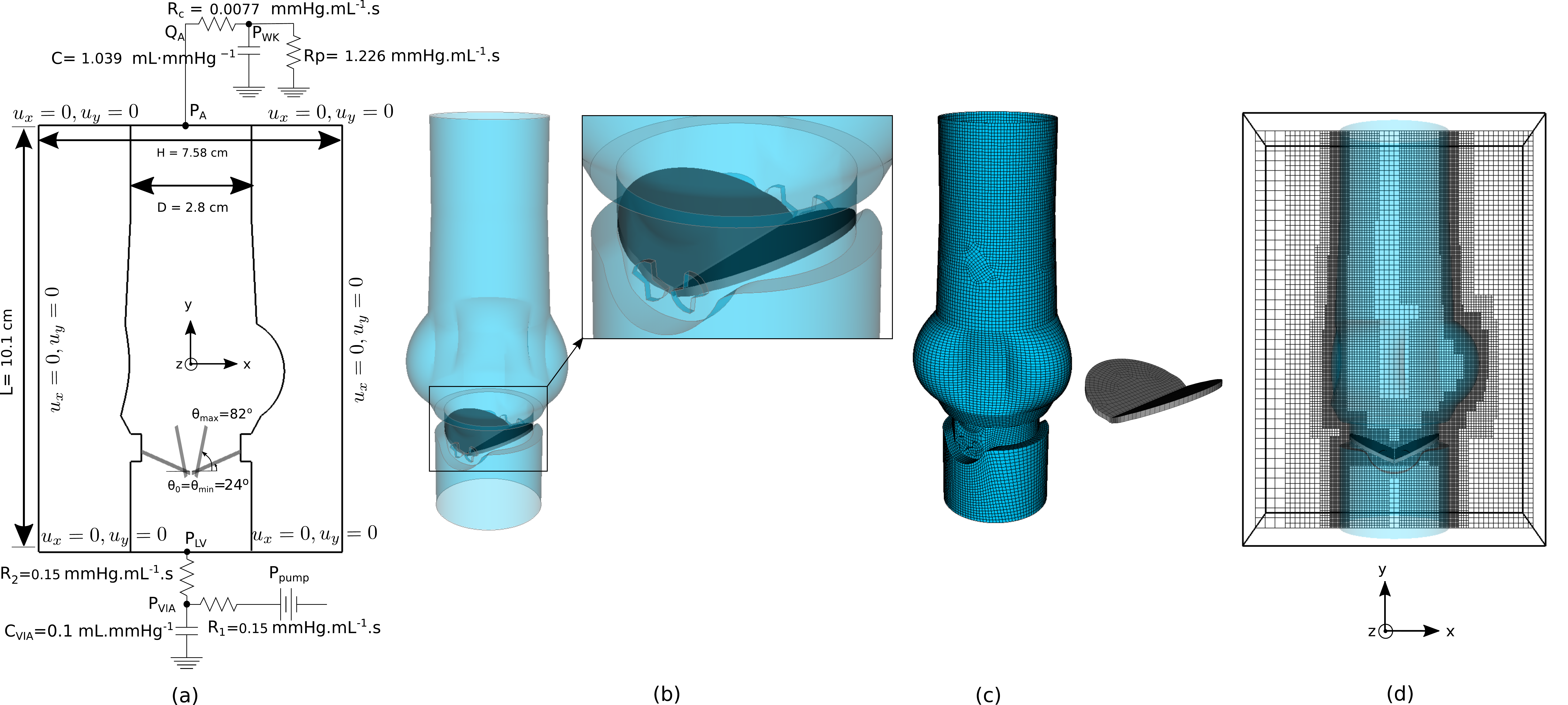}
\caption{ Model setup of the bileaflet mechanical heart valve (Sec.~\ref{subsec:MHV}). 
(a) Schematic cross-section view of the aortic test section on the $x$-$y$ mid-plane plane showing the dimensions and boundary conditions.
A three-element Windkessel model is used at the downstream and upstream of the test section while all other boundaries are set to solid wall boundary.
(b) Assembly of the aortic test section including the bileaflet valve and their position around the hinges.
(c) Computational mesh of the aortic test section and the valve. (d) The computational domain in which the aortic test section is embedded in. The 
block-structured adaptively refined Cartesian grid is shown on the $x$-$y$ mid-plane.} 
\label{fig:MHV3D_mesh}
\end{figure}
\begin{figure}[t!!]
		\centering
			\includegraphics[width=0.68\textwidth]{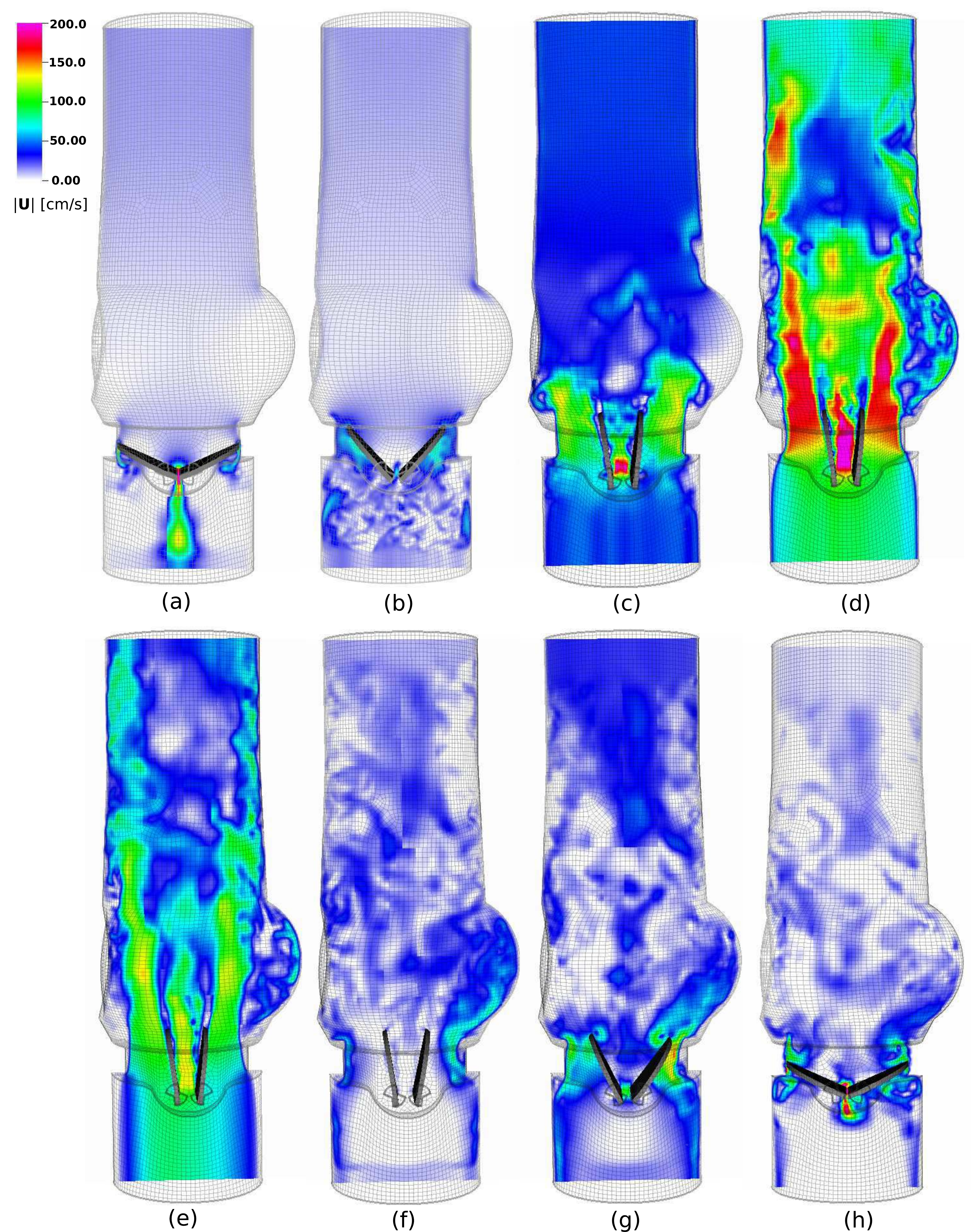}
\caption{Velocity magnitudes of bileaflet mechanical heart valve (Sec.~\ref{subsec:MHV}) on the $x$-$y$ mid-plane bisecting the valve, given at times (a) $t=0.0225$~s, (b) $t=0.1050$~s, (c) $t=0.1275$~s, (d) $t=0.1875$~s, (e) $t=0.2925$~s, (f) $t=0.3550$~s, (g) $t=0.3675$~s, and (h) $t=0.375$~s.} 
\label{fig:MHV3D_velocity}
\end{figure} 
Bileaflet mechanical heart valves remain commonly used in heart valve replacement because of their durability.
In this section, we consider the simulation of fluid flow through a geometrically realistic model of a 25 mm St Jude Medical Regent\textsuperscript{TM} bileaflet 
mechanical heart valve in the aortic test section of an experimental pulse duplicator platform \cite{Scotten2004,scotten2011importance}.
These systems are used in academia, 
industry, and regulatory agencies to assess the performance of prosthetic heart valves.
This simulation aims to demonstrate the ability of the present approach in modeling a realistically complex problem 
that involves multiple moving and stationary parts, 
including the small-scale feature of the hinge geometry.
These parts come in close contact with each other and undergo substantial pressure loading when closed.
In particular, we use experimental pressure and flow data obtained from a real pulse duplicator to establish realistic boundary models for the FSI model. 
These boundary models are calibrated in isolation from the rest of the system in that 
the motion of the leaflets, including the timing of valve opening and closing, is not prescribed and the simulated pressures, flow rates, and 
leaflet kinematics all emerge from integrating these three model components.
The experiments conducted to calibrate the model used saline as a test fluid, 
which we model as a Newtonian fluid with uniform density $\rhof =1.0$~$\textrm{g}\cdot\textrm{cm}^{-3}$ and uniform dynamic viscosity $\mu=1.0\,\textrm{cP}$.
The computational domain is 
$\Omega = [0, \, L_x]\times[0, \, L_y]\times[0, \, L_z]$ a rectangular cuboid of size
 $L_x\times L_y\times L_z = 7.07 \, \textrm{cm} \times 10.1 \, \textrm{cm} \times 7.07\, \textrm{cm}$.
Fig.~\ref{fig:MHV3D_mesh} shows the geometrically detailed three-dimensional model including the valve leaflets and the aortic test section of the pulse duplicator, 
the computational mesh, and relevant boundary conditions.
The Eulerian domain is discretized using $N=4$ nested grid levels, with coarse grid 
spacing $h_{\textrm{coarsest}}=L_x/34 \approx 0.208 \, \textrm{cm}$ and refinement ratio $r = 2$ between levels, 
leading to $h_{\textrm{finest}} \approx 0.026 \, \textrm{cm}$. 
The volumetric mesh of each valve leaflet consists of hexahedral elements leading to a surface representation composed of bilinear quadrilateral elements with $\Mfac \approx 2$.
The test section is a stationary surface described by bilinear quadrilateral elements with $\Mfac=3$.
The thickness of the valve leaflets is about $0.07$~cm and the density is $\rhos=1.8$~$\textrm{g}\cdot\textrm{cm}^{-3}$. The gap distance between the two valves is measured to be $0.0275$~cm. 
The penalty spring constant associated with the valves is set to $\kappa=5 \times 10^{6}$~$\textrm{g}\cdot(\textrm{cm}\cdot \textrm{s})^{-2}$.
The motion of each leaflet is constrained to consist only of rotation about predefined hinge axes.
A restoring spring-damper tortional force is used to keep the right valve leaflet restricted in 
its rotation between $\theta_{\textrm{min}}=24^{\circ}$ and $\theta_{\textrm{min}}=82^{\circ}$.
The same rule is applied to keep the left valve leaflet restricted with mirrored angle. 
The test section is kept stationary by means of spring-type penalty forces with $\kappa=1.5 \times 10^{6}$~$\textrm{g}\cdot\textrm{cm}^{-2}\cdot \textrm{s}^{-2}$.
A fixed time step size of $\Delta t=2 \, \mu\textrm{s}$ is used for the simulation.
Three-element Windkessel models establish downstream loading conditions and the upstream driving conditions
for the aortic test section.
A combination of normal traction and zero tangential velocity boundary conditions are used at the inlet and outlet to couple the 
reduced-order models to the detailed description of the flow within the test section. 
The values of the resistances and compliance for the upstream model are $C_{\textrm{VIA}}=0.1 \, \textrm{mL}\cdot \textrm{mmHg}^{-1}$ and
$R_1= R_2 = 0.15 \, \textrm{mmHg}\cdot\textrm{mL}^{-1}\cdot\textrm{s}$.
The values at the downstream model are $R_{\textrm{c}}=0.0077 \, \textrm{mmHg}\cdot\textrm{mL}^{-1}\cdot\textrm{s}$ ,
 $R_{\textrm{p}}=1.226 \, \textrm{mmHg}\cdot\textrm{mL}^{-1}\cdot\textrm{s}$, and $C =1.039 \, \textrm{mL}\cdot \textrm{mmHg}^{-1}$; 
see Griffith et al.~\cite{griffith2009simulating} and Lee et al.~\cite{lee2020fluid} for further discussion on the specification and parameterization of
the Windkessel models.
Solid wall boundary conditions are imposed on the remaining boundaries of the computational domain.

Fig.~\ref{fig:MHV3D_velocity} shows the velocity magnitudes on the plane bisecting the valves at different time points within the simulated cardiac cycle. At early times when the valves are closed
there exists a large pressure difference across the test section that forces jets of the fluid to escape from small gaps around the valves
with velocity magnitude that reaches about $450 \, \text{cm}\cdot\textrm{s}^{-1}$ at its peak. These hinge gap flows have been 
well characterized for bileaflet mechanical valves result from gaps between the leaflets and the housing and around the hinge areas 
in a fully closed position \cite{travis2001sensitivity,yoganathan2004fluid,scotten2011importance}.
These flows create a complex vortical interaction later in the diastole phase of the cardiac cycle, when the valve is closed and supporting a physiological pressure load.
Von Karman like vortex shedding 
clearly occurs during the systolic phase, when the valves are fully open. 
Many prior numerical simulations of mechanical heart valve models have been restricted to imposing either
experimental flow rates, or flow dynamics under prescribed leaflet motion.
Because the flow rate is not imposed in the present model, and because the time-dependent configuration of 
the valve determines the resistance of the aortic test section, this simulation demonstrates a non-trivial test of the numerical method.
 
 %%%%%%%%%%%%%%%%%%%%%%%%%%%%%%%%%%%%%%%%%%%%%%%%%%%%%%%%%%%%%%%%%%%%%%%%%%%%%%%%%%%%%%%%%%%%%%%%%%%%%%%%%%%%%%%%%%%%%%%%%%% 

\subsection{Transport of rigid blood clots in the inferior vena cava at exercise flow conditions}
\label{subsec:ivc_clot}

\begin{figure}[t!!]
		\centering
			\includegraphics[width=0.6\textwidth]{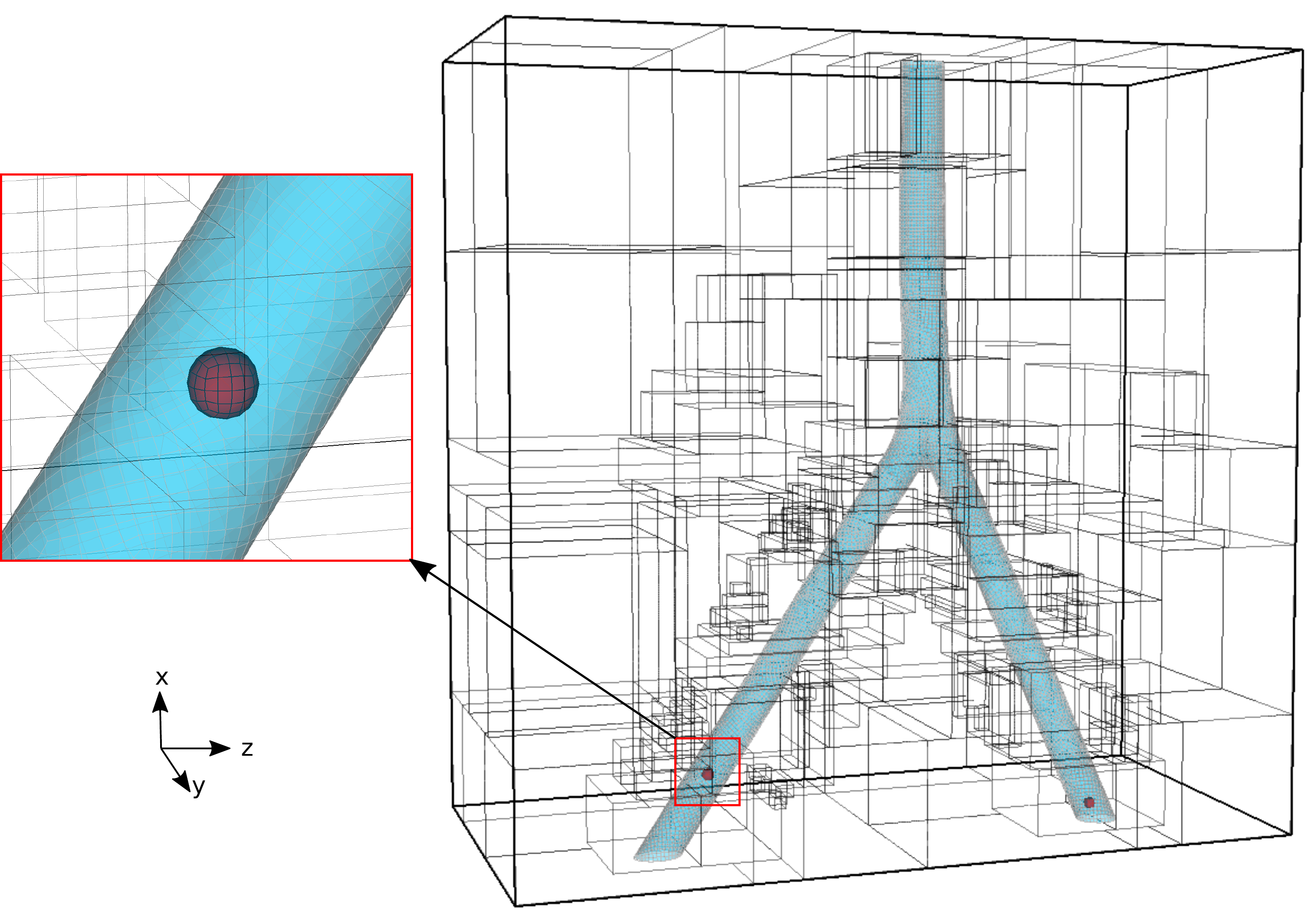}
\caption{Computational mesh of the patient-averaged model of the inferior vena cava (IVC), the locally refined Cartesian grid and rigid spherical clots insde IVC (Sec.~\ref{subsec:ivc_clot}).} 
\label{fig:Eulerian-mesh-ivc}
\end{figure} 
As a demonstration of the method's capability in modeling the motion of unconstrainted three-dimensional objects in a complex geometry, we simulate transport of rigid blood clots 
through the inferior vena cava (IVC) at a high flow rate of 100 cm$^{3}$/s, which corresponds to exercise flow conditions with a maximum Reynolds number of about $\Re=1500$.
The IVC is a large vein that transports deoxygenated blood from lower extremities of the body back to the right atrium of the heart.
The geometry of the IVC shown in Fig.~\ref{fig:Eulerian-mesh-ivc} is a modified version of the patient-averaged model by Rahbar et al.~\cite{rahbar2011three} 
that has been recently used in studies of the hemodynamics~\cite{craven2018steady,gallagher2018steady,kolahdouz2020immersed}. 
\begin{figure}[b!!]
		\centering
			\includegraphics[width=0.95\textwidth]{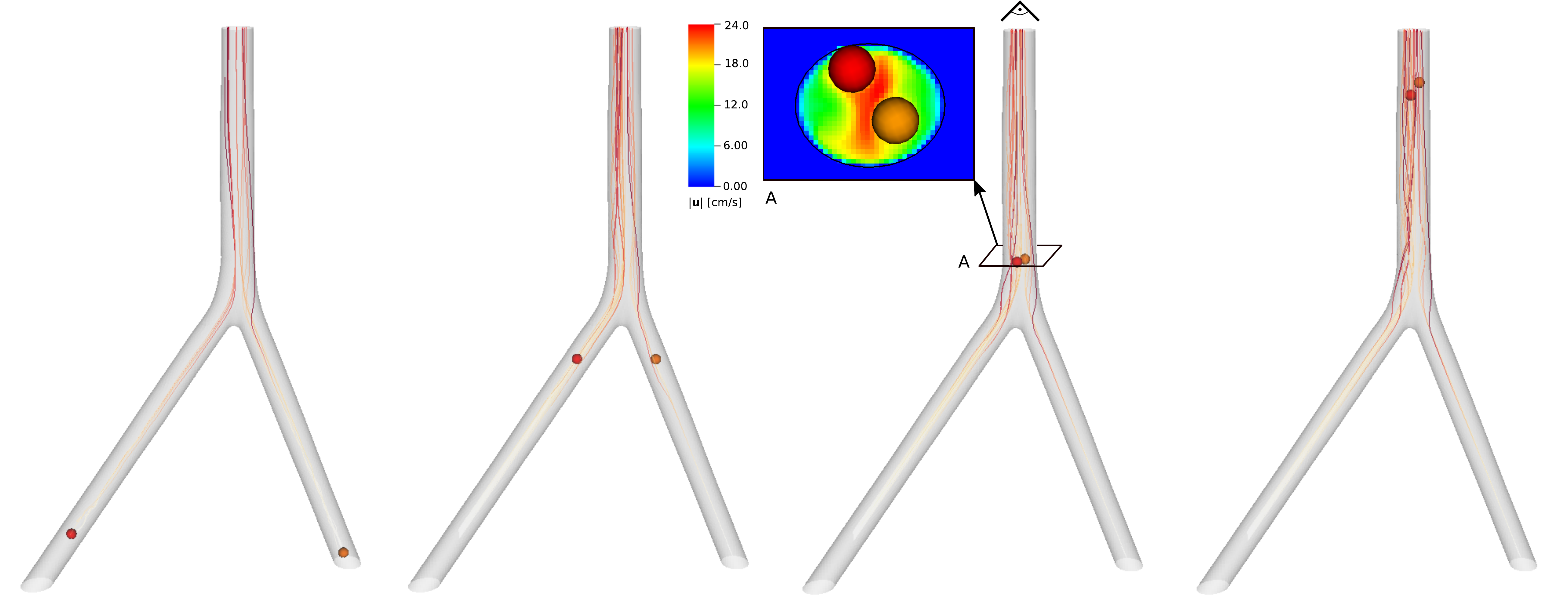}
\caption{Transport of spherical clots inside IVC under exercise condition. The snapshots are given from left to right at times 
$t=5.0$~s, $t=6.75$~s, $t=7.50$~s, and $t=8.50$~s. A top view of velocity magnitudes are shown for cross section A.} 

\label{fig:ivc_result}
\end{figure}  
The infrarenal IVC downstream of the iliac bifurcation has average hydraulic diameter $D_\text{h} =2.8 $ cm. 
The density of the fluid is $\rhof=1.817 $~$\textrm{g}\cdot\textrm{cm}^{-3}$, and the viscosity is $\muf=5.487\times 10^{-2} $ g/cm s.
Steady fully-developed parabolic velocity boundary conditions are imposed at the upstream inlets of the iliac veins.
The two inlets are circular with diameter $D=2.44 $ cm, which then transitions to an elliptical shape a short distance downstream.
The surface of the IVC is described using bilinear quadrilateral elements with $\Mfac \approx 2$. 
The IVC is embedded in a rectangular computational domain of size $\Omega=L_x\times L_y \times L_z = 50 \, \textrm{cm} \times 25 \, \textrm{cm} \times 50 \, \textrm{cm}$.
The Eulerian domain is discretized using a three-level locally refined grid with a refinement ratio of four between the grid levels, resulting in a grid spacing of 
$h_{\text{coarsest}}=\frac{25}{16}  \,\text{cm}=1.5625$~cm on the coarsest level and $h_{\text{finest}}=\frac{25}{16\times 4^2} \,\text{cm} \approx 0.098$~cm on the finest grid level.
At the outlet, the normal traction and tangential velocities are set to zero.
Solid-wall boundary conditions are imposed along the remainder of $\p\Omega$.
Once steady state condition of the flow has been reached at $t=5.0$~s, two neutrally buoyant rigid spherical clots ($\rhos=\rhof$) of diameter $D_i=1$~cm are 
released from positions close to the two inlets.
The volumetric mesh of each sphere consists of hexahedral elements leading to a surface mesh composed of  bilinear quadrilateral elements with $\Mfac=2$. 
A fixed time step size of $\dt=0.1$~ms is used. The penalty spring constant associated with both the IVC and the clots is 
$\kappa=1.05 \times 10^{5}$~$\textrm{g}\cdot(\textrm{cm}\cdot \textrm{s})^{-2}$. Fig.~\ref{fig:ivc_result} shows the transport of 
rigid blood clots under the unconstrained rigid-body model. A few streamlines are also plotted to show the direction of the flow.
Despite the moderately high Reynolds number of the flow and the complex pathway, the clots migrate towards the outlet through the iliac veins while 
being confined within the IVC structure. After passing the confluence of the veins, the clots come in very close contact with each other
in a more complex fluid environment where we expect the emergence of a pair of counter-rotating vortices after the confluence 
as a result of the two streams from iliac veins merging together. The clots continue migrating towards the outlet with a higher speed 
where towards the end, the clot from the right ilict vein tends to get ahead. This can be attributed to the complex nature of the flow in the region
above the confluence.

%%%%%%%%%%%%%%%%%%%%%%%%%%%%%%%%%%%%%%%%%%%%%%%%%%%%%%%%%%%%%%%%%%%%%%%%%%%%%%%%%%%%%%%%%%%%%%%%%%%%%%%%%%%%%%%%%%%%%%%%%%%
\section{Discussion and Conclusions}
\label{sec:discussion}

This work has introduced a numerical approach to simulating 
fluid-structure interaction that we refer to as an immersed Lagrangian-Eulerian (ILE) method.
Results from applying this ILE method to benchmark problems of rigid body fluid-structure interaction with increasing difficulty were presented.
In addition, we also show representative results from applications of this methodology to two biomedical FSI models: 
a bileaflet mechanical heart valve 
under physiological conditions in a model of the aortic test section of a commercial pulse duplicator, and transport of rigid blood clots inside 
a patient-averaged model of the inferior vena cava under exercise conditions. 
Unlike existing partitioned methods for FSI, the ILE formulation uses an immersed approach to couple the fluid and structure subdomains and thereby reduces or even eliminates
the need for grid regeneration during dynamic simulations.
In this formulation, it is crucial to deploy a coupling approach that provides the forces from \textit{only} the \textit{exterior} physical fluid region, and in 
this work, we use a coupling scheme
based on the immersed interface method, which enables us to evaluate these external fluid tractions along the fluid-structure interface.
 At least in principle, however, the present method could be used with any coupling strategy that determines the net exterior fluid force acting on the 
fluid-solid interface. Notice that this excludes conventional immersed boundary formulations using regularized delta functions because 
such formulations provide the total force from \textit{both} the 
the exterior~(physical) and interior~(nonphysical) fluid regions.
We remark that the jump conditions associated with the singular interfacial force are projected onto continuous Lagrangian basis functions, as described previously \cite{kolahdouz2020immersed}.
This requires the solution of linear systems of equations. 
We emphasize, however, that these solves involve only surface degrees of freedom, and thus are an order of magnitude 
smaller than the volumetric fluid equations that also are solved in each time step. In addition, the linear systems involve 
discrete $L^2$ projection equations that can be solved optimally using simple iterative methods (e.g.~by a Krylov 
method preconditioned with a diagonal or lumped mass matrix).
Overall, the dominant cost of each time step is in solving the incompressible Navier-Stokes equations, and 
the remaining computations are relatively inexpensive. A strength of the current approach is that it enables the use of 
fast Cartesian grid solvers for the incompressible Navier-Stokes equations.

The FSI coupling strategy allows fluid and solid subproblems to be solved in 
 a partitioned manner as independent, nonconforming discretizations that are coupled through interface conditions. 
In our discretization approach, there exist two Lagrangian representations of the fluid-solid interface, 
including the boundary of the volumetric mesh used in solving the equations of rigid body dynamics and a surface mesh 
that moves with the local fluid velocity.
These two representations are constrained to move together by
a Lagrange multiplier surface force. Exactly imposing the constraint would require the solution of a saddle-point system that couples the 
Eulerian and Lagrangian variables, but we develop a practical numerical scheme
that avoids the complex numerical linear algebra associated with such systems
by relaxing this constraint using a penalty formulation.
In the penalty formulation, the surface mesh moves according to the local fluid velocity but exerts force locally to the
fluid as a weak imposition of the no-slip condition.
Discrepancies in the positions of the boundary representations can be controlled by increasing the penalty parameter. In the present work, 
the maximum relative displacement is always less than $0.1$ of the Cartesian grid spacing.

An attractive feature of the present ILE method is that, at least for the specific examples considered herein,
it enables the use of a simple Dirichlet-Neumann coupling scheme 
\cite{burman2009stabilization,banks2014analysis,bukac2016stability} 
without requiring either strong coupling 
or subiterations to maintain stability. In particular, the motion of the solid mesh is driven by the exterior fluid traction, and the motion 
of the solid mesh drives the motion of the fluid-structure interface representation used to impose the no-slip condition.
%~ A substantial benefit of this approach is that we can 
Although we do not theoretically determine whether the ILE formulation suffers from added mass-related instabilities, 
the computational tests reported herein suggest that it can stabily treat a broad range of density ratios, 
including structures that are less dense than the fluid, more dense than the fluid, and neutrally buoyant. This is demonstrated for multiple benchmark problems. 
%%%%%%%%%%%%%%%%%%%%%%%%%%%%%%%%%%From Introduction %%%%%%%%%%%%%%%%%%%%%%%%%%%%%%%%%%%%%%%%%%%%%%%%%%%%%%%%%%%%%%%
%~ In a recent study \cite{mohaghegh2017comparison} it was empirically shown that for problems with high added mass effect, 
%~ a so-called diffuse interface IB method that uses regularized delta functions is more robust
%~ than a sharp interface method based on ghost fluid method \cite{fedkiw1999non,fedkiw2002coupling,tseng2003ghost}.
%~ An important difference between diffuse-interface IB methods and sharp-interface methods like the ghost fluid method 
%~ is that the diffuse-interface method solves the equations of fluid dynamics in both the fluid and solid subdomains, whereas 
%~ the ghost fluid method only solves the fluid equations in the fluid subdomain. Like Peskin's IB method, the IIM also solves 
%~ the equations of fluid dynamics throughout the entire computational domain.
%%%%%%%%%%%%%%%%%%%%%%%%%%%%%%%%%%%%%%%%%%%%%%%%%%%%%%%%%%%%%%%%%%%%%%%%%%%%%%%%%%%%%%%%%%%%%%%%%%%%%%%%%%%%%%%
For instance, in a 2-DOF model of the oscillation of a cylinder under vortex induced vibration, we obtain stable results for mass ratios up to 40 times smaller than 
smallest value reported in recent prior work ~\cite{kim2018weak}. In the literature, added mass instabilities are typically ascribed to the treatment of 
fluid regions that become ``uncovered" by the structure 
as it moves through the computational domain.
These include ALE approaches \cite{forster2007artificial,causin2005added, le2001fluid} and other sharp interface methods \cite{borazjani2008curvilinear,vigmostad2010fluid}.
Specialized approaches can be needed with these methods to avoid pressure fluctuations for cases involving low and near-equal density ratios.
Although the total mass and momentum of the fluid are conserved, there are local changes in the fluid mass and momentum in the regions that are both ``covered'' and ``uncovered'' by the moving structure.
Evidently, these localized changes in the fluid mass and momentum can induce temporal discontinuities in the velocity.
In the IIM approach introduced by Li and Lai \cite{li2001immersed}, the fluid equations are solved on the entire computational domain, including 
region occupied by the structure, similar to diffuse-interface formulations like those used in Peskin's IB method.
%~ Specifically, the fluid momentum equation itself accounts for the transport of momentum from a location that transitions between 
%~ being ``uncovered'' at time $t$ to ``covered'' at time $t+\Delta t$; that momentum is transported to other parts of the 
%~ physical fluid domain and diffused by viscous damping, 
%~ including damping associated with the no-slip boundary conditions imposed at fluid-structure interfaces. 
%~ The local fluid momentum at locations that transition 
%~ from covered at time $t$ to uncovered at time $t+\Delta t$ is 
%~ likewise determined by transport of fluid momentum from other parts of the physical fluid domain. 
%~ This can be interpreted as a physics-based extrapolation procedure. 
%~ Looking at this from the perspective of the IIM,
%~ one can realize that 
Using the IIM terminology, the continuity of the velocity field across the fluid-structure interface implies $\llbracket \u(\x,t) \rrbracket = 0$. 
A direct consequence of this
condition for the immersed interface approach is that the jump in the material derivative 
of the velocity is also zero, i.e. $\llbracket\frac{{\mathrm D} \u}{\mathrm {Dt}}(\x,t)\rrbracket =0$.
This means that fluid trajectories do not cross the moving interface and 
fluid locations that are ``uncovered'' by the motion of the structure automatically possess 
velocities that are consistent with the equations of motion.
In contrast, for methods in which the 
fluid domain only exists on one side of the interface, the temporal jump associated with the material derivative 
of the velocity may be non-zero going from one time step to the next as a fixed Eulerian grid point may switch sides between consecutive time steps.
This temporal discontinuity in the fluid acceleration could be important as added mass instabilities 
arise when the inertial effect due to fluid forces are dominant. This could result
in a fluid pressure that is out of phase with the fluid acceleration.
 Although this has yet to be proved rigorously, we hypothesize that the apparent robustness of 
 the method to artificial added mass instabilities results from its
consistent treatment of the momentum of the fluid near the fluid-structure interface. 
Added mass-type instabilities might also be suppressed by our use of rigid-body structure models, which
 are driven by the net fluid force acting on the immersed structure rather than by pointwise forces. 
 Rigorous stability analyses and analytical investigations of the 
methodology similar to those recently developed for overset 
grid methods \cite{burman2009stabilization,guidoboni2009stable,banks2013stable,banks2018stable} may reveal stability criteria or clarify the absence of added mass-type instabilities in this formulation.
Finally, we note that although the formulation presented here assumes the use of 
a rigid body structural model, it is natural to extend this approach to immersed elastic bodies.

\section*{Acknowledgements}

We acknowledge research support through NIH Awards HL117063 and HL143336, NSF Awards DMS 1664645, CBET 175193, OAC 1450327, OAC 1652541, and OAC 1931516, 
and the U.S.~FDA Center for Devices and Radiological Health (CDRH) Critical Path program.
This research was supported in part by an appointment to the Research Participation Program at the U.S.~FDA administered by 
the Oak Ridge Institute for Science and Education through an interagency agreement between the U.S.~Department of Energy and FDA.
A.P.S.B also acknowledges support from NSF award OAC 1931368.
Computations were performed using facilities provided by University of North Carolina at Chapel Hill through the 
Research Computing division of UNC Information Technology Services and the high-performance computing clusters at the U.S.~FDA.
The findings and conclusions in this article have not been formally disseminated by the FDA and should not be construed to represent
 any agency determination or policy.
The mention of commercial products, their sources, or their use in connection with material reported herein is not to be construed as either an actual or implied endorsement of such products by the Department of Health and Human Services.
We also thank Kenneth Aycock, Saad Qadeer, Jianhua Qin, and Simone Rossi for their constructive comments to improve the manuscript.

\bibliography{paper-IIM-RBD}

\end{document}